\documentclass[a4paper,reqno, 10pt]{amsart}   	

\usepackage{amssymb}
\usepackage{amsmath}
\usepackage{amsthm}
\usepackage{amscd,mdwlist}
\usepackage[margin=1in]{geometry}

\usepackage[english]{babel}

\usepackage{xcolor}
\usepackage{cite}
\usepackage{float}
\usepackage{dsfont} 
\usepackage{enumitem}
\setlist[enumerate]{itemsep=0mm}
\usepackage{graphicx}

\usepackage{csquotes} 

\usepackage{tikz}

\usepackage[colorlinks=true, linkcolor=black, citecolor=black]{hyperref}
\usepackage[nameinlink]{cleveref}

\Crefname{subsection}{Subsection}{Subsections}

\addto\extrasenglish{

}

\theoremstyle{plain}
\newtheorem{theorem}{Theorem}[section]
\newtheorem{proposition}[theorem]{Proposition}
\newtheorem{corollary}[theorem]{Corollary}
\newtheorem{lemma}[theorem]{Lemma}
\newtheorem{definition}[theorem]{Definition}

\newtheorem*{theorem*}{Theorem}
\newtheorem*{corollary*}{Corollary}

\theoremstyle{remark}
\newtheorem{remark}[theorem]{Remark}

\AddToHook{env/proposition/begin}{\crefalias{theorem}{proposition}}
\AddToHook{env/lemma/begin}{\crefalias{theorem}{lemma}}
\AddToHook{env/corollary/begin}{\crefalias{theorem}{corollary}}
\AddToHook{env/remark/begin}{\crefalias{theorem}{remark}}
\AddToHook{env/exampes/begin}{\crefalias{theorem}{examples}}
\AddToHook{env/subsection/begin}{\crefalias{section}{subsection}}

\DeclareMathOperator{\supp}{supp}

\DeclareMathOperator{\divergence}{div}
\DeclareMathOperator{\im}{Im}

\DeclareMathOperator{\sgn}{sgn}

\DeclareMathOperator{\SG}{SG}

\DeclareMathOperator{\id}{Id}

\newcommand{\ntrh}[1]{\vec{n}\cdot{#1}}

\newcommand{\norm}[1]{\left\| #1 \right\|}

\def\Xint#1{\mathchoice
	{\XXint\displaystyle\textstyle{#1}}%
	{\XXint\textstyle\scriptstyle{#1}}%
	{\XXint\scriptstyle\scriptscriptstyle{#1}}%
	{\XXint\scriptscriptstyle\scriptscriptstyle{#1}}%
	\!\int}
\def\XXint#1#2#3{{\setbox0=\hbox{$#1{#2#3}{\int}$}
		\vcenter{\hbox{$#2#3$}}\kern-.5\wd0}}

\def\dashint{\Xint-}

\begin{document}
	\author{Waldemar Schefer}
	\address{Fakult\"at f\"ur Mathematik, Universit\"at Bielefeld, Postfach 100131, 33501 Bielefeld, Germany}
	\email{wschefer@math.uni-bielefeld.de}

\title{On the domains of first order differential operators on the Sierpi\'{n}ski gasket}

\date{\today}

\begin{abstract}
	We study the first order structure of $L^2$-differential one-forms on the Sierpi\'{n}ski gasket. We consider piecewise energy finite functions related to one-forms, normal parts, and show self-similarity properties of one-forms as in the case of energy finite functions on the Sierpi\'{n}ski gasket. We introduce first order differential operators, taking functions into functions, that can be understood as total derivatives with respect to some reference function or form. The main result is a pointwise representation result for certain elements in the domain of said first order differential operator by a ratio of normal derivatives. At last we provide three classes of examples in the domain of said first order differential operator with different types of discontinuities.
	\tableofcontents
\end{abstract}
\maketitle

\section{Introduction}
We study first order differential operators, taking functions into functions, on the Sierpi\'{n}ski gasket. These operators can be understood as total derivatives with respect to certain reference functions or forms. In particular, the operators extend the total derivatives introduced in \cite{hino10}. They are also used in the article \cite{hs24} within the context of continuity and transport-type equations on fractals. Our main interest here is to investigate the domain of said first order differential operators and prove pointwise properties for functions in their domain. This also relates to the idea of \enquote{limited continuity} for Radon-Nikodym derivatives studied in \cite{bhs14}. As a key tool we use a basis for divergence-free one forms which relate to the loop structure of the underlying space as in \cite{cgis13}.

First derivatives of univariate functions are a fundamental operation of classical analysis. In the analysis on fractals counterparts of second order differential operators are well established. On the Sierpi\'{n}ski gasket they were first introduced stochastically by constructing diffusion processes, \cite{gold87, kus87, bp88}, and later analytically constructing Dirichlet forms, \cite{kig89, kus89, kig93}. Both approaches extend to various other classes of fractals and come with reliable tools to study the associated second order differential operator, see for example \cite{aof, def}. Much less is known in the first order case. Early concepts were introduced in \cite{kus89, kig93harmonicmetric, str00, tep00}. For p.c.f.\ fractals carrying resistance forms of index one (see \cite{hino10, kus93}), another notion of first order derivatives $\frac{df}{dg}$ of energy finite functions $f$ with respect to an energy-dominant reference function $g$ were introduced in \cite[Theorem 5.4]{hino10}. Related concepts were used in \cite[Section 4]{bk19} and \cite{hm20}. In terms of differential one-forms (see \cite{cs03, cs09, irt12}), these operators relate to the differential $\partial f$ of $f$ by the formula
\begin{equation}\label{E:intro differential formula}
	\partial f = \frac{df}{dg} \partial g,
\end{equation}
which may be seen as a generalization of the familiar identity $df = \frac{df}{dx} dx$ for differentiable functions on the real line. 

On the Sierpi\'{n}ski gasket $\SG$, we use the vector calculus introduced in \cite{cs03} (and further developed in \cite{cs09, irt12}) based on the standard self-similar resistance form and the usual three point boundary $V_0$. The crucial difference of our differential operators to the derivation introduced by the vector calculus, which maps functions into forms, is the use of the Hodge star operator. In the case of an $n$-dimensional riemmanian manifold, the \emph{Hodge star operator} is an isomorphism between $k$-forms and $(n-k)$-forms giving the \emph{Poincar\'{e} duality}. In our case, an important observation is that the usual energy form on $\SG$ is \enquote{one-dimensional}. This was first made in \cite{kus89} and later extended to other classes of fractals in \cite{hino10, hino13upperestimates}. Under this observation we define an analogue of the Hodge star operator between functions (zero-forms) and one-forms as in \cite{bk19}. Thus, the combination of the derivation $\partial$, its adjoint $\partial^\ast$ that represents a negative divergence, and the Hodge star operator gives us two first order differential operators that both take functions into functions.

The effect of the Hodge star operator can be explained as follows. One can choose a one-form $\omega$ that represents the Hodge star via multiplication with functions. Then, the mentioned operators are analogues of
\begin{equation}
	u\mapsto \langle \nabla u, \omega\rangle \quad\text{and}\quad u\mapsto \divergence(u\omega).
\end{equation}
The main difference of those two operators, in particular if $\omega$ is divergence-free, are the satisfied boundary conditions in their respective domain. More precisely, in the case of a finite metric graph, the domain of the operator build from $\partial$ contains functions that are continuous in each vertex, while the domain of the operator build from $\partial^\ast$ contains functions that satisfy some weighted Kirchhoff condition in each vertex. This observation was first made in \cite[Section 5]{bk19} in an unweighted case and extended to more general cases in \cite[Lemma 6.9]{diss}.

In general, the domain of the \enquote{divergence-type} operator is larger than the domain of energy finite functions, and it is not easy to understand. In particular, it may contain elements that are discontinuous. In the article \cite{hs24} this operator with respect to a given divergence-free one-form $\omega$ is used as a key tool to obtain integration by parts formulas and well-posedness results for continuity and transport-type equations with velocity field $\omega$. Some first insight into the structure of the domain of said operator from a functional analysis point of view is provided there. In this article, we aim to improve our understanding of the domain both from a functional analysis point of view as well as in terms of related pointwise results.

Let us also note that the restriction to the Sierpi\'{n}ski gasket is made purely for simplicity. Similar results are true for any p.c.f.\ self-similar fractal.

\subsection{Main results}
We consider the usual notation for the Sierpi\'{n}ski gasket $K = \SG$, its three point boundary $V_0 = \{q_0, q_1, q_2\}$, and its iterated function system $\{F_i\}_{i\in\{0, 1, 2\}}$. Often $w$ is an arbitrary word $w= w_1w_2\cdots w_n$ of length $n = |w|$ and $K_w = F_w(K)$. On $\SG$ we consider the standard self-similar resistance form $(\mathcal E, \mathcal F)$. We denote the space of all square integrable differential one-forms as $\mathcal H$, the derivation by $\partial : \mathcal F \to \mathcal H$, taking the role of a gradient, and $\partial^\ast$ denotes the adjoint of the derivation, taking the role of a negative divergence. For precise definitions see \Cref{S:prelim}. We denote the Hodge star operator induced by the one-form $\omega\in\mathcal H$ by $\star_{\omega}$, its induced energy measure by $\nu_{\omega}$, and define the first order differential operator
\begin{equation*}
	\partial^\perp_{V_0} : D_{\nu_\omega}(\partial^\perp_{V_0}) \to L^2(K, \nu_{\omega}), \quad \partial^\perp_{V_0} f = -\partial^\ast_{V_0}(\star_{\omega}f) = -\partial^\ast_{V_0}(f\omega).
\end{equation*}
By $\partial^\ast_{V_0}$ we mean the dual operator of the derivation restricted to functions with zero boundary data on $V_0$. For the precise definitions, see \Cref{S:first order operator}. The focal point of this article is the following result.
\begin{theorem*}[See \Cref{T:pointwise representation}]
	Let $\omega\in \ker\partial^\ast$ be a divergence-free and minimal energy-dominant one-form and $f\in D_{\nu_\omega}(\partial^\perp_{V_0})$. If $f$ and $\omega$ have a finite loop expansion, then 
	\begin{equation}\label{E:pointwise representation intro}
		\lim_{m\to\infty} \dashint_{K_{wi^m}} f \ d\nu_{\omega} = \frac{\vec n \cdot (f\omega)(F_wq_i)}{\vec n\cdot \omega(F_wq_i)} \quad\text{for any junction point $F_wq_i\notin\ker\vec n \cdot\omega$}.
	\end{equation}
	Here, $\dashint_{K_w} f \ d_{\nu_\omega}$ represents the mean integral $\nu_{\omega}(K_w)^{-1} \int_{K_w} f \ d\nu_{\omega}$, $K_{wi^m}$ is the usual cell approximation of the junction point $F_wq_i$, and $\vec n \cdot \eta(q)$ denotes the normal part of the one-form $\eta$ at $q$.
\end{theorem*}
The absence of $V_0$ in the notation $\partial^\ast$ and $\partial^\perp$ means that one-forms or functions in their respective domain have zero boundary values. By the \enquote{loop expansion} we mean that there exists a countable basis for the subspace of divergence-free one-forms $\ker\partial^\ast$ which is characterized by the loop structure of the space. Similar for $\ker\partial^\perp$. Thus, we say an element has a finite loop expansion if its projection onto one of these kernels has a finite basis expansion. In the case when the loop expansion of $\omega$ is infinite, \eqref{E:pointwise representation intro} does not hold in general. This case is discussed in the more general version \Cref{T:general pointwise representation}.

The statement in \eqref{E:pointwise representation intro} can be interpreted as, such functions $f$ in $D_{\nu_{\omega}}(\partial^\perp_{V_0})$ are good enough to claim that the change of flow through a junction point $q$ after a modulation of $\omega$ by $f$ is described by \enquote{the value of $f$ at $q$}. For $f$ from the resistance form domain this may be clear, but for discontinuous elements $f\in D_{\nu_\omega}(\partial^\perp_{V_0})$ it is not. Moreover, if $u$ is from the resistance form domain, then $\partial^\perp_{V_0} u$ equals the unique element in $L^2(K, \nu_\omega)$, that we are going to denote by $\star_{\omega}^{-1}\partial u$ for reasons becoming apparent in the text, such that 
\begin{equation}\label{E:intro differential formula general}
	\partial u = (\star_{\omega}^{-1}\partial u) \omega.
\end{equation}
The identity \eqref{E:intro differential formula general} generalizes \eqref{E:intro differential formula} and is the analytic counterpart of a well known martingale representation theorem in which the right-hand side can be reinterpreted as a stochastic integral, see \cite[Section 6]{kus93} and \cite[Corollary 3.9 ($ii$)]{hino10}. Explicit information on the \enquote{integrand} $\star_{\omega}^{-1}\partial u$ is notoriously difficult to obtain. In the Gaussian case, one famous result is the Clark-Ocone formula from Malliavin calculus. If $u$ is in the domain of the infinitesimal generator of the resistance form, then $\star_{\omega}^{-1}\partial u$ is an element of $D_{\nu_\omega}(\partial^\perp_{V_0})$. Thus, \Cref{T:pointwise representation} is applicable for $f = \star_{\omega}^{-1}\partial u$ and \eqref{E:pointwise representation intro} can be viewed as a piece of information about the unknown \enquote{integrand} $\star_{\omega}^{-1}\partial u$, at least at junction points. This relation is explained in more detail in \Cref{R:SG pontwise representative}.

Another use of \Cref{T:pointwise representation} is discussed in \Cref{S:examples}. There, we use \eqref{E:pointwise representation intro} to study discontinuity properties of functions in $D_{\nu_\omega}(\partial^\perp_{V_0})$. We provide three different types of discontinuities for functions $f\in D_{\nu_\omega}(\partial^\perp_{V_0})$ that may happen at a junction point $q$. Intuitively, these can be summarized as follows:
\begin{enumerate}
	\item \Cref{S:discontinuity on cut}: $f$ can be discontinuous at $q$ if $q\in\ker \vec n\cdot \omega$.
	\item \Cref{S:discontinuity for loops}: $f$ is discontinuous at $q$ if the coefficients of the loops close to $q$ in the loop expansion shrink slower than $(3\slash 5)^n$.
	\item \Cref{S:discontinuity in general}: In several cases, $f$ is discontinuous at $q$ when approaching along vertical lines.
\end{enumerate}
The examples in \Cref{S:examples} are related to the \enquote{limited continuity} result in \cite[Theorem 5.5]{bhs14} in which the authors show that the Radon-Nikodym derivatives $d\nu_{f, g}\slash d\nu$, for energy finite functions $f, g$ and Kusuoka measure $\nu$, are continuous when restricted to any finite metric graph approximation of $\SG$. Note that the left hand-side of \eqref{E:pointwise representation intro} coincides with the Radon-Nikodym derivative $d\nu_{u, \omega}\slash d\nu_{\omega}$ for $f = \star_{\omega}^{-1}\partial u$. We think that an extension of \cite[Theorem 5.5]{bhs14} to more general one-forms and reference measures $\nu$ is possible, although we do not pursue this here. However, the examples in \Cref{S:examples} show some limits of such a generalization.

On the way to \Cref{T:pointwise representation}, we made various notable new observations. One very helpful property for the analysis on the Sierpi\'{n}ski gasket is the self-similarity property for the (standard) resistance form $(\mathcal E, \mathcal F)$, that says, $u\circ F_w\in \mathcal F$ for every $u\in\mathcal F$ and 
\begin{equation*}
	\mathcal E(u) =  \sum_{|w|=n} r^{-n}\mathcal E(u\circ F_w).
\end{equation*}
Here, $r=3\slash 5$ denotes the resistance scaling. We extend the self-similarity property to the space of all square integrable one-forms as follows:
\begin{theorem*}[See \Cref{T:omega Fw}]
	Let $\omega \in \mathcal H$ with Hodge decomposition $\omega = \partial\eta + \Theta$, for $\eta\in\mathcal F$ and divergence-free one-form $\Theta\in\ker\partial^\ast$. Moreover, let $\{\partial^{(|w|+1)} \psi_w \}_{w}$ denote the basis for $\ker\partial^\ast$ as in \Cref{T:onb} and write $\Theta = \sum_w \Theta_w \partial^{(|w|+1)}\partial\psi_w$. For any word $w$ we define
	\begin{equation*}
		\omega\circ F_w := \partial(\eta\circ F_w) + \sum_{w'} \Theta_{w'} \partial^{(|w'|-|w|+1)} (\psi_{w'}\circ F_w),
	\end{equation*}
	where $\partial^{(n)} = \partial$ for $n\leq 0$. Then, $\omega\circ F_w\in \mathcal H$ is a square integrable one-form with Hodge decomposition
	\begin{equation*}
		\begin{aligned}
			P_{\im\partial}(\omega\circ F_w) &= \partial \left(\eta\circ F_w + \sum_{w' : K_w \subsetneqq K_{w'}} \Theta_{w'}(\psi_{w'}\circ F_w)\right),\\
			P_{\ker\partial^\ast} (\omega\circ F_w) &= \sum_{\widetilde w} \Theta_{w\widetilde w} \partial^{(|\widetilde w|+1)}\psi_{\widetilde w},
		\end{aligned}	
	\end{equation*}
	where $P_{\im\partial}$ and $P_{\ker\partial^\ast}$ are the projections of $\mathcal H$ onto $\im\partial$ and $\ker\partial^\ast$, respectively. Moreover, the associated energy measures satisfy the scaling relation
	\begin{equation*}
		(\nu_{\omega, \omega'})_{F_w^{-1}} = r^{-|w|}\nu_{\omega\circ F_w, \omega'\circ F_w} \quad\text{for all $\omega, \omega'\in\mathcal H$}.
	\end{equation*}
	Here, $(\nu_{\omega, \omega'})_{F_w^{-1}}(A) = \nu_{\omega, \omega'}(F_wA)$ is the pushforward measure.
\end{theorem*}
The self-similarity property in \Cref{T:omega Fw} allows us to derive further self-similar identities for the divergence operator $\partial^\ast$, the normal parts $\vec n\cdot \omega$, and tangential parts $\vec t \cdot \omega$. Moreover, it simplifies several proofs as it allows us to localize and reduce the proof to the situation of the boundary vertex $q_0$ and the cells close to $q_0$.

For the normal parts we study a property related to the limited continuity, it is the approximation of the normal part at some junction point $F_wq_i$ when approaching along the sides of the cell $K_w$. More precisely, we show a general version of the following:
\begin{theorem*}[See \Cref{T:convergence of ntrh}]
	Let $\nu$ be a non-atomic Radon measure on $K$ with full support and $\omega\in D_{\nu}(\partial^\ast_{V_0})$ with finite loop expansion. Then
	\begin{equation}\label{E:convergence of ntrh statement intro}
		\lim_{m\to\infty} \ntrh \omega(F_{0^m}q_1)= \lim_{m\to\infty} \ntrh \omega(F_{0^m}q_2) = -\frac 12 \ntrh \omega(q_0).
	\end{equation}
\end{theorem*}
This is a generalization of \cite[Lemma 4.3]{ls14}, in which they consider the Hausdorff measure and energy finite functions. In \Cref{T:convergence of ntrh} we are in particular interested in the case when $\omega$ has an infinite loop expansion. In this case we see that \eqref{E:convergence of ntrh statement intro} stays true if the loop coefficients shrink faster than $\left(3\slash 5\right)^m$ close to $q_0$. With the help of the self-similarity property in \Cref{T:omega Fw}, equation \eqref{E:convergence of ntrh statement intro} extends to arbitrary junction points.

At last, we have a small result on the existence of the tangential part of one-forms.
\begin{theorem*}[See \Cref{T:tangential part}]
	Let $\omega \in \mathcal H$ with Hodge decomposition $\omega = \partial\eta + \Theta$, for $\eta\in\mathcal F$ and divergence-free one-form $\Theta\in\ker\partial^\ast$, and write $\Theta = \sum_{w} \Theta_w \partial^{(|w|+1)}\psi_w$. If the tangential derivative $\partial_T \eta (q_0)$ for $\eta$ at $q_0$ exists and
	\begin{equation}\label{E:condition tangential part intro}
		\sum_{m=0}^\infty |\Theta_{0^m}|5^{m} < \infty,
	\end{equation}
	then the tangential part $\vec t \cdot\omega(q_0)$ of $\omega$ at $q_0$ exists and we have
	\begin{equation*}
		\vec t \cdot\omega(q_0) =\partial_T(\eta\circ F_w)(q_0) + 10 \sum_{m=0}^{\infty}\Theta_{0^m} 5^m.
	\end{equation*}	
\end{theorem*}
Although tangential derivatives do not play a dominant role in our results, they appear in various calculations, see for example \Cref{R:pointwise bound tangential part} and \Cref{S:discontinuity in general}.

\subsection{Overview}
In \Cref{S:prelim} we recall basic objects for the analysis on the Sierpi\'{n}ski gasket and we review the the construction of the space of \emph{differential one-forms} $\mathcal H$, the \emph{derivation} $\partial$, \emph{coderivation} $\partial^\ast$, and discuss the \emph{energy-measures} $\nu_\omega$ and minimal energy-dominant energy measures. The backbone for the analysis in this article are the \emph{piecewise energy finite functions} $P_n\mathcal F$, which are introduced and studied in \Cref{S:peff}, and the loop basis $\{\partial^{(|w|+1)}\psi_w\}_w$ for $\ker\partial^\ast$, which is discussed in \Cref{S:basis}.

Then, we investigate self-similarity properties in \Cref{S:self similarity}, \emph{normal parts} $\vec n \cdot \omega$ in \Cref{S:normal parts}, and \emph{tangential parts} $\vec t\cdot \omega$ in \Cref{S:tangential parts}. In \Cref{S:first order operator} we introduce the \emph{Hodge star operator} $\star_{\omega}$ and the first order differential operators $\star_{\omega}\partial$ and $\partial^\perp_{V_0}$ and recover some results for piecewise energy finite functions. Then, \Cref{S:pointwise representative} contains the main result on the pointwise representation of functions in the domain of $\partial^\perp_{V_0}$ and \Cref{S:examples} discusses the examples for discontinuous elements.

\subsection{Acknowledgements}
I would like to thank Michael Hinz for valuable comments on the manuscript. Financial support by the German Research Foundation (GRK 2235 - 282638148) is gratefully acknowledged.

\section{Preliminaries: Basics on the Sierpi\'{n}ski gasket and differential one-forms}\label{S:prelim}

\subsection{Basic notation on the Sierpi\'{n}ski gasket}
Throughout this text we use the following notation similar to \cite{aof, def}. Let $K= \SG \subset \mathbb R^2$ be the standard Sierpi\'{n}ski gasket, defined as the attractor of the iterated function system
\begin{equation*}
	F_i : \mathbb R^2\to\mathbb R^2, \quad F_i(x) := \frac 12(x-q_i) + q_i
\end{equation*}
where $q_0 = (1\slash 2, \sqrt 3\slash 2)$, $q_1=(0, 0)$ and $q_2=(1, 0)$. For any word $w = w_1 w_2 \ldots w_n$, $w_i\in \{0, 1, 2\}$, we denote $F_w = F_{w_1} \circ \ldots \circ F_{w_n}$ and $K_w = F_w(K)$. As usual, for another word $\eta = \eta_1 \ldots \eta_m$, we have $w\eta = w_1\ldots w_n \eta_1 \ldots \eta_m$ and $F_{w\eta} = F_w \circ F_\eta$. The length $|w|$ of $w = w_1 w_2 \ldots w_n$ is defined as $|w|=n$. The empty word $w = \emptyset$ is of length 0 and $F_\emptyset = \id$ is the identity. We will always use $w$ as a notation for a word and if nothing else is written, $w$ will be an arbitrary word of arbitrary length including the empty word.

We use the standard graph approximation on $\SG$, that is, the vertex sets are 
\begin{equation*}
	V_0 = \{q_0, q_1, q_3\}, \quad V_m := \bigcup_i F_i (V_{m-1}), \quad V_\ast = \bigcup_{m=0}^\infty V_m.
\end{equation*}
For $p, q \in V_n$ we say $p\sim_n q$ if there exists a word $w$ of length $n$ such that $p, q\in K_w$. Then, the standard self-similar resistance form $(\mathcal E, \mathcal F)$, with resistance scaling $r = 3\slash 5$, is given by
\begin{equation*}
	\mathcal E(u, v) = \lim_{m\to\infty} \mathcal E_m(u, v) = \lim_{m\to\infty} r^{-m} \sum_{q\sim_m p} (u(q) - u(p))(v(q) - v(p)) \quad\text{for $u, v\in\mathcal F$}.
\end{equation*}
Here the sum contains every edge $q\sim_m p$, $q, p\in V_n$, only once and $\mathcal F$ is defined as the set of all continuous $u:K\to\mathbb R$ where the supremum of $\mathcal E_m(u, u)$ exists. For the precise construction see \cite{aof, def}. We denote $\mathcal E(u) = \mathcal E(u, u)$ and similar for other bilinear expressions. The self-similarity property means for every $u\in\mathcal F$ we have $u\circ F_w\in\mathcal F$ and
\begin{equation}\label{E:energy self similar}
	\mathcal E(u, v) = \sum_{w\in \mathcal P} r^{-|w|}\mathcal E(u\circ F_w, v\circ F_w),
\end{equation}
for any partition $\mathcal P$ of $K$, that is, $\mathcal P$ is a collection of words such that $K$ is the union of $K_w$ for $w\in\mathcal P$ and the intersection of different $K_w$ contains at most a single point. For any word $w$ and all $u, v : K \to\mathbb R$ such that $u\circ F_w, v\circ F_w \in \mathcal F$ we set
\begin{equation*}
	\mathcal E_{K_w}(u, v) := r^{-|w|} \mathcal E(u\circ F_w, v\circ F_w).
\end{equation*}

For each $u, v\in \mathcal F$ the energy measure $\nu_{u, v}$ exists and is uniquely determined by the property
\begin{equation}\label{E:energy measure}
	\nu_{u, v}(F_wK) = \mathcal E_{K_w}(u, v) = r^{-|w|} \mathcal E(u\circ F_w, v\circ F_w) = \lim_{m\to\infty} r^{-m}\sum_{\substack{q\sim_m p \\ q, p\in K_w}}(u(q) - u(p))(v(q) - v(p)).
\end{equation}
We denote $\nu_u = \nu_{u, u}$.

Moreover, there exists a natural metric $R$, called \emph{resistance metric}, assoaciated with $\mathcal E$ such that $(K, R)$ is a compact metric space and every $u\in\mathcal F$ is $1\slash 2$-Hölder continuous. 

Set $\mathcal F_{V_0} := \{u\in\mathcal F\mid u = 0 \text{ on $V_0$} \}$. From \cite{kig12} it is well known that for any finite Borel measure $\nu$ on $(K, R)$ with full support, $(\mathcal E, \mathcal F)$ is an $L^2(K, \nu)$ strongly local regular Dirichlet form. Thus, for the generator we set
\begin{equation*}
	D_{\nu}(\Delta_{V_0}) := \{u\in\mathcal F \mid \varphi \mapsto \mathcal E(u, \varphi) \text{ is a bounded functional on $(\mathcal F_{V_0}, \norm{\cdot}_{L^2(K, \nu)})$} \}
\end{equation*}
and $\Delta_{V_0}u \in L^2(K, \nu)$ is defined via the Riesz representation theorem, satisfying
\begin{equation*}
	\mathcal E(u, \varphi) = -(\Delta_{V_0}u, \varphi)_{L^2(K, \nu)} \quad\text{for every $\varphi\in\mathcal F_{V_0}$}.
\end{equation*}
In particular, a function $u\in\mathcal F$ is called harmonic if
\begin{equation*}
	\mathcal E(u, \varphi) = 0 \quad\text{for every $\varphi\in\mathcal F_{V_0}$}.
\end{equation*}
For harmonicity, no measure $\nu$ is needed. Moreover, there is a simple way of calculating harmonic functions $h$ using the $\frac 25 - \frac 15$-rule, see \cite[Equation (1.3.18)]{def}.

The normal derivative $\partial_nu$ of $u\in \mathcal F$ on $V_0$ is defined by
\begin{equation}\label{E:normal deriv on boundary}
	\partial_n u (q_i) := \lim_{n\to\infty} r^{-n}(2u(q_i) - u(F_{i^n}q_{i+1}) - u(F_{i^n}q_{i+2}))
\end{equation}
if the limit exists. Here, and in similar expressions, we consider the index $i$ to be modulo 3. For arbitrary junctions points $q = F_wq_i$ one defines
\begin{equation}\label{E:normal deriv on junction}
	\partial_n u(F_wq_i) := r^{-|w|}\partial_n (u\circ F_w)(q_i).
\end{equation}
It is important to note that the normal derivative $\partial_n u(q)$ in $q\in V_\ast\setminus V_0$ depends on the orientation, that is, the cell $K_w$ that one chooses to calculate above limit. If $u\in D_{\nu}(\Delta_{V_0})$ then we know that the normal derivative exists at every junction point $q\in V_\ast$ and for different cells $q = F_wq_i = F_{w'}q_j$ the normal derivatives satisfy the \emph{matching condition}
\begin{equation}\label{E:matching condition}
	\partial_n u (F_wq_i) + \partial_n u(F_{w'}q_j) = 0.
\end{equation}
In formulas containing different normal derivatives we occasionally indicate the orientation by an arrow, for example for $x=F_1q_0 = F_0q_1$ we may write $\partial_n^{\nearrow} h(x) = \partial_n h(F_1q_0)$.

If $C$ is a finite union of cells with boundary $\partial C$, then the following \emph{local Gauss-Green formula} holds for any $u\in D_{\nu}(\Delta_{V_0})$ and $v\in\mathcal F$
\begin{equation}\label{E:local gauss green}
	\mathcal E_C(u, v) = -\int_C (\Delta_{V_0} u) v \ d\nu + \sum_{\partial C} (\partial_n u) v,
\end{equation}
see \cite[Theorem 2.4.1]{def}. Here, $\mathcal E_C$ denotes the sum of $\mathcal E_{K_w}$ for cells $K_w$ in $C$.

\subsection{Differential one-forms}\label{S:one forms}
We present the construction as in \cite{hm20}, however the result is the same as in \cite{ebe99, cs03, hrt13}. Let $C_a(K\times K)$ denote the space of continuous antisymmetric real valued functions on $K\times K$. For $f\in C(K)$ we write $\bar f(x, y) := \frac 12(f(x) + f(y))$. Then, we define an action of $C(K)$ on $C_a(K\times K)$ by
\begin{equation}\label{E:action}
	(\bar f \omega)(x, y) := \bar f(x, y) \omega(x, y) \quad\text{for $x, y\in K$},
\end{equation}
for any $f\in C(K)$ and $\omega\in C_a(K\times K)$, so $C_a(K\times K)$ becomes a module. The universal derivation $d : \mathcal F\to C_a(K\times K)$ is given by
\begin{equation*}
	du(x, y) := u(x) - u(y) \quad\text{for $x, y\in K$}.
\end{equation*}
It is immediate from \eqref{E:energy measure} that
\begin{equation}
	\int_{K} fg \ d\nu_{u, v} = \lim_{n\to\infty} r^{-n}\sum_{p\sim_n q} \bar f(p, q) \bar g(p, q) d u(p, q) d v(p, q).
\end{equation}
On $\Omega_a^1(K)$, the submodule of $C_a(K\times K)$ containing all finite linear combinations $\sum_i \bar u_i dv_i$, for $u_i\in C(K)$, $v_i\in\mathcal F$, the bilinear extension of
\begin{equation}\label{E:H seminorm}
	\langle \bar f d u, \bar g d v\rangle_{\mathcal H} := \int_{K} fg \ d\nu_{u, v} \quad\text{for $f, g\in C(K)$, $u, v\in\mathcal F$},
\end{equation}
defines a symmetric nonnegative definite bilinear form (cf.\ \cite[Remark 2.1]{hrt13}).

Let $\norm{\cdot}_{\mathcal H} := \sqrt{\langle \cdot, \cdot\rangle_{\mathcal H}}$ be the Hilbert seminorm on $\Omega_a^1(K)$ induced by \eqref{E:H seminorm}. Then $\mathcal H$ denotes the completion of the quotient space $\Omega_a^1(K)\slash \ker\norm{\cdot}_{\mathcal H}$ with respect to $\norm{\cdot}_{\mathcal H}$, becoming a Hilbert space (cf.\ \cite{cs03, cs09, hm20, hm22, hrt13}). 

By approximation and \eqref{E:H seminorm}, the action \eqref{E:action} extends to an action of the bounded Borel functions $\mathcal B_b(K)$ on $\mathcal H$ such that
\begin{equation}\label{E:action on H}
	\norm{f \omega}_{\mathcal H} \leq \norm{f}_{\sup} \norm{\omega}_{\mathcal H} \quad\text{for $f\in \mathcal B_b(K)$, $\omega\in\mathcal H$}.
\end{equation}
For any $u\in\mathcal F$ we denote its $\mathcal H$-equivalence class of the universal derivation by $\partial u := [d u]_\mathcal H$. Then,
\begin{equation*}
	\partial : \mathcal F\to \mathcal H
\end{equation*}
defines a linear operator, called \emph{derivation}, satisfying $\norm{\partial u}_{\mathcal H}^2 = \mathcal E(u)$, for any $u\in\mathcal F$, and the Leibniz rule
\begin{equation}\label{E:Leibniz rule}
	\partial(uv) = v\partial u + u\partial v \quad\text{for $u, v\in\mathcal F$}.
\end{equation}
By construction, the space of finite linear combinations $\sum_i f_i \partial u_i$, for $f_i, u_i\in\mathcal F$, is dense in $\mathcal H$. 

Let $\mathcal F^\ast$ denote the topological dual of $\mathcal F$. Then we denote by $\partial^\ast$ the dual operator of $\partial$, given by
\begin{equation}\label{E:distributional coderivation}
	\partial^\ast : \mathcal H\to\mathcal F^\ast, \quad \partial^\ast \omega(u) := \langle \omega, \partial u\rangle_{\mathcal H}.
\end{equation}
We call $\partial^\ast$ the \emph{coderivation}. The kernel $\ker\partial^\ast$ of $\partial^\ast$ is a closed linear subspace of $\mathcal H$. Using that $1\in\mathcal F$ and the fact that $(\mathcal F\slash \mathbb R, \mathcal E)$ is a Hilbert space, the image $\im \partial$ of $\partial$ under $\mathcal F$ is a closed subspace of $\mathcal H$, see \cite[section 4]{hkt15}. Then, $\mathcal H$ satisfies the orthogonal decomposition
\begin{equation}\label{E:Hodge decomposition}
	\mathcal H = \im\partial \oplus \ker\partial^\ast,
\end{equation}
called \emph{Hodge decomposition}.

Throughout this text, we consider $V_0$ as the boundary, denote
\begin{equation*}
	\mathcal F_{V_0} = \{u\in\mathcal F\mid u = 0 \text{ on $V_0$}\},
\end{equation*}
and write $\partial_{V_0}$ for the linear operator $\partial$ restricted to $\mathcal F_{V_0}$. For a finite Borel measure $\nu$ on $K$ with full support, $\partial$ and $\partial_{V_0}$ can be seen as densely defined closed unbounded linear operators on $L^2(K, \nu)$. Then, let $(\partial^\ast_{V_0}, D_{\nu}(\partial^\ast_{V_0}))$ be the adjoint of $(\partial_{V_0}, \mathcal F_{V_0})$, i.e., $\partial^\ast_{V_0} : \mathcal H \to L^2(K, \nu)$ with domain
\begin{equation*}
	D_{\nu}(\partial^\ast_{V_0}) := \{\omega \in \mathcal H \mid \text{there is $f\in L^2(K, \nu)$ s.t.\ }\langle \omega, \partial_{V_0} u \rangle_{\mathcal H} = (f, u)_{L^2(K, \nu)} \text{ for every $u\in\mathcal F_{V_0}$} \}
\end{equation*}
and $\partial^\ast_{V_0} \omega := f$. Clearly, $\partial^\ast_{V_0}$ agrees with an $L^2(K, \nu)$-representative of the distributional coderivation \eqref{E:distributional coderivation}. For ease of notation, the $\nu$ dependance appears only in the domain $D_\nu(\partial^\ast_{V_0})$. Analogously, let $(\partial^\ast, D_{\nu}(\partial^\ast))$ be the adjoint of $(\partial, \mathcal F)$. It is easy to see that $\partial^\ast_{V_0}$ is an extension of $\partial^\ast$. For $u\in D_{\nu}(\Delta_{V_0})$ we have $\partial u\in D_{\nu}(\partial^\ast_{V_0})$ with $\partial^\ast_{V_0} \partial u = -\Delta_{V_0} u$.

For $f, g\in C(K)$ and $u, v\in \mathcal F$ let us define $\nu_{f\partial u, g\partial v} := fg d\nu_{u, v}$. This map extends, using the density of finite linear combinations of the form $\sum_i f_i\partial u_i$, for $f_i\in C(K)$ and $u_i\in\mathcal F$, to a positive definite symmetric bilinear map $(\omega, \eta)\mapsto \nu_{\omega, \eta}$ from $\mathcal H\times\mathcal H$ into the space of finite signed Radon measures, see \cite[Lemma 2.1]{hrt13}. This map is uniquely determined by the property
\begin{equation}\label{E:energy measure of cell}
	\int_{K} \varphi \ d\nu_{\omega, \eta} = \langle \varphi\omega, \eta\rangle_{\mathcal H} \quad\text{for every $\varphi\in \mathcal B_b(K)$ and $\omega, \eta\in\mathcal H$.}
\end{equation}
 Note that $\nu_{\partial u, \partial v} = \nu_{u, v}$, for $u, v\in\mathcal F$, and $\nu_{\omega} := \nu_{\omega, \omega}$ is nonnegative for every $\omega\in\mathcal H$.

As in \cite[Definition 2.1]{hino10}, a \emph{minimal energy-dominant measure} $\nu$ on $K$ is a $\sigma$-finite Borel measure on $K$ s.t.\ $\nu_u \ll \nu$ for all $u\in\mathcal F$ (energy-dominant), and for any other $\sigma$-finite energy-dominant Borel measure $\nu'$ on $K$ we have $\nu \ll \nu'$ (minimality). Moreover, we will call an element $\omega\in\mathcal H$ minimal energy-dominant, if its induced energy measure $\nu_{\omega}$ is minimal energy-dominant.

\section{Piecewise energy finite functions}\label{S:peff}
We define \emph{piecewise energy finite functions} similar to \emph{edgewise energy finite functions} on metric graphs:
\begin{definition}
	A \emph{piecewise energy finite function} $u$ of level $n\in\mathbb N_0$ consists of the family $\{u^{(w)}\}_{|w|=n}$ of functions $u^{(w)} : K_w \to \mathbb R$ such that $u^{(w)} \circ F_w \in \mathcal F$. We denote the set of all such $u$ by $P_n\mathcal F$. For each $u\in P_n\mathcal F$ we associate a real-valued function as follows:
	\begin{enumerate}
		\item For $u\in P_n\mathcal F$ we associate a real-valued function on $K\setminus V_n$, that we also denote by $u$, defined by
		\begin{equation}\label{E:local representation 1}
			u(x) := \sum_{|w|=n} u^{(w)}(x) \mathds 1_{K_w}(x) \quad\text{for any $x\in K\setminus V_n$.}
		\end{equation}
		We call \eqref{E:local representation 1} the \emph{local representation} and the functions $u^{(w)}$ the \emph{local representatives} of $u$.
		\item We extend the definition of $u$ in \eqref{E:local representation 1} to boundary points $V_0$ by the same formula. 
		\item We call $u\in P_n\mathcal F$ \emph{continuous} in $q\in V_n\setminus V_0$ if the local representation is continuous in $q$, by which we mean that
		\begin{equation*}
			u^{(w)}(q) = u^{(w')}(q) \quad\text{for all $|w|=|w'| = n$ such that $K_w \cap K_{w'} = \{q\}$.}
		\end{equation*}
		In this case, we extend the function $u$ from \eqref{E:local representation 1} to $q$ by
		\begin{equation*}
			u(q) := u^{(w)}(q) = u^{(w')}(q).
		\end{equation*}
	\end{enumerate}
	For $n\leq 0$ let us set $P_n\mathcal F = \mathcal F$. We extend the definition of $\partial$ to $P_n\mathcal F$ by
	\begin{equation}\label{E:partial for PEFF}
		\partial^{(n)}:P_n\mathcal F \to \mathcal H, \quad\partial^{(n)} u := \sum_{|w|=n} \mathds 1_{K_w}\partial u^{(w)},
	\end{equation}
	where $u^{(w)}$ is an arbitrary extension in $\mathcal F$ of the eponymous local representative.
\end{definition}
Clearly, $P_n\mathcal F \subset P_{n+1}\mathcal F$ for any $n\in \mathbb N_0$ in the sense that any $u\in P_n\mathcal F$ has obvious local representatives in $P_{n+1}\mathcal F$ with respect to the local representation \eqref{E:local representation 1}. Thus, for $u\in P_n\mathcal F$ and $m\geq n$, its local representatives $u^{(w)}$ in $P_{m}\mathcal F$ are uniquely defined. Moreover, $u$ is continuous in $V_{m} \setminus V_n$ due to its local representation being continuous on higher level junction points. In particular, the local representation can be seen as a function on $K$ that is continuous up to a finite set of points at which the function may be defined arbitrarily. This makes the local representation \eqref{E:local representation 1} a Borel measurable function for any arbitrary extension of \eqref{E:local representation 1} to all of $K$. Thus, in the presence of an atomless Borel measure $\nu$ on $K$ we may regard $P_n\mathcal F$ as a set of $\nu$-equivalence classes with respect to the local representation \eqref{E:local representation 1}. We often use this fact as energy measures of elements in $\mathcal H$ are atomless.

\begin{remark}
	Given a metric graph, that is, $V$ is the set of vertices, $E$ the set of edges, and $I_e$ is the interval representing $e\in E$. Then, the set of edgewise energy finite functions (w.r.t.\ the standard Dirichlet energy) is $\bigoplus_{e\in E} H^1(I_e)$, where $H^1(I)$ denotes the usual Sobolev space on the interval $I$. In our case,  the cells $K_w$ replace the edges $I_e$ and the junction points $V_n$ the vertices $V$. In both cases, the piecewise and edgewise energy finite functions are in general discontinuous along junction points and vertices.
\end{remark}

Next, let us discuss the definition of the derivation $\partial^{(n)}$ in \eqref{E:partial for PEFF}. For every local representative of $u\in P_n\mathcal F$ there always exists an extension in $\mathcal F$. More precisely, if $v : K_w \to \mathbb R$ such that $v\circ F_w \in\mathcal F$ then we can construct an extension of $v$ in $\mathcal F$ by choosing a harmonic function on $K\setminus K_w$ with respect to the values $v(q)$ for $q\in \partial K_w$, see \cite[Section 2.5]{def}. Thus, for any $u\in P_n\mathcal F$ we can choose local representatives in $\mathcal F$ and, since $\mathds 1_{K_w}\partial v \in\mathcal H$ for $v\in\mathcal F$, \eqref{E:partial for PEFF} can be defined. However, the chosen local representatives $u^{(w)}\in\mathcal F$ are not unique outside of $K_w$. The locality of $\mathcal E$ ensures that $\partial^{(n)}$ is well-defined, meaning it is independent of the chosen extensions of the local representatives: Given $u^{(w)}, v^{(w)}\in \mathcal F$ such that $u^{(w)} = v^{(w)}$ on $K_w$. Then, from \eqref{E:H seminorm} we have
\begin{equation}\label{E:local representation in H}
	\norm{\mathds 1_{K_w}\partial u^{(w)} - \mathds 1_{K_w}\partial v^{(w)}}_{\mathcal H}^2 = \norm{\mathds 1_{K_w}\partial(u^{(w)} - v^{(w)})}_{\mathcal H}^2 = \nu_{u^{(w)} - v^{(w)}}(K_w) = 0.
\end{equation}
Thus, \eqref{E:partial for PEFF} is well-defined.

The relation \eqref{E:energy measure} implies that
\begin{equation*}
	\norm{\partial^{(n)}u}_{\mathcal H} = \sqrt{\sum_{|w|=n}\mathcal E_{K_w}(u^{(w)}) } \quad\text{for $u\in P_n\mathcal F$}.
\end{equation*}
Using the Leibniz rule \eqref{E:Leibniz rule} for $\partial$, one verifies that $\partial^{(n)}$ satisfies the Leibniz rule
\begin{equation}\label{E:Leibniz for peff}
	\partial^{(n\vee m)}(uv) = v\partial^{(n)}u + u \partial^{(m)}v \quad\text{for $u\in P_n\mathcal F$, $v\in P_m\mathcal F$}.
\end{equation}
Here, $\vee$ denotes the maximum and $u$, $v$ on the right-hand side of \eqref{E:Leibniz for peff} are seen as Borel measurable extensions of their local representations \eqref{E:local representation 1}. Moreover, \eqref{E:energy measure of cell} implies for any $\varphi\in\mathcal B_b(K)$ and $u\in P_n\mathcal F, v\in P_m\mathcal F$ that
\begin{equation}\label{E:SG peff in H scalar prod}
	\langle \varphi\partial^{(n)}u, \partial^{(m)}v\rangle_{\mathcal H} = \sum_{|w|=n\vee m}\langle \varphi \mathds 1_{K_w}\partial u^{(w)}, \partial v^{(w)}\rangle_{\mathcal H}.
\end{equation}
Here $u^{(w)}$ and $v^{(w)}$ are the local representatives of $u$ and $v$ in $P_{n\vee m}\mathcal F$.

Similar to the identity $\partial^\ast_{V_0} \partial u = -\Delta_{V_0} u$ for $u\in D_{\nu}(\Delta_{V_0})$, we prove a version for piecewise energy finite functions. For this, we define an extension of $\Delta_{V_0}$ to piecewise energy finite functions.
\begin{definition}
	Let $\nu$ be a finite Borel measure on $K$ with full support. For $n\in \mathbb N_0$ we define
	\begin{equation*}
		P_nD_{\nu}(\Delta_{V_0}) := \{u\in P_n\mathcal F \mid \text{there are local representatives such that } u^{(w)} \in D_{\nu}(\Delta_{V_0}) \text{ for all $|w|=n$} \}.
	\end{equation*}
	For any $u\in P_nD_{\nu}(\Delta_{V_0})$ we define
	\begin{equation*}
		\Delta_{V_0}^{(n)} u := \sum_{|w|=n} \mathds 1_{K_w} \Delta_{V_0} u^{(w)},
	\end{equation*}
	where $u^{(w)}$ are chosen representatives in $D_{\nu}(\Delta_{V_0})$.
\end{definition}
Similar to \eqref{E:local representation in H}, while the local representatives $u^{(w)}$ are only unique on $K_w$, if $u\in P_nD_{\nu}(\Delta_{V_0})$ then any local representative $u^{(w)} \in \mathcal F$ can be replaced by an extension of $u^{(w)}\mathds 1_{K_w}$ that is in $D_{\nu}(\Delta_{V_0})$. This can be seen, for example, by a harmonic extension to $K\setminus K_w$ where the normal derivatives sum up to 0 at the junction points $\partial K_w$, which is sufficient, see for example \cite[Theorem 2.5.1]{def}. A locality argument as in \eqref{E:local representation in H} implies that $P_nD_{\nu}(\Delta_{V_0})$ and $\Delta_{V_0}^{(n)}$ are well-defined and independent of the chosen local representation.
\begin{remark}\label{R:gluing of peff}
	As pointed out before, functions in $P_n\mathcal F$ are discontinuous on $V_n$ in general. Additionally, functions in $P_nD_{\nu}(\Delta_{V_0})$ may be discontinuous on $V_n$ and the matching condition of normal derivatives  \eqref{E:matching condition} may not be true on $V_n$. From the \enquote{gluing} of energy finite functions in \cite[Section 2.5]{def} the following are true:
	\begin{enumerate}
		\item For $u\in P_n\mathcal F$ we have $u\in\mathcal F$ if and only if $u$ is continuous on $V_n$.
		\item Let $\nu$ be a finite Borel measure and $u\in P_nD_\nu(\Delta_{V_0})$. Then $u\in D_{\nu}(\Delta_{V_0})$ if and only if $u$ is continuous on $V_n$ and the matching condition
		\begin{equation*}
			\partial_n u^{(w)}(F_wq_i) + \partial_n u^{(w')}(F_{w'}q_j) = 0 
		\end{equation*}
		holds for every $q=F_wq_i = F_{w'}q_j \in V_n$.
	\end{enumerate}
\end{remark}
Since the normal derivative is itself local, see \eqref{E:normal deriv on boundary} and \eqref{E:normal deriv on junction}, it is well defined for any piecewise energy finite function in $P_nD_{\nu}(\Delta_{V_0})$. Thus, the local Gauss-Green formula \eqref{E:local gauss green} extends to functions in $P_n\mathcal F$ as follows.
\begin{proposition}\label{P:SG peff gauss green}
	Let $\nu$ be a finite Borel measure on $K$ with full support. For $u\in P_nD_{\nu}(\Delta_{V_0})$ and $v\in P_m\mathcal F$ we have
	\begin{equation*}
		\langle \partial^{(n)} u, \partial^{(m)} v\rangle_{ \mathcal H} = -\int_K (\Delta_{V_0}^{(n)} u ) v \ d\nu + \sum_{|w|=n\vee m} \sum_{\partial K_w} (\partial_n u^{(w)}) v^{(w)}.
	\end{equation*}
\end{proposition}
\begin{proof}
	It follows immediately from \eqref{E:SG peff in H scalar prod}, \eqref{E:H seminorm}, \eqref{E:energy measure}, and the local Gauss-Green formula \eqref{E:local gauss green}.
\end{proof}
Next, we extend the identity $\partial^\ast_{V_0} \partial u = -\Delta_{V_0} u$ for piecewise energy finite functions in $P_nD_\nu(\Delta_{V_0})$.
\begin{proposition}\label{P:partial u in partial ast}
	Let $\nu$ be a finite Borel measure with full support on $K$ and $u\in P_n\mathcal F$. Then, $\partial^{(n)}u \in D_{\nu}(\partial^\ast_{V_0})$ if and only if $u\in P_nD_{\nu}(\Delta_{V_0})$ and for all $q\in V_n\setminus V_0$ we have
	\begin{equation}\label{E:normal derivatives sum to 0}
		\sum_{|w|=n} \mathds 1_{K_w}(q)\partial_n u^{(w)}(q) = 0.
	\end{equation}
	Also, $\partial^{(n)}u \in D_{\nu}(\partial^\ast)$ if and only if $u\in P_nD_{\nu}(\Delta_{V_0})$ and \eqref{E:normal derivatives sum to 0} holds for every $q\in V_n$. In both cases, we have
	\begin{equation*}
		\partial^\ast_{V_0} \partial^{(n)} u = -\sum_{|w|=n} \mathds 1_{K_w}\Delta_{V_0} u^{(w)} = -\Delta_{V_0}^{(n)}u.
	\end{equation*}
\end{proposition}
\begin{remark}
	The condition \eqref{E:normal derivatives sum to 0}, which is essentially the matching condition \eqref{E:matching condition}, was already observed in \cite{irt12} for piecewise harmonic functions.
\end{remark}
\begin{proof}
	Assume $u\in P_nD_{\nu}(\Delta_{V_0})$ satisfies \eqref{E:normal derivatives sum to 0} for all $q\in V_n\setminus V_0$. Then, the Gauss-Green formula from \Cref{P:SG peff gauss green} implies
	\begin{equation*}
		\langle\partial^{(n)} u, \partial \varphi\rangle_{ \mathcal H} =   -(\Delta_{V_0}^{(n)} u, \varphi)_{L^2(K, \nu)} + \sum_{|w|=n}\sum_{\partial K_w} \partial_n u^{(w)} \varphi \quad\text{for every $\varphi\in \mathcal F_{V_0}$.}
	\end{equation*}
	The continuity of $\varphi$, $\varphi\vert_{V_0} = 0$, and the matching condition \eqref{E:normal derivatives sum to 0} imply $\sum_{|w|=n}\sum_{\partial K_w} \partial_n u^{(w)} \varphi = 0$. Thus, $\partial^{(n)} u \in D_{\nu}(\partial^\ast_{V_0})$ with $\partial^\ast_{V_0} \partial^{(n)} u = -\Delta_{V_0}^{(n)} u$.
	
	Next, assume $\partial^{(n)} u \in D_{\nu}(\partial^\ast_{V_0})$. We are going to show that $u$ is locally in the domain of $\Delta_{V_0}$ with respect to a pushforward of $\nu$. Fix some $w$. For any $\psi\in\mathcal F_{V_0}$ there is $\varphi\in \mathcal F$ such that $\varphi\circ F_w = \psi$ and $\varphi \circ F_{w'} = 0$ for every $|w'| = n$, $w'\neq w$. This implies $\varphi = 0$ on $K\setminus K_w$ and in particular $\varphi\in\mathcal F_{V_0}$. Thus,
	\begin{equation*}
		\begin{aligned}
			\int_K \left(\left(\partial^\ast_{V_0}\partial^{(n)}u \right)\circ F_w\right) \psi \ d\nu_{F_w^{-1}} &= \int_{K_w} \left(\partial^\ast_{V_0}\partial^{(n)}u \right) \varphi \ d\nu = (\partial^\ast_{V_0}\partial^{(n)}u, \varphi)_{L^2(K, \nu)} = \langle \partial^{(n)} u, \partial\varphi \rangle_{\mathcal H} \\
			&= \sum_{|w'|=n} \langle\mathds 1_{K_{w'}}\partial u^{(w')}, \partial\varphi\rangle_{\mathcal H} = \langle \mathds 1_{K_w} \partial u^{(w)}, \partial\varphi\rangle_{\mathcal H} \\
			&= \mathcal E_{K_w}(u^{(w)}, \varphi) = r^{-n}\mathcal E(u^{(w)}\circ F_w, \psi).
		\end{aligned} 	
	\end{equation*}
	Hence, $u^{(w)} \circ F_w \in D_{\nu_{F_w^{-1}}}(\Delta_{V_0})$ and $-r^{-n} \Delta_{V_0} \left(u^{(w)}\circ F_w\right) = \left(\partial^\ast_{V_0}\partial^{(n)}u \right)\circ F_w$. This is true for any $w$.
	
	For general $\varphi\in\mathcal F$, the global Gauss-Green formula applied to each $w$ then implies
	\begin{equation*}
		\begin{aligned}
			\left(\partial^\ast_{V_0}\partial^{(n)} u, \varphi\right)_{L^2(K, \nu)} &= \sum_{|w|=n} \left((\partial^\ast_{V_0}\partial^{(n)} u) \circ F_w, \varphi\circ F_w\right)_{L^2(K, \nu_{F_w^{-1}})} \\
			&=\sum_{|w|=n} \left(r^{-n}\mathcal E(u^{(w)}\circ F_w, \varphi\circ F_w) - r^{-n}\sum_{\partial K} \partial_n (u^{(w)}\circ F_w) \varphi \circ F_w\right) \\
			&=\sum_{|w|=n} \mathcal E_{K_w}(u^{(w)}, \varphi) - \sum_{|w|=n} \sum_{\partial K_w} \partial_n u^{(w)} \varphi \\
			&= \langle \partial^{(n)} u, \partial\varphi\rangle_{\mathcal H} - \sum_{|w|=n} \sum_{\partial K_w} \partial_n u^{(w)} \varphi.
		\end{aligned}	
	\end{equation*}
	Note, by definition $r^{-n}\partial_n(u^{(w)}\circ F_w) = \partial_n u^{(w)} \circ F_w$. In particular, the normal derivatives $\partial_n u^{(w)}$ on $\partial K_w$ exist due to $u^{(w)} \circ F_w \in D_{\nu_{F_w^{-1}}}(\Delta_{V_0})$. Now, using the fact that $\partial^{(n)} u \in D_{\nu}(\partial^\ast_{V_0})$ one deduces \eqref{E:normal derivatives sum to 0} for all $q\in V_n\setminus V_0$ by testing with $\varphi\in\mathcal F_{V_0}$ that is 1 in $q$ and $0$ on any other level-$n$ junction point.
	
	It remains to show the existence of local representatives in $D_{\nu}(\Delta_{V_0})$. We replace the local representative $u^{(w)}$ (which were chosen abritrary in $\mathcal F$) by a harmonic extension of $u^{(w)}\vert_{K_w}$ that satisfies the matching condition \eqref{E:matching condition}. We denote this extension also by $u^{(w)}$. Then $u^{(w)} \in D_{\nu}(\Delta_{V_0})$, see \Cref{R:gluing of peff}, with $\Delta_{V_0} u^{(w)} = \mathds 1_{K_w} \partial^\ast_{V_0}\partial^{(n)} u$  since for any $\varphi\in\mathcal F_{V_0}$ we have
	\begin{equation*}
		\mathcal E(u^{(w)}, \varphi) = \sum_{|w'|=n} r^{-n} \mathcal E(u^{(w)}\circ F_{w'}, \varphi\circ F_{w'}) = (\partial^\ast_{V_0} \partial^{(n)} u \circ F_w , \varphi\circ F_w)_{L^2(K, \nu_{F_w^{-1}})} +  \sum_{|w'|=n } \sum_{\partial K_{w'}} \partial_n u^{(w)} \varphi. 
	\end{equation*}
	Note, we used the harmonicity of $u^{(w)}$ outside of $K_w$, which implies $\Delta_{V_0} (u^{(w)}\circ F_{w'} ) = 0$ for every $w' \neq w$. Due to the continuity of $\varphi$ and $u^{(w)}$ satisfying the matching condition \eqref{E:normal derivatives sum to 0}, only the integral on the right-hand side remains, which equals $(\mathds 1_{K_w}\partial^\ast_{V_0}\partial u ,\varphi)_{L^2(K, \nu)}$. Thus, $u\in P_nD_{\nu}(\Delta_{V_0})$.
	
	The proof is analogous in the case $\partial^{(n)}u\in D_{\nu}(\partial^\ast)$.
\end{proof}
Note, from the proof one sees that the definition of $u\in P_nD_{\nu}(\Delta_{V_0})$ via extended local representatives is equivalent to  the (local) assumption $u^{(w)}\circ F_w \in D_{\nu_{F_w^{-1}}}(\Delta_{V_0})$ for every $|w|=n$. At last, let us briefly discuss piecewise harmonic functions and their relation to minimal energy-dominant measures.
\begin{definition}
	Let $P_n \mathbb H_{V_0}$ denote the set of all $h\in P_n\mathcal F$ such that $h^{(w)}\circ F_w$ is harmonic for every $|w|=n$, that is,
	\begin{equation*}
		\mathcal E(h^{(w)}\circ F_w, \varphi) = 0 \quad\text{for every $\varphi\in \mathcal F_{V_0}$.}
	\end{equation*}
\end{definition}
\begin{remark}
	We sometimes call elements in $h \in P_n\mathbb H_{V_0}$ \emph{piecewise $n$-harmonic}. In the literature, for example in \cite{aof}, a \emph{(piecewise) $n$-harmonic function} or \emph{level-$n$ piecewise harmonic function} is assumed to be continuous. Note that our definition does not need any continuity on junction points between $n$-cells. However, we are often in the case that $\partial^{(n)} h  \in D_{\nu} (\partial^\ast_{V_0})$. Thus, by \Cref{P:partial u in partial ast} the matching condition \eqref{E:normal derivatives sum to 0} is satisfied, which can be thought of a \enquote{continuity of the derivative of $h$}.
\end{remark}
\begin{proposition}[See Theorem 5.6, \cite{hino10}]\label{P:piecewise harmonic are med}
	Every element in $h\in P_n\mathbb H_{V_0}$ that is nonconstant on each $n$-cell, that is, $h^{(w)} \circ F_w$ is a nonconstant harmonic function for every $|w|=n$, is minimal energy-dominant in the sense that $\partial^{(n)} h \in \mathcal H$ is minimal energy-dominant.
\end{proposition}
\begin{proof}
	In \cite[Theorem 5.6]{hino10} the case $n = 0$ is proven. By localization to each $n$-cell, the general statement follows.
\end{proof}

\section{A basis for solenoidal one-forms}\label{S:basis}
The set $\ker\partial^\ast$ and its basis are studied in \cite{irt12} and \cite{cgis13} using different approaches. In \cite{irt12} the authors consider finitely ramified cell structures and construct a basis for $\ker\partial^\ast$ using piecewise $n$-harmonic functions. In particular, they showed the relation of $\ker\partial^\ast$ to the various levels of holes in the cell structure of $X$. The goal of \cite{cgis13} was to construct and study path integrals with respect to one-forms. Here, they focused on the case of the standard Sierpi\'{n}ski gasket. They construct a basis $\{dz_{w}\}_w$ for $\ker\partial^\ast$ whose relation to the holes of $\SG$, \emph{lacunas} in their terminology, is clear by definition. The element $dz_{w}$ is defined to be the $(|w|+1)$-exact one-form with minimal energy and normalized path integral along the lacuna $\ell_w$ of the cell $K_{w}$. We use the result of \cite{irt12} to calculate the basis. In \Cref{R:loops} we compare this basis to the one in \cite{cgis13}.
\begin{figure}[H]
	\centering
	\includegraphics[width=0.4\linewidth]{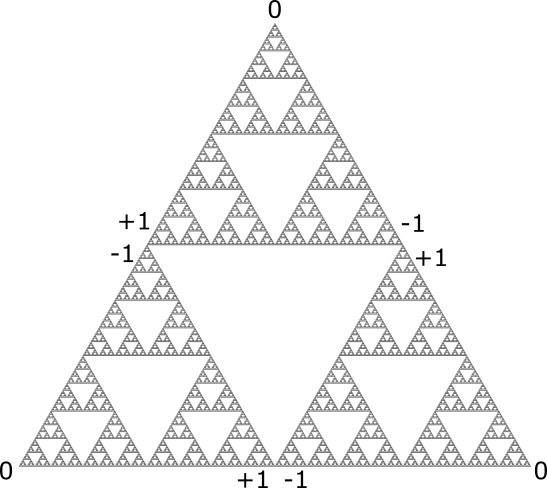}
	\caption{The distribution of boundary values on $V_1\cap K_i$ of the local representatives $\psi^{(i)}$. }
	\label{fig:psi}
\end{figure}
In \cite[Theorem 5.6.\ and 5.11.]{irt12} the authors prove for finitely ramified sets with a finitely ramified cell structure the existence of an orthogonal basis for $\ker\partial^\ast$ that consists of the \enquote{gradient} of piecewise nonconstant $n$-harmonic functions. Let us show that the following piecewise energy finite function $\psi$ generates this basis. Let us define $\psi \in P_1\mathcal F$ with local representative $\psi^{(i)}$ on $K_i$, being the continuous piecewise $1$-harmonic function given by the boundary values $\psi^{(i)}(q_i) = 0$, $\psi^{(i)}(F_iq_{i+1}) = +1$ and $\psi^{(i)}(F_iq_{i+2}) = -1$, where the index is seen mod 3. See also \Cref{fig:psi} for a representation of $\psi$ by the values of its local representatives $\psi^{(i)}$ on $V_1$. In particular, $\psi \in P_1 \mathbb H_{V_0}$.
\begin{theorem}\label{T:onb}
	For any $w$ define $\psi_w \in P_{(|w|+1)}\mathcal F$ by its local representatives
	\begin{equation}
		\psi_w^{(w)} := \psi \circ F_w^{-1} \quad\text{and}\quad \psi_w^{(w')} = 0 \quad\text{for every $w'\neq w$,}
	\end{equation}
	and if $w = \emptyset$ we set $\psi_{\emptyset} = \psi$. Then $\{\partial^{(|w|+1)}\psi_w \}_w = \bigcup_{n=0}^\infty\bigcup_{|w| = n} \{\partial^{(|w|+1)}\psi_w \}$ is an orthogonal basis for $\ker\partial^\ast$ with $\norm{\partial^{(|w|+1)}\psi_w}_{\mathcal H}^2 = 30\cdot (5\slash 3)^{|w|}$. Moreover, for every $\Phi \in \ker\partial^\ast$, there are $\Phi_w\in\mathbb R$ such that
	\begin{equation*}
		\sum_{n=0}^\infty \sum_{|w|=n} \Phi_w^2 \left(\frac 53 \right)^{|w|}< \infty
	\end{equation*}
	and
	\begin{equation*}
		\Phi = \sum_{n=0}^\infty\sum_{|w|=n} \Phi_w \partial^{(n+1)}\psi_w.
	\end{equation*}
	We often write $\Phi = \sum_w \Phi_w \partial^{(|w|+1)}\psi_w$. We say an element $\omega\in \mathcal H$ consists of \emph{finitely many loops} or has a \emph{finite loop expansion}, if the coefficients of its $\ker\partial^\ast$-part, with respect to the Hodge decomposition \eqref{E:Hodge decomposition}, are almost all zero with regard to the basis $\{\partial^{(|w|+1)}\psi_w \}_w$.
\end{theorem}

\begin{remark}\label{R:loops} \
	\begin{enumerate}
		\item The subset of all elements $\omega\in\mathcal H$ that consists of finitely many loops is dense in $\mathcal H$. This follows immediately from \cite[Theorem 5.6]{irt12}. We use this fact to state simpler versions of our results for such elements.
		\item In \cite{cgis13} the authors construct path integrals on $\SG$ for differential one-forms in $\mathcal H$. Crucial for their analysis are the, so called \emph{lacunas} $\ell_{w}$, that is the elementary path related to the first removed triangle in the cell $K_w$, see  \cite[Definition 2.1]{cgis13}. The lacunas are oriented clockwise. Associated with the lacuna $\ell_w$ is the $(|w|+1)$-exact one-form $dz_w$, that minimizes   among all $(|w|+1)$-exact forms with path integral 1 along $\ell_w$, see \cite[Definition 2.11]{cgis13}. In \cite[Theorem 2.27]{cgis13} the authors prove that the $dz_w$'s give an orthogonal basis for the space of harmonic forms, that is $\ker\partial^\ast$. Comparing $dz_w$ to the basis in \Cref{T:onb}, one can verify that $dz_w = - (1\slash 6)\partial^{(|w|+1)}\psi_w$. Here the minus sign is due to $\psi_w$ being oriented counterclockwise. Moreover, the local potentials associated with $dz_w$, see \cite[Proposition 2.10]{cgis13}, agree with the local representatives of $\psi_w$, more precisely with $- (1\slash 6)(\psi_w\circ F_i)$ on $K_{wi}$ for $i=0, 1, 2$.
	\end{enumerate}
\end{remark}

\begin{proof}[Proof of \Cref{T:onb}]
	Let $\widetilde{P_n\mathbb H_{V_0}}$ be the space of piecewise $n$-harmonic functions modulo additive constants on each $n$-cell that satisfy the matching condition \eqref{E:normal derivatives sum to 0} on $V_n$. Due to \cite[Theorem 5.6]{irt12}, we have
	\begin{equation*}
		\bigcup_n \partial^{(n)}(\widetilde{P_n\mathbb H_{V_0}})\subset \ker\partial^\ast \quad\text{dense}.
	\end{equation*}
	As in \cite[Theorem 5.11]{irt12}, we calculate the dimension of $\partial^{(n)}(\widetilde{P_n\mathbb H_{V_0}})$: Since $\ker\partial = \mathbb R$, we see that $\dim \partial^{(n)}(\widetilde{P_n\mathbb H_{V_0}}) = \dim \widetilde{P_n\mathbb H_{V_0}}$. The space of harmonic functions modulo constants on $\SG$ is 2-dimensional. Due to the self-similarity, this is also true on each cell and there are $3^n$ level-$n$ cells. Moreover, condition \eqref{E:normal derivatives sum to 0} for every junction point in $V_n$ reduces the dimension by $|V_n| -1$. Thus,
	\begin{equation*}
		\dim \widetilde{P_n\mathbb H_{V_0}} = 2\cdot 3^n - \left(|V_n| - 1\right) = \frac{3^n-1}{2}.
	\end{equation*}
	On the other hand, we have
	\begin{equation*}
		\left|\bigcup_{k=0}^{n-1}\bigcup_{|w|=k}\{\partial^{(k+1)}\psi_w \} \right| = \sum_{k=0}^{n-1} \left|\bigcup_{|w|=k}\{\partial^{(k+1)}\psi_w \}\right| = \sum_{k=0}^{n-1} 3^k = \frac{3^n -1}{2}.
	\end{equation*}
	It is easy to verify that $\psi \in \widetilde{P_1\mathbb H_{V_0}}$ and the $\psi_w$ are linearly independent for $|w|= 1$. Due to the self-similarity of $\SG$, this extends to all words. Thus the one-forms $\partial^{(|w|+1)}\psi_w$, for $|w|\leq n$, are a basis for $\partial^{(n+1)}(\widetilde{P_{n+1}\mathbb H_{V_0}})$. At last, the density implies that $\left\{\partial^{(|w|+1)}\psi_w\right\}_w$ is a basis for $\ker\partial^\ast$.
	
	For the orthogonality, take arbitrary $w$ and $w'$ and assume $|w|\geq |w'|$. Clearly, if $K_w\cap K_{w'}$ is empty or a single vertex, then  $\langle\partial^{(|w|+1)}\psi_w, \partial^{(|w'|+1)}\psi_{w'}\rangle_{\mathcal H} = 0$ by locality. Now, assume $K_w \subsetneqq K_{w'}$. Then there is a word $\eta$, $|\eta| \geq 1$, such that $w = w' \eta$. We calculate
	\begin{align*}
		\langle \partial^{(|w|+1)}\psi_w, \partial^{(|w'|+1)}\psi_{w'}\rangle_{ \mathcal H} &= \sum_{i}\mathcal E_{K_{wi}}(\psi_w, \psi_{w'}) = r^{-(|w|+1)} \sum_{i} \mathcal E(\psi_w\circ F_{wi}, \psi_{w'}\circ F_{wi}) \\
		&= r^{-(|w|+1)} \sum_{i} \sum_{V_0}\partial_n(\psi\circ F_i) \cdot \psi\circ F_{\eta i} = r^{-|w|} \sum_i \sum_{V_0} \partial_n \psi \circ F_i \cdot \psi \circ F_{\eta i} \\
		&= r^{-|w|} \sum_{i} \sum_{\partial K_i} \partial_n \psi \cdot \psi\circ F_\eta.
	\end{align*}
	In the last term we sum over $V_1$, while $V_1\setminus V_0$ appears twice (from both sides). On $V_0$, we have $\partial_n\psi = 0$. On $V_1\setminus V_0$ we use the continuity of $\psi\circ F_{\eta}$ (as $|\eta| \geq 1$) and the matching condition \eqref{E:normal derivatives sum to 0} for $\psi$. Thus, $\partial \psi_{w}$ and $\partial \psi_{w'}$ are orthogonal. 
	
	Similarly, we calculate
	\begin{equation*}
		\norm{\partial\psi_w}_{\mathcal H}^2 = r^{-(|w|+1)} \sum_{i}\mathcal E(\psi\circ F_i) = 18r^{-(|w|+1)} = 30\cdot (5\slash 3)^{|w|}.
	\end{equation*}
\end{proof}

\section{Self-similarity}\label{S:self similarity}

Due to the Hodge decomposition, every $\omega\in\mathcal H$ can be represented by $\omega = \partial\eta + \Theta$, where $\eta\in\mathcal F$ and $\Theta\in\ker\partial^\ast$. By the self-similarity of $(\mathcal E, \mathcal F)$, it is clear that $\partial(\eta\circ F_w)$ is well-defined for any word $w$. For a piecewise energy finite function $u\in P_n\mathcal F$ the composition $u\circ F_w$, with its local representatives being characterized by
\begin{equation*}
	\left\{u^{(w')}\circ F_w : |w'| = n \text{ and } K_{w'} \subset K_w \right\},
\end{equation*}
is also a piecewise energy finite function of level $n-|w|$. Now that $\Theta$ can be represented by \Cref{T:onb} with a basis consisting of derivations of piecewise energy finite functions, we can define $\Theta\circ F_w$ in $\mathcal H$ using derivations of $\psi_{w'}\circ F_w$. The self-similarity of $(\mathcal E, \mathcal F)$ in \eqref{E:energy self similar} is then recovered on $\mathcal H$ in terms of the respective energy measures.
\begin{theorem}[Self-similarity of $\mathcal H$]\label{T:omega Fw}
	Given $\omega\in \mathcal H$ with $\omega = \partial\eta + \Theta$ being its Hodge decomposition \eqref{E:Hodge decomposition} and write $\Theta = \sum_w \Theta_w \partial^{(|w|+1)}\psi_w$ as in \Cref{T:onb}. For any $w$ we define
	\begin{equation*}
		\omega\circ F_w := \partial(\eta\circ F_w) + \sum_{w'} \Theta_{w'} \partial^{(|w'|-|w|+1)} (\psi_{w'}\circ F_w),
	\end{equation*}
	where $\partial^{(n)} = \partial$ for $n\leq 0$. Then, $\omega\circ F_w\in \mathcal H$ with Hodge decomposition
	\begin{equation}\label{E:w circ Fw Hodge}
		\begin{aligned}
			P_{\im\partial}(\omega\circ F_w) &= \partial \left(\eta\circ F_w + \sum_{w' : K_w \subsetneqq K_{w'}} \Theta_{w'}(\psi_{w'}\circ F_w)\right),\\
			P_{\ker\partial^\ast} (\omega\circ F_w) &= \sum_{\widetilde w} \Theta_{w\widetilde w} \partial^{(|\widetilde w|+1)}\psi_{\widetilde w},
		\end{aligned}	
	\end{equation}
	where $P_{\im\partial}$ and $P_{\ker\partial^\ast}$ are the projections of $\mathcal H$ onto $\im\partial$ and $\ker\partial^\ast$, respectively. Moreover, the associated energy measures satisfy the scaling relation
	\begin{equation}\label{E:measure scaling}
		(\nu_{\omega, \omega'})_{F_w^{-1}} = r^{-|w|}\nu_{\omega\circ F_w, \omega'\circ F_w} \quad\text{for all $\omega, \omega'\in\mathcal H$}.
	\end{equation}
	Here, $(\nu_{\omega, \omega'})_{F_w^{-1}}(A) = \nu_{\omega, \omega'}(F_wA)$ is the pushforward measure.
\end{theorem}
\begin{proof}
	As $\supp \psi_{w'} = K_{w'}$, we can write
	\begin{equation}\label{E:omega Fw proof 1}
		\omega\circ F_w = \partial(\eta\circ F_w) + \sum_{w': K_{w}\subsetneqq K_{w'}} \Theta_{w'} \partial (\psi_{w'}\circ F_w) + \sum_{w': K_{w'}\subset K_{w}} \Theta_{w'} \partial^{(|w'|-|w|+1)} (\psi_{w'}\circ F_w).
	\end{equation}
	For $w'$ such that $K_{w}\subsetneqq K_{w'}$, there is $\widetilde w$, $1\leq |\widetilde w| \leq |w|-|w'|$, such that $w = w'\widetilde w$. In this case, $\psi_{w'}\circ F_w = \psi\circ F_{\widetilde w}$, and since $|\widetilde w|\geq 1$, we have $\psi\circ F_{\widetilde w}\in \mathcal F$. Thus, the sum of the first two terms in \eqref{E:omega Fw proof 1} is in $\im\partial$.
	
	In the contrary case, $w'$ such that $K_{w'}\subset K_{w}$, there is $\widetilde w$ such that $w' = w\widetilde w$, where $\widetilde w$ might be the empty word. Then $\psi_{w'}\circ F_w = \psi_{\widetilde w}$ and therefore
	\begin{equation*}
		\sum_{w': K_{w'}\subset K_{w}} \Theta_{w'} \partial^{(|w'|-|w|+1)} (\psi_{w'}\circ F_w) = \sum_{w': w' = w\widetilde w} \Theta_{w\widetilde w} \partial^{(|\widetilde{w}|+1)} (\psi_{w\widetilde w}\circ F_w) = \sum_{\widetilde w} \Theta_{w\widetilde w} \partial^{(|\widetilde{w}|+1)} \psi_{\widetilde w}.
	\end{equation*}
	Moreover, 
	\begin{equation*}
		\sum_{\widetilde w}\Theta_{w\widetilde w}^2 \left(\frac{5}{3} \right)^{|\widetilde w|} = \left(\frac{5}{3} \right)^{-|w|} \sum_{\widetilde w}\Theta_{w\widetilde w}^2 \left(\frac{5}{3} \right)^{|w\widetilde w|} \leq \left(\frac{5}{3} \right)^{-|w|} \sum_{w'}\Theta_{w'}^2 \left(\frac{5}{3} \right)^{|w'|} < \infty.
	\end{equation*}
	Thus, \eqref{E:omega Fw proof 1} implies that $\omega\circ F_w\in\mathcal H$ with \eqref{E:w circ Fw Hodge} as its Hodge decomposition \eqref{E:Hodge decomposition}.
	
	Now, let us prove the measure scaling relation \eqref{E:measure scaling}. For arbitrary $u, v\in \mathcal F$, the self-similarity \eqref{E:energy self similar} implies the energy measure scaling
	\begin{equation}\label{E:omega Fw proof 2}
		(\nu_{u, v})_{F_w^{-1}} = r^{-|w|} \nu_{u\circ F_w, v\circ F_w}.
	\end{equation}
	This follows from the following calculation: For any cell $K_{w'}$ we have
	\begin{align*}
		(\nu_{u, v})_{F_w^{-1}}(K_{w'}) &= \nu_{u, v}(F_{ww'}K) = r^{-|ww'|} \mathcal E(u\circ F_{ww'}, v\circ F_{ww'}) \\
		&= r^{-|w|} r^{-|w'|}\mathcal E((u\circ F_w)\circ F_{w'}, (v\circ F_w)\circ F_{w'} )\\
		&= r^{-|w|}\nu_{u\circ F_w, v\circ F_w}(K_{w'}).
	\end{align*}
	For $u, v\in P_n\mathcal F$, \eqref{E:SG peff in H scalar prod} gives 
	\begin{equation*}
		\nu_{\partial^{(n)} u, \partial^{(n)} v} = \sum_{|w|=n} \mathds 1_{K_w} d\nu_{u^{(w)}, v^{(w)}}.
	\end{equation*}
	By the definition of the pushforward measure, and since energy measures do not charge finite sets, \eqref{E:omega Fw proof 2} implies
	\begin{equation*}
		(\nu_{\partial^{(n)} u, \partial^{(n)} v})_{F_w^{-1}} = r^{-|w|} \nu_{\partial(u^{(w)}\circ F_w), \partial(v^{(w)}\circ F_w)}.
	\end{equation*}
	Since each of the functions $\eta$ and $\psi_w$ is a piecewise energy finite function, we conclude \eqref{E:measure scaling} by using \eqref{E:energy measure of cell} and the linearity in $\mathcal H$.
\end{proof}
As an immediate application, we have a self-similarity property for $\partial^\ast_{V_0}$. 
\begin{proposition}[Self-similarity of $\partial^\ast_{V_0}$]\label{P:partial ast omega Fw}
	Let $\nu$ be a finite Borel measure on $K$ with full support. If $\omega\in D_\nu(\partial^\ast_{V_0})$, then $\omega\circ F_w \in D_{\nu_{F_w^{-1}}}(\partial^\ast_{V_0})$ with
	\begin{equation*}
		\partial^\ast_{V_0}(\omega\circ F_w) = r^{|w|}(\partial^\ast_{V_0}\omega)\circ F_w \quad\text{$\nu_{F_w^{-1}}$-almost everywhere}.
	\end{equation*}
\end{proposition}
\begin{remark}
	Note, if $\omega = \partial u$, then $\omega\circ F_w = \partial(u\circ F_w)$ by the definition in \Cref{T:omega Fw}. Thus, \Cref{P:partial ast omega Fw} implies
	\begin{equation*}
		\Delta_{V_0}(u\circ F_w) = r^{|w|}(\Delta_{V_0} u) \circ F_w \quad\text{$\nu_{F_w^{-1}}$-almost everywhere.}
	\end{equation*}
	This coincides with the self-similarity of the Laplacian, for example as in \cite{def}, as follows: For this sake, let us change the notation of the Laplacian as in \cite{def}, that is, we denote $\Delta_{\nu}$ for the Laplacian with respect to the measure $\nu$. Then, \Cref{P:partial ast omega Fw} reads
	\begin{equation*}
		\Delta_{\nu_{F_w^{-1}}} (u\circ F_w) = r^{|w|}(\Delta_{\nu} u) \circ F_w \quad\text{$\nu_{F_w^{-1}}$-almost everywhere.}
	\end{equation*}
	If $\nu$ is a self-similar measure, lets denote it by $\mu$, with probabilities $\mu_i$, then it is easy to see that
	\begin{equation*}
		\Delta_{\mu_{F_w^{-1}}} v = \mu_w^{-1} \Delta_{\mu} v \quad\text{for any $v\in D(\Delta_{\mu})$.}
	\end{equation*}
	Thus, we recover
	\begin{equation*}
		\Delta_{\mu}(u\circ F_w) = r^{|w|}\mu_w (\Delta_{\mu} u) \circ F_w,
	\end{equation*}
	in an appropriate $\mu$-almost everywhere sense, which coincides with the self-similarity of the Laplacian as in \cite[Equation (2.1.8)]{def}.
\end{remark}
\begin{proof}
	Given $\varphi\in \mathcal F_{V_0}$. Let $\widetilde\varphi := \varphi \circ F_w^{-1}$ on $K_w$ and $0$ else. Note that $\widetilde\varphi\in\mathcal F$ as $\widetilde\varphi = 0$ on $\partial K_w$. As $\varphi = \widetilde\varphi \circ F_w$, the energy measure scaling \eqref{E:measure scaling} in \Cref{T:omega Fw} implies
	\begin{equation*}
		r^{-|w|}\langle \omega\circ F_w, \partial\varphi\rangle_{\mathcal H} = r^{-|w|}\langle \omega\circ F_w, \partial\widetilde\varphi \circ F_w\rangle_{\mathcal H} = (\nu_{\omega, \partial\widetilde\varphi})_{F_w^{-1}}(K) = \langle \omega, \mathds 1_{K_w}\partial\widetilde\varphi\rangle_{\mathcal H}.
	\end{equation*}
	Since $\supp \widetilde\varphi = K_w$, $\mathds 1_{K_w} \partial\widetilde\varphi = \partial\widetilde\varphi$, and $\omega\in D_{\nu}(\partial^\ast_{V_0})$ we have
	\begin{equation*}
		r^{-|w|}\langle \omega\circ F_w, \partial\varphi\rangle_{\mathcal H}  = (\partial^\ast_{V_0}\omega, \widetilde\varphi)_{L^2(K_w, \nu)} = (\partial^\ast_{V_0}\omega\circ F_w, \varphi)_{L^2(K, \nu_{F_w^{-1}})}.
	\end{equation*}
	The claim follows.
\end{proof}

\section{Normal parts}\label{S:normal parts}
Let us first define the \emph{normal part} $n_{V_0}\omega$ on the boundary $V_0$ for a one-form $\omega$ as in \cite{hs24} and extend its definition to the normal part $\vec n \cdot \omega$, defined on any junction point, the same way one extends the normal derivative with the help of the self-similarity.
\begin{definition}
	Let $\nu$ be a finite Borel measure on $K$ with full support and $\omega\in D_{\nu}(\partial^\ast_{V_0})$.
	\begin{enumerate}
		\item The normal part $n_{V_0} \omega$ of $\omega$ on $V_0$ is defined as the distribution
		\begin{equation*}
			n_{V_0}\omega(\varphi) := \langle \omega, \partial\varphi\rangle_{\mathcal H} - \int_{K} (\partial^\ast_{V_0} \omega) \varphi \ d\nu \quad\text{for $\varphi\in\mathcal F$.}
		\end{equation*}
		\item The normal part $n_{V_0}\omega(q)$ of $\omega$ at $q\in V_0$ is defined as
		\begin{equation*}
			n_{V_0}\omega(q) := n_{V_0}\omega(\psi_q),
		\end{equation*}
		where $\psi_q \in \mathcal F$ is the unique harmonic function with boundary values $\psi_q(q) = 1$ and $\psi_q = 0$ on $V_0\setminus \{q\}$.
		\item The normal part $\vec n \cdot \omega (F_wq)$ of $\omega$ at $F_wq$ from $K_w$, for any word $w$ and $q\in V_0$, is defined as
		\begin{equation*}
			\vec n \cdot \omega(F_wq) := r^{-|w|}n_{V_0}(\omega\circ F_w)(q).
		\end{equation*}
	\end{enumerate}
\end{definition}
Note, $\vec n \cdot\omega$ is well defined because of \Cref{P:partial ast omega Fw}, which gives $\omega\circ F_w \in D_{\nu_{F_w^{-1}}}(\partial^\ast_{V_0})$. Moreover, it coincides with the normal derivative for energy finite functions:
\begin{equation*}
	\vec n \cdot (\partial u)(F_wq) = \partial_n u (F_wq) \quad\text{for every $u\in D_{\nu}(\Delta_{V_0})$, $w$, and $q\in V_0$.}
\end{equation*}
For $\omega\in\mathcal H$ let $\omega = \partial\eta + \Theta$ be its Hodge decomposition \eqref{E:Hodge decomposition} and write $\Theta = \sum_w \Theta_w \partial^{(|w|+1)}\psi_w$. Clearly,  $\omega\in D_{\nu}(\partial^\ast_{V_0})$ if and only if $\eta\in D_{\nu}(\Delta_{V_0})$. In this case, $n_{V_0}\omega = n_{V_0}(\partial\eta) = \partial_n \eta$ as $\Theta\in \ker\partial^\ast$. Thus, if $\omega\in D_{\nu}(\partial^\ast_{V_0})$, then \Cref{T:omega Fw} and \Cref{P:partial ast omega Fw} imply
\begin{equation}\label{E:normal part on H}
	\vec n \cdot \omega(F_wq) = \partial_n \eta(F_wq) + \sum_{w'}\Theta_{w'} \partial_n\psi_{w'}(F_wq) \quad\text{for any $w$ and $q\in V_0$}.
\end{equation}
Some simple observations on the normal part:
\begin{enumerate}
	\item The sum in \eqref{E:normal part on H} is finite, more precisely, we only sum over cells $K_{w'}$ that contain $F_wq$ and up to words of length $|w|$. This is because $\partial_n\psi_{w'} = 0$ on $\partial K_{w'}$ and outside of $K_{w'}$. 
	\item $\ntrh\omega(q)$ for $q\in V_n$ is unique up to a sign, or more precisely, $\ntrh\omega$ satisfies the matching condition
	\begin{equation}\label{E:SG matching condition normal part}
		\ntrh\omega(F_{w}q_i) + \ntrh \omega(F_{w'}q_j) = 0 \quad\text{for every $q = F_wq_i = F_{w'}q_j$}.
	\end{equation} 
	This follows since we only sum over $|w'|\leq |w|$ and \eqref{E:normal derivatives sum to 0}.
	\item On $V_0$ we have $\ntrh\omega = \partial_n\eta$.
\end{enumerate}
\begin{remark}
	By \eqref{E:normal part on H}, one can extend the definition of $\vec n \cdot \omega$ to less regular $\omega\in \mathcal H$, by only requiring that the normal derivative of $\eta\in\mathcal F$ exists, as in \cite[Definition 2.3.1]{def}.
\end{remark}
Next, \Cref{T:omega Fw} and \Cref{P:partial ast omega Fw} imply the following local Gauss-Green formula for $\partial^\ast_{V_0}$ and piecewise energy finite functions:
\begin{proposition}[Local Gauss-Green formula for $\partial^\ast_{V_0}$]\label{P:gauss green for div}
	Let $\nu$ be a finite Borel  measure on $K$ with full support. For $\omega\in D_{\nu}(\partial^\ast_{V_0})$ and $\varphi\in P_m\mathcal F$ we have
	\begin{equation*}
		\langle \omega, \partial^{(m)}\varphi\rangle_{ \mathcal H} = (\partial^\ast_{V_0}\omega, \varphi)_{L^2(K, \nu)}  + \sum_{|w|=m}\sum_{\partial K_w} (\ntrh\omega)\varphi^{(w)} .
	\end{equation*}
\end{proposition}
\begin{remark}
	The \enquote{locality} in this statement is given since we are free to test with $\mathds 1_{K_w}\partial\varphi$ for $\varphi\in\mathcal F$ and any $w$, as $\mathds 1_{K_w}\varphi \in P_{|w|}\mathcal F$.
\end{remark}
\begin{proof}
	Let $\omega=\partial\eta + \Theta \in \partial(\mathcal F)\oplus \ker\partial^\ast$ be its Hodge decomposition \eqref{E:Hodge decomposition}. If $\varphi\in\mathcal F$, then the standard Gauss-Green formula of $(\mathcal E, \mathcal F)$ implies
	\begin{equation}\label{E:guass green for div proof 1}
		\langle \omega, \partial\varphi\rangle_{ \mathcal H} = \mathcal E(\eta, \varphi) = (-\Delta_{V_0} \eta, \varphi)_{L^2(K, \nu)} + \sum_{V_0} (\partial_n\eta) \varphi.
	\end{equation}
	The statement follows as $\partial^\ast_{V_0} \omega = -\Delta_{V_0} \eta$ and $\vec n \cdot \omega = \vec n \cdot(\partial \eta) = \partial_n\eta$ on $V_0$.
	
	Now, for general $\varphi \in P_m\mathcal F$, the energy measure self-similarity \eqref{E:measure scaling} in \Cref{T:omega Fw} implies
	\begin{equation*}
		\langle \omega, \partial^{(m)}\varphi\rangle_{ \mathcal H}  = \sum_{|w|=m}\langle \omega, \mathds 1_{K_w}\partial \varphi^{(w)}\rangle_{ \mathcal H} = \sum_{|w|=m} r^{-|w|}\langle \omega\circ F_w, \partial( \varphi^{(w)} \circ F_w)\rangle_{ \mathcal H}.
	\end{equation*}
	Then, \Cref{P:partial ast omega Fw} allows us to use \eqref{E:guass green for div proof 1} for $\varphi^{(w)}\circ F_w\in \mathcal F$:
	\begin{equation*}
		\langle \omega, \partial^{(m)}\varphi\rangle_{ \mathcal H}  =\sum_{|w|=m}r^{-|w|} (\partial^\ast_{V_0}(\omega\circ F_w), \varphi^{(w)}\circ F_w)_{L^2(K, \nu_{F_w^{-1}})} + \sum_{|w|=m}\sum_{V_0} r^{-|w|}\left(\vec n\cdot(\omega\circ F_w)\right) \left(\varphi^{(w)}\circ F_w\right).
	\end{equation*}
	Using \Cref{P:partial ast omega Fw} to rewrite $\partial^\ast_{V_0}(\omega\circ F_w)$ into $\partial^\ast_{V_0} \omega$ and the definition of $\vec n \cdot \omega$ on junction points, we conclude the statement.
\end{proof}
The normal part $\ntrh \omega$ shares many similarities to the normal derivative on $\SG$. Next, we show that the well known approximation along the sides of the triangle in $\SG$ (see \cite[Lemma 4.3]{ls14}):
\begin{equation}\label{E:simple partialn u approximation}
	\lim_{n\to \infty} \partial_n u(F_{0^n}q_1) = \lim_{n\to \infty} \partial_n u(F_{0^n}q_2) = -\frac 12 \partial_n u(q_0) \quad\text{for $u\in D_{\nu}(\Delta_{V_0})$},
\end{equation}
is also satisfied by the normal part. However, in \cite[Lemma 4.3]{ls14}, the authors prove the statement in the case of the Hausdorff measure $\nu$. The proof for more general measures is the same but requires a little bit more effort. We prove the following generalization of \cite[Lemma 4.3]{ls14}:
\begin{proposition}\label{P:simple partialn u approximation}
	If $\nu$ is a non-atomic Radon measure on $K$ with full support, then \eqref{E:simple partialn u approximation} holds.
\end{proposition}
From the pointwise representation \eqref{E:normal part on H} of $\vec n \cdot\omega$ and \eqref{E:simple partialn u approximation}, the following is immediately implied by \Cref{P:simple partialn u approximation}:
\begin{corollary}\label{C:convergence of ntrh baby version}
	If $\nu$ is a non-atomic Radon measure on $K$ with full support, then
	\begin{equation*}
		\lim_{n\to\infty} \vec n\cdot \omega(F_{0^n}q_1) = \lim_{n\to \infty} \vec n\cdot\omega(F_{0^n}q_1) = -\frac 12 \vec n \cdot \omega(q_0)
	\end{equation*}
	for every $\omega\in D_{\nu}(\partial^\ast_{V_0})$ that consists of finitely many loops
\end{corollary}
In the case when $\omega$ has infinitely many loops, such a convergence does not hold in general. We study this in \Cref{T:convergence of ntrh}. Now, let us prove \Cref{P:simple partialn u approximation}. The strategy is the same as in the proof of \cite[Lemma 4.3]{ls14}. First, we estimate $u$, more precisely $u(q_0) - u$, at junction points close to $q_0$ in \Cref{lem:u(pn) estimate}. Although not necessary for the proof of \Cref{P:simple partialn u approximation}, we extend the rate of convergence in \Cref{lem:u(pn) estimate} to all $x$ in the cell $K_{0^n}$ in \Cref{lem:local behaviour}. This local behaviour generalizes \cite[Theorem 4.3]{bbst99}, see also \cite[Theorem 2.7.4]{def}. Then, in the proof of \Cref{P:simple partialn u approximation} we separate $u$ into its symmetric and skew-symmetric part with respect to $q_0$ and prove \eqref{E:simple partialn u approximation} for both individually. For the skew-symmetric part we need the estimates in \Cref{lem:u(pn) estimate} and \Cref{lem:local behaviour}.

Now, let us start with an estimate of $u(q_0) - u$ at junction points close to $q_0$.
\begin{lemma}\label{lem:u(pn) estimate}
	Let $\nu$ be a Radon measure on $K$ with full support. Then, for every $u\in D_{\nu}(\Delta_{V_0})$ with $u(q_0) = 0$ we have the following estimate for any $n\in\mathbb N$ and $i=1, 2$: 
	\begin{equation}\label{eq:u(pn) estimate}
		\begin{aligned}
			|u(F_{0^n}q_i)| &\leq r^n|\partial_n u(q_0)| +  \left(\frac 15\right)^n \left|u(q_1) - u(q_2)\right| +  \frac 12 r^n \nu(K_{0^n})^{1\slash 2}\norm{\Delta_{V_0}u}_{L^2(K_{0^n}, \nu)} \\
			&\qquad + 2 \sum_{k=0}^{n-1}\left(\frac{1}{5}\right)^{k} r^{n-k} \nu(K_{0^{n-1-k}})^{1\slash 2} \norm{\Delta_{V_0} u}_{L^2(K_{0^{n-1-k}}, \nu) }.
		\end{aligned}
	\end{equation}
	Note, for words we use the convention $w^0 = \emptyset$.
\end{lemma}
\begin{proof}
	We will first prove upper bounds on $|u(p_n) + u(p_n')|$ and $|u(p_n) - u(p_n')|$, where $p_n = F_{0^n}q_1$, $p_n' = F_{0^n}q_2$. Let $h^{(n)}_{q}$ be the unique continuous piecewise $n$-harmonic function that is 1 on $q$ and $0$ on all $n$-neighbours of $q$. Using Gauss-Green on $h^{(n)}_{q_0}$ we get
	\begin{equation}\label{eq green gauss on tent function corner}
		\int_{K_{0^n}} h^{(n)}_{q_0} \Delta_{V_0} u \ d\nu = r^{-n}(2u(q_0) - u(p_n) - u(p_n')) - \partial_n u(q_0).
	\end{equation}
	Since $u(q_0) = 0$, we have
	\begin{equation}\label{eq representation of u_s(p_n)}
		u(p_n) + u(p_n') = - r^n\partial_n u(q_0) - r^n \int_{K_{0^n}} h^{(n)}_{q_0} \Delta_{V_0} u \ d\nu.
	\end{equation}
	Using Hölders inequality, and the fact that $0\leq h^{(n)}_q \leq 1$, we conclude
	\begin{equation}\label{eq u+u estimate}
		|u(p_n) + u(p_n')| \leq r^{n}|\partial_n u(q_0)| + r^n \nu(K_{0^n})^{1\slash 2} \norm{\Delta_{V_0} u}_{L^2(K_{0^n}, \nu)}.
	\end{equation}
	Now we want to establish a similar bound on $u(p_n) - u(p_n')$. For this, we apply the classic local Gauss-Green formula \eqref{E:local gauss green} to $h_q^{(n+1)}$ for each level-$(n+1)$ junction point $q$ in the cell $K_{0^n}$:
	\begin{equation}\label{eq green gauss on tent function junction}
		-\int_K h^{(n+1)}_{q} \Delta_{V_0} u \ d\nu = r^{-{(n+1)}}\sum_{p\sim_{n+1} q}(u(q) - u(p)).
	\end{equation}
	Thus we get the following system of equations:
	\begin{equation*}
		\begin{pmatrix}
			4 & -1 & -1 \\
			-1 & 4 & -1 \\
			-1 & -1 & 4 
		\end{pmatrix}
		\cdot
		\begin{pmatrix}
			u(p_{n+1})  \\
			u(p_{n+1}')  \\
			u(z_{n+1})  
		\end{pmatrix}
		= \begin{pmatrix}
			u(q_0) + u(p_{n}) - r^{n+1} \int_K h^{(n+1)}_{p_{n+1}} \Delta_{V_0} u \ d\nu\\
			u(q_0) + u(p_{n}') - r^{n+1} \int_K h^{(n+1)}_{p_{n+1}'} \Delta_{V_0} u \ d\nu \\
			u(p_n) + u(p_{n}') - r^{n+1} \int_K h^{(n+1)}_{z_{n+1}} \Delta_{V_0} u \ d\nu
		\end{pmatrix}.
	\end{equation*}
	We can solve this system to get a representation formula for $u(p_{n+1})$ and $u(p_{n+1}')$. This implies the following formula for the difference:
	\begin{equation}\label{eq representation of u_a(p_n)}
		u(p_{n+1}) - u(p_{n+1}') = \frac 15(u(p_n) - u(p_n')) - 2 r^{n+1}\int_K (h^{(n+1)}_{p_{n+1}} - h^{(n+1)}_{p_{n+1}'} )\Delta_{V_0} u \ d\nu.
	\end{equation}
	Since $h^{(n+1)}_{p_{n+1}} - h^{(n+1)}_{p_{n+1}'}$ is supported on $K_{0^n}$, we can estimate the integral similar as above. Hence,
	\begin{equation}\label{eq u-u estimate}
		|u(p_{n}) - u(p_{n}')| \leq \left(\frac 15 \right)^n|u(q_1) - u(q_2)| + 4 \sum_{k=0}^{n-1}\left(\frac{1}{5}\right)^{k} r^{n-k} \nu(K_{0^{n-1-k}})^{1\slash 2} \norm{\Delta_{V_0} u}_{L^2(K_{0^{n-1-k}}, \nu) }.
	\end{equation}
	Recall that for words we use the convention $w^0  = \emptyset$. Thus, combining \eqref{eq u+u estimate} and \eqref{eq u-u estimate} we conclude the claimed estimate in \eqref{eq:u(pn) estimate}.
\end{proof}
With the help of \Cref{lem:u(pn) estimate} we can study the local behaviour of functions in $D_{\nu}(\Delta_{V_0})$ similar to \cite[Theorem 4.3]{bbst99}.
\begin{lemma}\label{lem:local behaviour}
	Let $\nu$ be a non-atomic Radon measure on $K$ with full support. Then there exists $C>0$ such that every $u\in D_{\nu}(\Delta_{V_0})$ satisfies 
	\begin{equation}\label{eq:local behaviour general}
		\sup_{x\in K_{0^n}}|u(q_0) - u(x)| \leq C\left(r^n |\partial_n u(q_0)| +  \left(\frac 15\right)^n \left|u(q_1) - u(q_2)\right| + o(r^n)\norm{\Delta_{V_0} u}_{L^2(K, \nu)} \right)
	\end{equation}
	for sufficiently large $n$. Here, $o(r^n)$ denotes a function such that $\lim_{n\to\infty} r^{-n}o(r^n) = 0$. If additionally there exists $N\in\mathbb N_{0}$, $c_0>0$ and $\alpha \geq 0$ such that
	\begin{equation}\label{eq:measure growth}
		\nu(K_{0^n}) \leq r^{\alpha n}c_0 \quad\text{for every word $n\geq N$,}
	\end{equation}
	then every $u\in D_{\nu}(\Delta_{V_0})$ satisfies the following behaviour at $q_0$:
	\begin{enumerate}
		\item If $\partial_n u(q_0) \neq 0$ then for any $n\geq N$
		\begin{equation}\label{eq:local behaviour under growth 1}
			\sup_{x\in K_{0^n}} |u(q_0) - u(x)| \lesssim r^n. 
		\end{equation}
		\item If $\partial_n u(q_0) = 0$ then for any $n\geq N$
		\begin{equation}\label{eq:local behaviour under growth 2}
			\sup_{x\in K_{0^n}} |u(q_0) - u(x)| \lesssim \begin{cases}
				r^{(1+\alpha\slash 2)n} &\text{if }\alpha <  \frac{\ln 9}{\ln 5 - \ln 3},\\
				r^{(1+\alpha\slash 2)n} n = \left(\frac 15\right)^n n &\text{if }\alpha =  \frac{\ln 9}{\ln 5 - \ln 3},\\
				\left(\frac 15\right)^n &\text{if }\alpha >  \frac{\ln 9}{\ln 5 - \ln 3}.
			\end{cases}
		\end{equation}
	\end{enumerate}
	The constants in \eqref{eq:local behaviour under growth 1} and \eqref{eq:local behaviour under growth 2} heavily depend on $u$ and \eqref{eq:measure growth}.
\end{lemma}
\begin{remark}
	Note that \eqref{eq:measure growth} allows $\alpha = 0$, which is satisfied by every Radon measure by local finiteness.
\end{remark}
\begin{proof}
	Without loss of generality let us assume that $N=0$. Else we can simply replace $K$ by $\widetilde K := K_{0^N}$ as then $\nu(\widetilde K_{0^n}) = \nu(K_{0^{N+n}})$. Moreover, we assume that $u(q_0) = 0$, else we take $u - u(q_0)$ instead.
	
	First, we show the claimed estimates on the boundary of the cell $K_{0^n}$ using the previous result in \Cref{lem:u(pn) estimate}. For the sum in \eqref{eq:u(pn) estimate} we use the estimate
	\begin{equation*}
		\begin{aligned}
			\sum_{k=0}^{n-1} \left(\frac{1}{5} \right)^k r^{n-k} \nu(K_{0^{n-k-1}})^{1\slash 2} &= r^n \left(\sum_{k=0}^{\lfloor \frac{n-1}{2} \rfloor} \left(\frac{1}{3} \right)^k\nu(K_{0^{n-k-1}})^{1\slash 2}   + \sum_{k=\lfloor \frac{n-1}{2} \rfloor +1}^{n-1 } \left(\frac{1}{3} \right)^k\nu(K_{0^{n-k-1}})^{1\slash 2}  \right) \\
			&\leq  r^n \left(\nu(K_{0^{\lfloor \frac{n-1}{2}\rfloor }})^{1\slash 2}\sum_{k=0}^{\lfloor \frac{n-1}{2} \rfloor} \left(\frac{1}{3} \right)^k   + \left(\frac{1}{3} \right)^{\lfloor \frac{n-1}{2} \rfloor} \sum_{k=\lfloor \frac{n-1}{2} \rfloor +1}^{n-1 } \nu(K_{0^{n-k-1}})^{1\slash 2}  \right)\\
			&\leq \frac 32 r^n \nu(K_{0^{\lfloor \frac{n-1}{2} \rfloor}})^{1\slash 2} + \nu(K)^{1\slash 2} r^n \left(\frac{1}{3} \right)^{\lfloor \frac{n-1}{2} \rfloor} n = o(r^n).
		\end{aligned}
	\end{equation*}
	Since $\nu$ is a non-atomic Radon measure, $\nu(K_{0^m})$ converges to 0 for $m\to\infty$, explaining the function $o(r^n)$. Under the additional growth condition \eqref{eq:measure growth} the sum in \eqref{eq:u(pn) estimate} is estimated by
	\begin{equation}\label{eq:local behaviour proof 1}
		\sum_{k=0}^{n-1} \left(\frac 15\right)^k r^{n-k} \nu(K_{0^{n-1-k}})^{1\slash 2} \leq c_0^{1\slash 2} r^{(1+\alpha\slash 2)n - 1\slash 2}\sum_{k=0}^{n-1}\left(\frac{1}{3r^{\alpha\slash 2}}\right)^k.
	\end{equation}
	Note that
	\begin{equation*}
		3r^{\alpha\slash 2} = 1 \quad\text{iff}\quad \alpha = -\frac{\ln 9}{\ln r} = \frac{\ln 9}{\ln 5 - \ln 3} \approx 4.3.
	\end{equation*}
	Thus, if $\alpha < \frac{\ln 9}{\ln 5 - \ln 3}$, then $3r^{\alpha\slash 2} < 1$ and the sum in \eqref{eq:local behaviour proof 1} is bounded by $(1-3r^{\alpha\slash 2})^{-1}$. If $\alpha = \frac{\ln 9}{\ln 5 - \ln 3}$, the sum is exactly $n$. If $\alpha > \frac{\ln 9}{\ln 5 - \ln 3}$ the sum is bounded by $(3r^{\alpha\slash 2})^{-n}((3r^{\alpha\slash 2})^{-1} -1)^{-1}$. Thus, together with \eqref{eq:u(pn) estimate} we have the following estimate
	\begin{equation*}
		|u(p_n)| \lesssim r^n|\partial_n u(q_0)| + \left(\frac{1}{5}\right)^{n}|u(q_1) - u(q_2)| + \norm{\Delta_{V_0}u}_{L^2(K, \nu)}\left(r^{(1+\alpha\slash 2)n} + \begin{cases}
			r^{(1+\alpha\slash 2)n} &\text{if }\alpha <  \frac{\ln 9}{\ln 5 - \ln 3}\\
			r^{(1+\alpha\slash 2)n}n &\text{if }\alpha =  \frac{\ln 9}{\ln 5 - \ln 3}\\
			r^{(1-\frac{\alpha \ln 3}{2\ln r})n} &\text{if }\alpha >  \frac{\ln 9}{\ln 5 - \ln 3}
		\end{cases}\right). 
	\end{equation*}
	Comparing the rates, we conclude the statements at the junction point $x = p_n$. The same follows for $x = p_n'$.
	
	For the general case we write $u\circ F_{0^n} = h + g$, where $h\in\mathcal F$ is harmonic in $K$ with boundary values $u\circ F_{0^n}$ on $V_0$ and $g \in \mathcal F_{V_0}$ is given by
	\begin{equation*}
		g(x) := \int_{K} G(x, y) (-\Delta_{V_0}(u\circ F_{0^n})) \ d\nu_{F_{0^n}^{-1}},
	\end{equation*}
	where $G$ denotes the Green function, see for example \cite[Section 2.6]{def}. Note, from \Cref{P:partial ast omega Fw} we know that $u\circ F_{0^n} \in D_{\nu_{F_{0^n}^{-1}}}(\Delta_{V_0})$ with $\Delta_{V_0}(u\circ F_{0^n}) = r^n (\Delta_{V_0} u \circ F_{0^n})$. Thus, the boundedness of $G$, by a constant $C_G>0$, implies the uniform bound
	\begin{equation*}
		|g(x)| = r^n\left|\int_{K_{0^n}} G(x, F_{0^n}^{-1}(y)) \Delta_{V_0} u  \ d\nu \right| \leq r^n \nu(K_{0^n})^{1\slash 2} \norm{\Delta_{V_0} u}_{L^2(K_{0^n}, \nu)} C_G.
	\end{equation*}
	Now, using the maximum principle on the harmonic part $h$ and above bound on $g$ we have
	\begin{equation*}
		\sup_{x\in K_{0^n}} |u(x)| \leq \sup_{\partial K_{0^n}} |u(x)| + r^{(1+\alpha\slash 2)n} \norm{\Delta_{V_0} u}_{L^2(K, \nu)}C_G.
	\end{equation*}
	In combination with the previously established estimates on the boundary points $\partial K_{0^n}$, we conclude the proof.
\end{proof}

\begin{proof}[Proof of \Cref{P:simple partialn u approximation}]
	First, let us assume that $\nu$ is symmetric with respect to the reflection trough the line between $q_0$ and $F_1q_2$. More precisely, if $R:K\to K$ denotes said reflection, then we assume that $\nu_R = \nu$.
	
	Let $u = u_s + u_a$, where $u_s$ is the symmetric and $u_a$ the skew-symmetric part of $u$ with respect to $q_0$, that is,  
	\begin{equation*}
		u_s(x) = \frac 12 \left(u(x) + u(Rx)\right) \quad\text{and}\quad u_a(x) = \frac 12 \left(u(x) - u(Rx)\right).
	\end{equation*}
	In particular, $u_s, u_a\in D_{\nu}(\Delta_{V_0})$. This is because $\varphi\circ R\in\mathcal F$ for any $\varphi\in \mathcal F$ and
	\begin{equation*}
		\mathcal E(u\circ R, \varphi) = \mathcal E(u, \varphi\circ R) = -\int_{K}  \Delta_{V_0} u \cdot \varphi\circ R \ d\nu = -\int_{K}  \Delta_{V_0} u \circ R \cdot \varphi\ d\nu_R = -\int_{K}  \Delta_{V_0} u \circ R \cdot\varphi\ d\nu
	\end{equation*}
	for every $\varphi\in\mathcal F_{V_0}$. This then implies $\Delta_{V_0} (u\circ R) = \Delta_{V_0} u \circ R$ and $\norm{\Delta_{V_0} (u\circ R)}_{L^2(K, \nu)} \leq \norm{\Delta_{V_0} u}_{L^2(K, \nu)}$.
	\begin{figure}[H]
		\centering
		\begin{tikzpicture}[scale=4]
			\draw (0, 0) -- (1, 0) -- (1/2,0.866) -- (0, 0);
			\draw (1/2, 0) -- (3/4, 0.433) -- (1/4, 0.433) -- (1/2, 0);
			\node[above] at (1/2, 0.866) {$q_0$};
			\node[left] at (0, 0) {$p_n$};
			\node[right] at (1, 0) {$p_n'$};
			\node[below] at (1/2, 0) {$z_{n+1}$};
			\node[right] at (3/4, 0.433) {$p_{n+1}'$};
			\node[left] at (1/4, 0.433) {$p_{n+1}$};
		\end{tikzpicture}
		\caption{The level-$n$ cell $K_{0^n}$ at $q_0$.}
		\label{fig:q0cell}
	\end{figure}
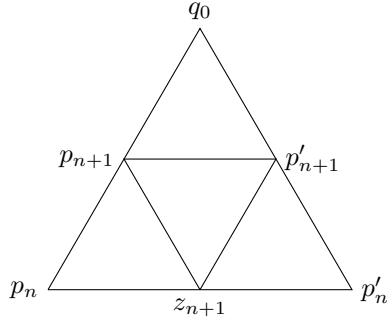
	Now let us prove \eqref{E:simple partialn u approximation} for $u_s$ and $u_a$.	Let us denote the vertices $V_{n+1}$ in the cell $K_{0^n}$ as in \Cref{fig:q0cell}. For the skew-symmetric part we have
	\begin{equation*}
		\partial_n u_a (q_0) = \lim_{n\to\infty} r^{-n}(2u_a(q_0) - u_a(p_n) - u_a(p_n')) = 0,
	\end{equation*}
	as $u_a(q_0) = 0$ and $u_a(p_n) = - u_a(p_n')$. Thus, we need to show that $\lim_{n\to\infty}\partial_n u_a(p_n) = 0$ since the skew-symmetry implies $\partial_n u_a(p_n) = - \partial_n u_a(p_n')$. The local Gauss-Green formula in $K_{0^n}$ with respect to $u_a$ and $h^{(n)}_{p_n}$, the continuous piecewise $n$-harmonic function that is $1$ on $p_n$ and zero on any other vertex in $V_n$, gives
	\begin{equation*}
		\partial_n u_a(p_n) = r^{-n}(2u_a(p_n) - u_a(q_0) - u_a(p_n')) + \int_{K_{0^n}} \Delta_{V_0} u_a h_{p_n}^{(n)} \ d\nu.
	\end{equation*}
	The integral term clearly converges to $0$ for $n\to\infty$. For the other term we use \eqref{eq:local behaviour general}, which gives
	\begin{equation*}
		r^{-n}(2u_a(p_n) - u_a(q_0) - u_a(p_n')) \leq Cr^{-n} o(r^n),
	\end{equation*}
	for some constant $C>0$ that depends on the boundary values of $u_a$ on $V_0$ and $\norm{\Delta_{V_0}u_a}_{L^2(K, \nu)}$. Hence, $\lim_{n\to\infty} \partial_n u_a(p_n) = 0$.
	
	Next, we consider the symmetric part $u_s$. Due to the symmetry we have $\partial_n u_s(p_n) = \partial_n u_s(p_n')$. Thus, 
	\begin{equation*}
		\int_{K_{0^n}} \Delta_{V_0} u_s \ d\nu = \partial_n u_s(q_0) + 2\partial_n u_s(p_n).
	\end{equation*}
	From here, we easily conclude that $2\lim_{n\to\infty}\partial_n u_s(p_n) = - \partial_n u_s(q_0)$. Hence, we have shown the statement in the case of symmetric measures.
	
	Now, let us assume that $\nu$ is no longer symmetric. We define the symmetrization of $\nu$ by $\nu_S$ simply by
	\begin{equation*}
		\nu_S := \frac 12(\nu + \nu_R).
	\end{equation*}
	Clearly, $\nu, \nu_R \ll \nu_S$. Moreover, we have
	\begin{equation*}
		\frac{d\nu}{d\nu_S} = 2 - \frac{d\nu_{R}}{d\nu_S}.
	\end{equation*}
	Thus, the range of the Radon-Nikodym derivatives $d\nu\slash d\nu_S$ and $d\nu_{R}\slash d\nu_S$ is contained in $[0, 2]$. It is easy to check that $u \in D_{\nu_S}(\Delta_{V_0})$ with 
	\begin{equation*}
		\Delta_{\nu_S} u = \frac{d\nu}{d\nu_S}\Delta_{\nu} u \quad\text{$\nu_S$-almost everywhere.}
	\end{equation*}
	Here, we use the notation $\Delta_{\mu}$ for the Laplacian $\Delta_{V_0}$ with respect to the measure $\mu$. Moreover, the measure $\nu_S$ stays a non-atomic Radon measure with full support as it is true for $\nu_R$. Hence, by our previous arguments, \eqref{E:simple partialn u approximation} holds for any function in $D_{\nu_S}(\Delta_{V_0})$ and therefore we conclude the proof of \Cref{P:simple partialn u approximation} as $D_{\nu}(\Delta_{V_0}) \subset D_{\nu_S}(\Delta_{V_0})$.
\end{proof}

At last, let us state the version of \Cref{C:convergence of ntrh baby version} for one-forms with infinite loop expansion. Intuitively, we are about to show that the same statement from \Cref{C:convergence of ntrh baby version} holds if the loops of $\omega \in D_{\nu}(\partial^\ast_{V_0})$ do not \enquote{pile up too much} close to the limit vertex. More precisely, the critical scaling of the loop coefficients $\Theta_w$, of the Hodge decomposition $\omega = \partial\eta + \Theta$, is $r^{|w|}$. Note, for $\Theta$ to be an element of $\mathcal H$, the critical scaling is $r^{|w|\slash 2}$, see \Cref{T:onb}.
\begin{theorem}\label{T:convergence of ntrh}
	Let $\nu$ be a non-atomic Radon measure on $K$ with full support and $\omega\in D_{\nu}(\partial^\ast_{V_0})$ with $\omega = \partial\eta + \Theta$ being its  Hodge decomposition \eqref{E:Hodge decomposition}. For any $q = F_wq_i$, consider the approximation $F_{wi^m}q_j$ for $j\neq i$.
	
	If there exists $N\in\mathbb N_0$, $a>0$ and $\theta \in [0, 3\slash 5)$ such that
	\begin{equation}\label{E:convergence of ntrh condition}
		|\Theta_{wi^m}| \leq a \theta^{|w|+m} \text{for all $m \in N$},
	\end{equation}
	then
	\begin{equation}\label{E:convergence of ntrh statement}
		\lim_{m\to\infty} \ntrh \omega(F_{wi^m}q_j) = -\frac 12 \ntrh \omega(F_wq_i).
	\end{equation}
	Note, $\ntrh \omega(F_{wi^m}q_j)$ and $\ntrh \omega(F_wq_i)$ are oppositely oriented. Moreover, \eqref{E:convergence of ntrh statement} generally fails if \eqref{E:convergence of ntrh condition} fails, that is, if $|\Theta_{wi^m}| \geq a (3\slash 5)^{|w|+m}$ for some $N\in\mathbb N_0$, $a>0$ and all $m\in N$.
\end{theorem}

\begin{proof}
	By definition
	\begin{equation*}
		\ntrh\omega(F_wq_i) = r^{-|w|}\ntrh{(\omega\circ F_w)}(q_i).
	\end{equation*}
	Thus, by \Cref{T:omega Fw} it is enough to proof \eqref{E:convergence of ntrh statement} for $q = q_0$ and approximating sequence $F_{0^m}q_j$, for $j = 1, 2$. Since \eqref{E:simple partialn u approximation} is a local property, by \Cref{P:simple partialn u approximation} it is true for any $u \in P_n\mathcal F$, in particular for the $\psi_w$. Due to \Cref{C:convergence of ntrh baby version} and $\vec n \cdot \omega(q_0) = \partial_n \eta(q_0)$, the general statement \eqref{E:convergence of ntrh statement} holds if and only if
	\begin{equation}\label{E:Convergenve of ntrh proof 1}
		\lim_{m\to\infty}\sum_{|w'|\geq n} \Theta_{w'} \partial_n\psi_{w'}(F_{0^m}q_j) = 0 \quad\text{for some $n\in\mathbb N_0$}.
	\end{equation}
	We  prove \eqref{E:Convergenve of ntrh proof 1}. Since $\supp \psi_{w'} = K_{w'}$, the sum in \eqref{E:Convergenve of ntrh proof 1} only considers $w' = 0^k$, for $k = n, n+1, \ldots, m$. One verifies that
	\begin{equation*}
		\partial_n \psi_{0^k}(F_{0^m}q_1) = - \partial_n \psi_{0^k}(F_{0^m}q_2) = 15\cdot \frac{5^k}{3^m}.
	\end{equation*}
	See also the proof of \Cref{P:discontinuity along sides using loops}. Thus,
	\begin{equation*}
		\sum_{|w'|\geq n} \Theta_{w'} \partial_n\psi_{w'}(F_{0^m}q_j) = \pm 15\cdot 3^{-m} \sum_{k=n}^m \Theta_{0^k}5^k.
	\end{equation*}
	Now, assumption \eqref{E:convergence of ntrh condition} for $n = N$, implies $|\Theta_{0^k}| \leq a\theta^{k}$. Then, if $\theta \neq 1\slash 5$, we have
	\begin{equation*}
		\left|3^{-m} \sum_{k=N}^m \Theta_{0^k}5^k\right| \leq a\cdot 3^{-m}\sum_{k=N}^m (5\theta)^k = \frac{3a}{5\theta -1}\cdot \left[\left(\frac{5}{3}\theta\right)^{m+1}- 3^{-(m+1)} (5\theta)^N\right].
	\end{equation*}
	For $\theta = 1\slash 5$ we get the bound $m3^{-m}$. Thus, \eqref{E:Convergenve of ntrh proof 1} is true for $n\geq N$, with $N$ coming from \eqref{E:convergence of ntrh condition}.
	
	To see that the convergence in \eqref{E:convergence of ntrh statement} in general fails if $\theta \geq 3\slash 5$, simply choose $\Theta_{w} = (3\slash 5)^{|w|}$. Then, from previous calculation we have
	\begin{equation*}
		\sum_{|w'|\geq n} \Theta_{w'} \partial_n\psi_{w'}(F_{0^m}q_j) = \pm 15\cdot 3^{-m} \sum_{k=n}^m \Theta_{0^k}5^k = \pm \frac{15}{2}\left( 3 - 3^{n-m}\right).
	\end{equation*}
	Thus, for $m\to\infty$ the term in \eqref{E:Convergenve of ntrh proof 1} converges $\pm 45\slash 2$. For $\Theta_w = \theta^{|w|}$ and $\theta > \frac 3 5$, the series in \eqref{E:Convergenve of ntrh proof 1} diverges to $\pm\infty$.
\end{proof}

\section{Tangential parts}\label{S:tangential parts}
The tangential derivative $\partial_Tu$ of $u\in\mathcal F$ on $V_0$ is defined as 
\begin{equation*}
	\partial_T u(q_i) = \lim_{n\to\infty} 5^n(u(F_{i^n}q_{i+1}) - u(F_{i^n}q_{i+2})),
\end{equation*}
if the limit exists. For $q = F_{w}q_i$ one then defines $\partial_T u(F_wq_i) := 5^{|w|}\partial_T(u\circ F_w)(q_i)$. Unlike for the normal derivative, $u\in D_{\nu}(\Delta_{V_0})$, for some  finite Borel measure with full support, is not sufficient for the existence of $\partial_T u$, see for example \cite[Theorem 2.7.8]{def} and the comment after. However, if $u$ is harmonic, then $\partial_T u$ exists and, similar to the normal derivative, can be calculated without the limit procedure. Thus, we can extend the definition of $\partial_T$ to any $u\in P_n\mathcal F$ and we know that $\partial_Tu$ exists for every piecewise harmonic $u$. 
\begin{definition}
	Let $\omega\in\mathcal H$ with $\omega = \partial\eta + \Theta$ being its Hodge decomposition \eqref{E:Hodge decomposition} and write $\Theta = \sum_{w} \Theta_w \partial^{(|w|+1)}\psi_w$. We define the \emph{tangential part} $\vec t \cdot \omega$ of $\omega$ in $F_wq_i$ as
	\begin{equation}\label{E:tangential part}
		\vec t \cdot \omega (F_wq_i) := \partial_T \eta(F_wq_i) + \sum_{w'} \Theta_{w'}\partial_T \psi_{w'} (F_w q_i),
	\end{equation}
	if $\partial_T \eta(F_wq_i)$ exists and the series $\sum_{w'} \Theta_{w'}\partial_T \psi_{w'} (F_w q_i)$ is absolutely convergent.
\end{definition}
We provide a sufficient condition that $\vec t\cdot \omega$ is well-defined.
\begin{theorem}\label{T:tangential part}
	Let $\omega\in \mathcal H$ with $\omega = \partial\eta + \Theta$ being its Hodge decomposition \eqref{E:Hodge decomposition} and $\Theta = \sum_{w} \Theta_w \partial^{(|w|+1)}\psi_w$. Let $q =  F_w q_i$. If $\partial_T \eta(F_wq_i)$ exists and
	\begin{equation}\label{E:condition tangential part}
		\sum_{m=0}^\infty |\Theta_{wi^m}|5^{m} < \infty,
	\end{equation}
	then $\vec t \cdot\omega(F_wq_i)$ is well-defined and we have
	\begin{equation*}
		\vec t \cdot\omega(F_w q_i) = 5^{|w|}\vec t\cdot(\omega\circ F_w)(q_i) = 5^{|w|}\left(\partial_T(\eta\circ F_w)(q_i) + 10 \sum_{m=0}^{\infty}\Theta_{wi^m} 5^m \right).
	\end{equation*}
\end{theorem}
\begin{proof}
	First, it is clear that $\vec t\cdot \omega(F_w q_i)$ exists if $\vec t\cdot(\omega\circ F_w)(q_i)$ exists. Then, the natural scaling for (piecewise) energy finite functions implies that $\vec t\cdot \omega(F_w q_i) = 5^{|w|}\vec t\cdot(\omega\circ F_w)(q_i)$. Moreover, the statement on $\eta$ is clear. Thus, it is enough to study the summability of $\sum_{w'} \Theta_{w'}\partial_T \psi_{w'} (F_w q_i)$.
	
	Let us prove the statement for $\eta = 0$ and $q = q_0$. Clearly, $\partial_T \psi_w(q_0) = 0$ unless $w = 0^m$, for some $m\in\mathbb N_0$. One verifies for any $k\in\mathbb N_0$ that
	\begin{equation*}
		\psi_{0^m}(F_{0^{m+1+k}}q_1)= (\psi\circ F_0)(F_{0^k}q_1) = 5^{-k}
	\end{equation*}
	and due to symmetry $\psi_{0^m}(F_{0^{m+1+k}}q_2) = -5^{-k}$. See also the proof of \Cref{P:discontinuity along sides using loops}. Thus,
	\begin{equation*}
		\partial_T\psi_{0^m}(q_0) = \lim_{k\to\infty} 5^{m+1+k}\left(\psi_{0^m}(F_{0^{m+1+k}}q_1) -\psi_{0^m}(F_{0^{m+1+k}}q_2) \right) = 10\cdot 5^m.
	\end{equation*}
	Hence, $\vec t \cdot\omega(q_0)$ is well-defined under condition \eqref{E:condition tangential part}.
\end{proof}

\section{Hodge star operator, Poincar\'{e} duality, and first order differential operators}\label{S:first order operator}
For any $\omega\in\mathcal H$ let us define the linear map
\begin{equation}\label{eq:hodge star}
	\star_\omega : L^2(K, \nu_\omega) \to \mathcal H, \quad \star_\omega f := f\omega.
\end{equation}
By definition of $\nu_\omega$, $\star_\omega$ is well-defined and isometric. We are now mostly interested in the case when $\star_\omega$ is an isomorphism. We have the known result:
\begin{corollary}[Proposition 4.1 in \cite{hs24}]
	For $\omega\in\mathcal H$ the linear map $\star_\omega$ is an isometric isomorphism if and only if $\omega$ is minimal energy-dominant.
\end{corollary}
Assume $\omega\in\mathcal H$ is minimal energy-dominant. In this case, we will call $\star_\omega$ \emph{Hodge star operator} as in \cite{bk19, hs24}, and denote the inverse $\star_\omega^{-1}$ again by $\star_\omega$. Similar to the  Hodge star operator on manifolds giving the Poincaré-Duality between harmonic $k$-forms and $(n-k)$-forms, $\star_\omega$ is an isometric isomorphism between 0-forms, i.e., scalar valued functions in $L^2(K, \nu_\omega)$, and one-forms in $\mathcal H$. Moreover, the bijectivity of $\star_{\omega}$ implies the following $L^2$-version of the Hodge decomposition \eqref{E:Hodge decomposition}:
\begin{equation}\label{eq L2 hodge}
	L^2(K, \nu_{\omega}) = \im \star_{\omega}\partial \oplus \ker\partial^\ast\star_{\omega}.
\end{equation}

The usual example for minimal energy-dominant $\omega$ is given by $\partial u$, where $u\in\mathcal F$ is minimal energy-dominant. From \cite[Proposition 2.7]{hino10}, the set of minimal energy-dominant $u\in\mathcal F$ is dense in $\mathcal F$. On $\SG$, we have the following nice class of minimal energy-dominant functions which is an immediate application of \cite[Theorem 5.6]{hino10}.
\begin{corollary}\label{C:harmonics are minimal energy dom}
	For $n\in\mathbb N_0$, every $u\in P_n\mathcal F$ that is piecewise $n$-harmonic and nonconstant on every cell is minimal energy-dominant, in the sense that $\partial^{(n)}u$ is minimal energy-dominant.
\end{corollary}
In particular, if additionally $\partial^{(n)} u$ satisfies the matching condition \eqref{E:normal derivatives sum to 0} on $V_n\setminus V_0$ then $\partial^{(n)} u \in \ker\partial^\ast_{V_0}$ by \Cref{P:partial u in partial ast}. This is a \enquote{simple} way to construct interesting examples of minimal energy-dominant $\omega \in \ker\partial^\ast_{V_0}$.

\begin{remark}
	In \cite[Proposition 2.7]{hino10} the author proves that the set of minimal energy-dominant $u\in\mathcal F$ is dense in $\mathcal F$ on general locally compact separable metric spaces with a regular symmetric Dirichlet form. The peculiarity of \cite[Theorem 5.6]{hino10} is that on Sierpi\'{n}ski gasket-type fractals every nonconstant harmonic function is minimal energy-dominant. This statement fails on many other p.c.f.\ fractals (e.g.\ Pentagasket, Vicsek set, and Hata's tree like set) as there are nonconstant harmonic functions that are constant on some cells. That means, in most cases only a subclass of the harmonic functions induce minimal energy-dominant measures. However, in the case that $X$ is a p.c.f.\ fractal one can simply remove the cells $C$ where the harmonic function is constant and study the differential operators on the restricted domain $X\setminus C$, see \cite[Remark 6.9]{hs24} and the examples mentioned therein.
\end{remark}

\begin{definition}
	Assume $\omega\in\mathcal H$ is minimal energy-dominant. Then we define $L^2(K, \nu_\omega)$-versions of $\partial$ and $\partial^\ast_{V_0}$ by
	\begin{equation*}
		\star_\omega \partial : \mathcal F \to L^2(K, \nu_\omega), \quad\star_\omega\partial u := \star_\omega(\partial u),
	\end{equation*}
	and
	\begin{equation*}
		\partial^\perp_{V_0} : D_{\nu_\omega}(\partial^\perp_{V_0}) \to  L^2(K, \nu_\omega), \quad \partial^\perp_{V_0}f := -\partial^\ast_{V_0}(\star_\omega f),
	\end{equation*}
	where
	\begin{equation*}
		D_{\nu_\omega}(\partial^\perp_{V_0}) := \{f \in L^2(K, \nu_\omega) \mid \star_\omega f\in D_{\nu_\omega}(\partial^\ast_{V_0})  \}.
	\end{equation*}
	Also, we consider
	\begin{equation*}
		D_{\nu_\omega}(\partial^\perp) := \{f\in L^2(K, \nu_w) \mid \star_\omega f\in D_{\nu_\omega}(\partial^\ast) \}, \quad \partial^\perp = (\partial^\perp_{V_0}, D_{\nu_\omega}(\partial^\perp)).
	\end{equation*}
\end{definition}

\begin{remark}
	The operators $\star_{\omega}\partial$ and $\partial^\perp_{V_0}$ were used in \cite{hs24} to study continuity and transport-type equations on fractals with respect to the velocity field $\omega$. There, it was important that $\omega\in \mathcal H$ is minimal energy-dominant and divergence-free, that is, $\omega\in\ker\partial^\ast$. Then, from \Cref{P:F in domain of div}, it follows that $\partial^\perp_{V_0}$ is an extension of $\star_\omega\partial$. It can then be shown that $\star_\omega\partial$ is skew-symmetric, which allows the use of \emph{boundary quadruples} to study maximally dissipative operators $\star_{\omega}\partial \subset A \subset \partial^\perp_{V_0}$.
\end{remark}

Let us discuss the relation of $\star_{\omega}\partial$ and $\partial^\perp_{V_0}$ to the total derivative in \cite{hino10}. In \cite[Theorem 5.4]{hino10} the author shows for index-one forms $(\mathcal E, \mathcal F)$ and minimal energy-dominant $g\in\mathcal F$ that for every $f\in \mathcal F$ there exists a total derivative $df\slash dg: K \to \mathbb R$ such that for $\nu_g$-almost every $x\in K$:
\begin{equation*}
	f(x) - f(y) = \frac{df}{dg}(x) (g(x) - g(y)) + R_x(y) \quad\text{for $y\in K$}.
\end{equation*}
The remainder term $R_x(y)$ converges to 0 of order $o(g(x) - g(y))$ for certain sequences $y\to x$. In particular, $df\slash dg$ induces a bounded linear operator 
\begin{equation*}
	\frac{d}{dg}: \mathcal F \to L^2(K, \nu_g),
\end{equation*}
with the property
\begin{equation*}
	\mathcal E(f, h) = \int_K \frac{df}{dg} \frac{dh}{dg}  \ d\nu_g \quad\text{for every $f, h\in\mathcal F$}.
\end{equation*}
From the proof of \cite[Theorem 5.4]{hino10} one also has for $f, h\in \mathcal F$
\begin{equation*}
	\frac{df}{dg} \frac{dh}{dg} = \frac{d\nu_{f, h}}{d\nu_g} \quad\text{$\nu_g$-almost everywhere.}
\end{equation*}
As clearly $dg\slash dg = 1$ $\nu_g$-almost everywhere, we have for any $\varphi\in C(K)$
\begin{equation*}
	\int_K \star_{\partial g} \partial f \cdot \varphi \ d\nu_g = \langle \partial f, \varphi\partial g\rangle_{ \mathcal H} = \int_K \varphi \ d\nu_{f, g} = \int_K \varphi \frac{d\nu_{f, g}}{d\nu_{g}} \ d\nu_{g} = \int_K \varphi \frac{df}{dg} \ d\nu_{g}.
\end{equation*}
Hence,
\begin{equation}\label{eq:hino derivative}
	\star_{\partial g}\partial f = \frac{df}{dg} \quad\text{$\nu_g$-almost everywhere for every $f\in\mathcal F$}
\end{equation}
and therefore
\begin{equation*}
	\partial^\perp_{V_0} = -\left(\frac{d}{dg} \right)^\ast \quad\text{for $\omega = \partial g$}.
\end{equation*}
This means, $\star_{\omega}\partial$ generalizes the total derivative $d\slash dg$ for differential one-forms $\omega\in\mathcal H$ instead of energy finite functions $g\in \mathcal F$ and $\partial^\perp_{V_0}$ is the negative adjoint operator of the total derivative.

Next, let us compare the operators $\star_{\omega}\partial$ and $\partial^\perp_{V_0}$ in the context of metric graphs to get an idea of their difference. Assume $X$ is a finite metric graph with vertex set $V$ and edge set $E$. Every edge shall be parametrized by $[0, 1]$. We define the energy form
\begin{equation*}
	\mathcal E(u, v) := \sum_{e\in E} \int_{0}^1 u_e' v_e' \ dx, \quad \mathcal F = \bigoplus_{e\in E} H^1(0, 1) \cap C(X).
\end{equation*}
Choose $\omega = \partial \eta$, for $\eta\in \mathcal F$ such that $\inf_{e\in E}\inf_{(0, 1)}\eta'_e >0$. Then $\omega$ will be minimal energy-dominant as $\nu_\omega= ((h_e')^2 dx)_{e\in E}$. Moreover, we will see that 
\begin{equation*}
	\star_{\omega}\partial u = \left(\frac{u_e'}{\eta_e'}\right)_{e\in E}.
\end{equation*}
We ignore the boundary here. Since $\partial^\perp$ is the negative adjoint of $\star_\omega\partial$, we see that
\begin{equation*}
	\partial^\perp f = \left( \frac{\left(f_e \eta_e'\right)'}{(\eta_e')^2} \right)_{e\in E}
\end{equation*}
with domain
\begin{equation*}
	D_{\nu_\omega}(\partial^\perp) = \{f \in L^2(X, \nu_\omega) : f_e \eta_e' \in H^1(0, 1) \text{ for all $e$ and } \sum_{e\in E} n_e(v)(f_e\eta_e')(v) = 0 \text{ for all $v$}\}.
\end{equation*}
Here $u_e(v)$, for $u = f\eta'$, is the appropriate value of $u_e$ when approaching $x=0$ or $1$, depending on the orientation corresponding to the vertex $v$, and $n_e(v)$ is the orientation of $e$ at $v$, i.e.,  $n_e(v) = +1$ if $v$ is identified with $1$ along $e$ and $n_e(v) = -1$ if $v$ is identified with $0$ along $e$. For $f\in D_{\nu_\omega}(\partial^\perp)$, the term $(f_e\eta_e')(v)$ is well defined as $f_e\eta_e' \in H^1(0, 1)$ and by embedding results, $f_e\eta_e'$ is continuous with continuous extension to $x=0$ and $1$. 

If $f$ and $\eta$ are more regular, for example $f_e \in H^1(0, 1)$ and $\eta_e\in H^2(0, 1)$, then
\begin{equation}\label{E:partial perp on graph}
	\partial^\perp f = \left(\frac{f_e'}{\eta_e'} + f_e \frac{\eta_e''}{(\eta_e')^2} \right)_{e\in E}.
\end{equation}
Thus, the operation $\partial^\perp$ can be understood as the operation $\star_{\omega}\partial f$ with an additional pertubation of multiplication-type that depends on the second derivative of $\eta$.

If $\eta$ is harmonic, then $\star_{\omega}\partial$ and $\partial^\perp$ agree as operations. However, for the domains we still have the change of vertex conditions. Harmonicity of $\eta$ implies that $\eta_e'$ is constant for each edge. Thus, the domain reduces to
\begin{equation}\label{E:domain partial perp on graph}
	D_{\nu_\omega}(\partial^\perp) = \{f \in \bigoplus_{e\in E} H^1(0, 1) : \sum_{e\in E} n_e(v)\eta_e'f_e(v) = 0 \text{ for all $v$}\}.
\end{equation}
Thus, $\mathcal F$ and $D_{\nu_\omega}(\partial^\perp)$ share the same differential regularity, but have different vertex conditions as $\mathcal F$ requires continuity in each vertex and $D_{\nu_\omega}(\partial^\perp)$ has weighted Kirchhoff-type vertex conditions.

If every vertex $v$ has an even number of edges attached to $v$, then $\mathcal F\subset D_{\nu_\omega}(\partial^\perp)$. In the contrary case, when $v$ is attached to an odd number of edges, any nontrivial function satisfying the Kirchhoff-vertex type condition at $v$ will not be continuous in $v$. Moreover, if $X$ contains loops, then one can construct functions in $D_{\nu_\omega}(\partial^\perp)$ associated with the loop, i.e., a function whose support coincides with the loop.

This example shows us that the difference between $\star_{\omega}\partial$ and $\partial^\perp_{V_0}$ can be characterized by:
\begin{enumerate}
	\item On a mutual domain, the operators $\star_{\omega}\partial$ and $\partial^\perp_{V_0}$ agree up to a multiplication operator that depends on $\partial^\ast_{V_0}\omega$ if $\omega$ is sufficiently regular.
	\item The domain of $\star_{\omega}\partial$ is characterized by a continuity condition in every junction point, while the domain $D_{\nu_{\omega}}(\partial^\perp_{V_0})$ is characterized by a weighted Kirchhoff-type condition in every junction point.
\end{enumerate}
The first statement is treated in  \Cref{P:F in domain of div}, on the mutual domain $\mathcal F$, and in \Cref{P:PnF in domain of div}, on the mutual domain given by certain piecewise energy finite functions $P_n\mathcal F$. The second statement can be seen from \Cref{P:partial u in partial ast}, which we write more clearly in \Cref{C:star partial u in partial perp}.

\subsection{Immediate consequences}\label{S:consequences}
Let us collect $L^2(K, \nu_\omega)$-versions of some results from the previous sections, which follow immediately from the bijectivity of $\star_\omega$. First, a variant of \Cref{P:partial u in partial ast} that says when the gradient $\star_{\omega}\partial^{(n)} u$ of a piecewise energy finite function $u\in P_n\mathcal F$ is in the domain $D_{\nu_{\omega}}(\partial^\perp_{V_0})$.
\begin{corollary}[\Cref{P:partial u in partial ast}]\label{C:star partial u in partial perp}
	Let $\omega\in \mathcal H$ be minimal energy-dominant. For $u\in P_n\mathcal F$ we have $\star_{\omega}\partial^{(n)}u\in D_{\nu_\omega}(\partial^\perp_{V_0})$ if and only if $u \in P_nD(\Delta_{\nu_\omega})$ and the normal derivatives of $u$ satisfy the matching condition \eqref{E:normal derivatives sum to 0} in $V_n\setminus V_0$, that is, $u$ satisfies \eqref{E:normal derivatives sum to 0} for all $q\in V_n\setminus V_0$. 
	
	Similarly, $\star_{\omega}\partial^{(n)}u\in D_{\nu_\omega}(\partial^\perp)$ if and only if $u\in P_nD(\Delta_{\nu_\omega})$ and $u$ satisfies \eqref{E:normal derivatives sum to 0} for every $q\in V_n$. In both cases
	\begin{equation*}
		\partial^\perp_{V_0}\star_{\omega} \partial^{(n)}u = - \partial^\ast_{V_0} \partial^{(n)}u = \Delta_{V_0}^{(n)} = \sum_{|w|=n} \mathds 1_{K_w}\Delta_{V_0}u^{(w)}.
	\end{equation*}
\end{corollary}
Next, the orthogonal basis for $\ker\partial^\ast$ in \Cref{T:onb} gives an orthogonal basis for $\ker\partial^\perp$.
\begin{corollary}[\Cref{T:onb}]\label{C:onb}
	Let $\omega\in \mathcal H$ be minimal energy-dominant. Then $\bigcup_{n=0}^\infty\bigcup_{|w| = n} \{\star_{\omega}\partial\psi_w \}$ is an orthogonal basis for $\ker\partial^\perp$ with $\norm{\star_{\omega}\partial\psi_w}_{L^2(K, \nu_\omega)}^2 = 18\cdot (5\slash 3)^{|w|}$. Moreover, for every $\Phi \in \ker\partial^\perp$, there are $\Phi_w\in\mathbb R$ such that
	\begin{equation*}
		\sum_{n=0}^\infty \sum_{|w|=n} \Phi_w^2 (5\slash 3)^{|w|} < \infty
	\end{equation*}
	and
	\begin{equation*}
		\Phi = \sum_{n=0}^\infty\sum_{|w|=n} \Phi_w \star_{\omega}\partial\psi_w.
	\end{equation*}
\end{corollary}
As a last consequence, we have a variant of the local Gauss-Green formula in \Cref{P:gauss green for div}. 
\begin{corollary}[Local Gauss-Green formula for $\partial^\perp_{V_0}$, see \Cref{P:gauss green for div}]\label{C:gauss green for div}
	Let $\omega\in\mathcal H$ be minimal energy-dominant. For $f\in D(\partial^\perp_{V_0})$ and $\varphi\in P_n\mathcal F$ we have
	\begin{equation*}
		(\partial^\perp_{V_0} f, \varphi)_{L^2(K, \nu_\omega)} = -(f, \star_{\omega}\partial^{(n)}\varphi)_{L^2(K, \nu_\omega)} + \sum_{|w|=n}\sum_{\partial K_w} \left[\vec n\cdot (f\omega)\right]  \varphi^{(w)}.
	\end{equation*}
\end{corollary}
Note, the normal parts of $f\omega$ exist for every $f\in D_{\nu_{\omega}}(\partial^\perp_{V_0})$ as, by definition, this is equivalent to $f\omega\in D_{\nu}(\partial^\ast_{V_0})$.

\subsection{Energy finite functions in $D_{\nu_\omega}(\partial^\perp_{V_0})$}\label{S:peff2}
Note that $\partial$ and $\partial^\ast_{V_0}$ behave similar to the gradient and divergence in a euclidean or riemannian setting, since they switch between scalar valued functions and \enquote{vector valued functions}, represented by the one-forms $\mathcal H$. Now that $\star_{\omega}\partial$ and $\partial^\perp_{V_0}$ are both $L^2(K, \nu_\omega)$-operators, we can compare them. We will see that, under some regularity of $\omega$, $\partial^\perp_{V_0}$ is just an extension of $\star_{\omega}\partial$. The following example shows the crucial difference between $\star_{\omega}\partial$ and $\partial^\perp_{V_0}$. We are in particular interested in the case when $\omega$ is a harmonic form.

In the following two propositions we want to show that a similar behaviour happens in the case of $\SG$. First, in \Cref{P:F in domain of div} we will show that $\partial^\perp_{V_0}$ behaves similar to the representation \eqref{E:partial perp on graph} if $\omega\in D_{\nu_\omega}(\partial^\ast_{V_0})$. Additionally, if  $\omega\in D_{\nu_\omega}(\partial^\ast_{V_0})$ then $\mathcal F\subset D_{\nu_\omega}(\partial^\perp_{V_0})$. This is due to the Gauss-Green formula \Cref{C:gauss green for div}, together with the continuity of energy finite functions and matching condition \eqref{E:SG matching condition normal part} for the normal parts of $\omega$. Thus, also here we can show that $\partial^\perp_{V_0}$ is an extension of $\star_{\omega}\partial$ up to a pertubation by a multiplication operator, which vanishes if $\omega$ is harmonic.

Then, in \Cref{P:PnF in domain of div} we study piecewise energy finite functions $P_n\mathcal F$ in $D_{\nu_\omega}(\partial^\perp_{V_0})$. We see that the continuity of $u\in P_n\mathcal F$ in all of the junction points besides $\ker\vec n\cdot\omega\cup V_0$ is sufficient and necessary. Here we want to add, that this agrees with the Kirchoff-type vertex conditions \eqref{E:domain partial perp on graph} in the metric graph example. In the metric graph example, the Kirchhoff-type vertex condition in a vertex $v$ are equivalent to continuity in $v$ if and only if the amount of edges at $v$ are exactly two. In our standard Sierpi\'{n}ski gasket setting, at each junction point are exactly two cells. In contrast, the corresponding result of \Cref{P:PnF in domain of div} for the \emph{level-3 Sierpi\'{n}ski gasket} is that for $u\in P_n\mathcal F$ to be in $D_{\nu_\omega}(\partial^\perp_{V_0})$ it needs to have continuity in junction points where two cells touch and \enquote{true} Kirchhoff-type conditions whenever three cells touch.
\begin{proposition}\label{P:F in domain of div}
	Let $\omega\in\mathcal H$ be minimal energy-dominant. Then $\mathcal F\subset D(\partial^\perp_{V_0})$ if and only if $\omega\in D_{\nu_\omega}(\partial^\ast_{V_0})$. In this case
	\begin{equation*}
		\partial^\perp_{V_0} u = \star_{\omega}\partial u - u\partial^\ast_{V_0} \omega \quad\text{for $u\in\mathcal F$}.
	\end{equation*}
	Similarly, we have $\mathcal F\subset D(\partial^\perp)$ iff $\omega\in D_{\nu_\omega}(\partial^\ast)$.
\end{proposition}
\begin{proof}
	First, for $u, \varphi\in\mathcal F$  the product rule for $\partial$ implies
	\begin{align*}
		(u, \star_{\omega}\partial\varphi)_{L^2(K, \nu_\omega)} &= \langle u\partial\varphi, \omega\rangle_{ \mathcal H} = \langle \partial(u\varphi), \omega\rangle_{ \mathcal H} - \langle \varphi\partial u, \omega\rangle_{ \mathcal H}\\
		&=\langle \partial(u\varphi), \omega\rangle_{ \mathcal H} - (\star_{\omega}\partial u, \varphi)_{L^2(K, \nu_\omega)}.
	\end{align*}
	Now, if $\omega\in D(\partial^\ast_{V_0})$, then for all $\varphi\in\mathcal F_{V_0}$ we have $u\varphi\in \mathcal F_{V_0}$ and thus
	\begin{equation*}
		(u, \star_{\omega}\partial\varphi)_{L^2(K, \nu_\omega)} = (u\partial^\ast_{V_0} \omega - \star_{\omega}\partial u, \varphi)_{L^2(K, \nu_\omega)}.
	\end{equation*}
	Hence, proving $\mathcal F\subset D(\partial^\perp_{V_0})$ with $\partial^\perp_{V_0} u = \star_{\omega}\partial u - u\partial^\ast_{V_0}\omega$.
	
	Now let us assume that $\mathcal F\subset D(\partial^\perp_{V_0})$. Then, we have for $u\in\mathcal F$ and $\varphi\in\mathcal F_{V_0}$
	\begin{equation*}
		-(\partial^\perp_{V_0} u, \varphi)_{L^2(K, \nu_\omega)} = (u, \star_{\omega}\partial\varphi)_{L^2(K, \nu_\omega)} = \langle \partial(u\varphi), \omega\rangle_{ \mathcal H} - (\star_{\omega}\partial u, \varphi)_{L^2(K, \nu_\omega)}.
	\end{equation*}
	Now, choosing $u = 1$ and using $\partial 1 = 0$, we conclude 
	\begin{equation*}
		-(\partial^\perp_{V_0} 1, \varphi)_{L^2(K, \nu_\omega)} = \langle \partial\varphi, \omega\rangle_{ \mathcal H} \quad\text{for all $\varphi\in\mathcal F_{V_0}$}.
	\end{equation*}
	Thus $\omega\in D(\partial^\ast_{V_0})$ with $\partial^\ast_{V_0} \omega = -\partial^\perp_{V_0} 1$. 
	
	The proof of the second equivalence is analogous.
\end{proof}

\begin{proposition}\label{P:PnF in domain of div}
	Let $\omega\in \mathcal H$ be minimal energy-dominant such that $\omega\in D_{\nu_\omega}(\partial^\ast_{V_0})$. Then $u\in P_n\mathcal F$ is an element of $D_{\nu_\omega}(\partial^\perp_{V_0})$ if and only if $u$ is continuous on $V_n\setminus \ker \ntrh\omega$. In this case
	\begin{equation*}
		\partial^\perp_{V_0} u = \star_{\omega}\partial^{(n)} u - u\partial^\ast_{V_0} \omega.
	\end{equation*}
	Similarly, $u\in D(\partial^\perp)$ iff $u$ is continuous on $V_n\setminus \ker\ntrh\omega$ and $u = 0$ on $V_0\setminus \ker\vec n\cdot\omega$.
\end{proposition}
\begin{proof}
	Let $u\in P_n\mathcal F$ and $\varphi\in\mathcal F$. Obviously $u\varphi\in P_n\mathcal F$. Using the product rule of $\partial^{(n)}$ and Gauss-Green for $\partial^\ast_{V_0}$ in \Cref{P:gauss green for div}, we derive
	\begin{align*}
		(u, \star_{\omega}\partial\varphi)_{L^2(K, \nu_\omega)} &= \langle \omega, u\partial\varphi\rangle_{ \mathcal H} =\langle \omega, \partial^{(n)}(u\varphi)\rangle_{ \mathcal H} - \langle \omega, \varphi \partial^{(n)} u\rangle_{ \mathcal H}\\
		&=(\partial^\ast_{V_0}\omega, u\varphi)_{L^2(K, \nu_\omega)} + \sum_{|w|=n} \sum_{\partial K_w} u^{(w)} \varphi \ntrh\omega - (\star \partial^{(n)}u, \varphi)_{L^2(K, \nu_\omega)}.
	\end{align*}
	If $u\in D(\partial^\perp_{V_0})$, then $(u, \star_{\omega}\partial\varphi)_{L^2(K, \nu_w)} = -(\partial^\perp_{V_0} u, \varphi)_{L^2(K, \nu_\omega)}$ for all $\varphi\in\mathcal F_{V_0}$. Thus, in above equation we can approximate $\mathds 1_{\{q\}}$, $q\in V_n\setminus V_0$, with $\varphi\in\mathcal F_{V_0}$ to get that $\sum_{|w|=n} \ntrh\omega (q) u^{(w)}(q)\mathds 1_{K_w}(q) = 0$. Since $\ntrh \omega$ satisfies the matching condition \eqref{E:SG matching condition normal part}, we conclude continuity of $u$ for each $q \in V_n\setminus\ker\ntrh\omega$.
	
	If $u$ is continuous on $V_n\setminus\ker\ntrh\omega$, then $\sum_{|w|=n} \sum_{\partial K_w} u^{(w)} \varphi \ntrh\omega = 0$ for each $\varphi\in\mathcal F_{V_0}$. Thus $u\in D(\partial^\perp_{V_0})$ with the stated equality.	The proof of the second equivalence is analogous.
\end{proof}
\begin{remark}\label{R:discontinuity on cut}
	\Cref{P:PnF in domain of div} implies that $D_{\nu_\omega}(\partial^\perp_{V_0})$ is generally larger than $\mathcal F$, as it allows for functions in $P_n\mathcal F$ that are discontinuous on $\ker\ntrh\omega \subset V_\ast$. For a precise example, see \Cref{S:discontinuity on cut}.
\end{remark}

\section{Pointwise representatives}\label{S:pointwise representative}
In this section we aim to prove our main result.
\begin{theorem}\label{T:pointwise representation}
	Let $\omega\in \ker\partial^\ast_{V_0}$ be minimal energy-dominant and $f\in D_{\nu_\omega}(\partial^\perp_{V_0})$. If $f$ and $\omega$ have a finite loop expansion, then
	\begin{equation}\label{E:pointwise representation}
		\lim_{m\to\infty} \dashint_{K_{wi^m}} f \ d\nu_{\omega} = \frac{\vec n \cdot (f\omega)(F_wq_i)}{\vec n\cdot \omega(F_wq_i)} \quad\text{for $F_wq_i\notin\ker\vec n \cdot\omega$}.
	\end{equation}
	Here, $\dashint_{K_w} f \ d_{\nu_\omega}$ represents the mean integral $\nu_{\omega}(K_w)^{-1} \int_{K_w} f \ d\nu_{\omega}$. Note that $\{F_wq_i\}  = \bigcap_m K_{wi^m}$ and $\lim_{m\to\infty}\nu_{\omega}(K_{wi^m}) = 0$.
\end{theorem}
We prove the more general version \Cref{T:general pointwise representation} in \Cref{S:proof}. There, we allow $\omega$ to have infinitely many loops. For now, let us discuss the result of \Cref{T:pointwise representation}.

Clearly, if $f \in D_{\nu_\omega}(\partial^\perp_{V_0})$ in \Cref{T:pointwise representation} is continuous, then \eqref{E:pointwise representation} implies that $f$ coincides with the ratio $\vec n \cdot (f\omega) \slash \vec n \cdot\omega$ on $V_\ast \setminus \ker \vec n\cdot\omega$. Equivalently, the ratio extends to a continuous function on $\SG\setminus \ker \vec n \cdot \omega$ as $V_\ast$ is dense. Note that $\ker\vec n \cdot \omega$ cannot be too large for divergence-free minimal energy-dominant $\omega$ due to \Cref{P:gauss green for div} giving the formula $\sum_{\partial K_w} (\vec n\cdot \omega) \varphi = \nu_{\omega, \partial\varphi}(K_w)$ for any $\varphi\in\mathcal F$ and cell $K_w$.

On the other hand, if the ratio in \eqref{E:pointwise representation} is discontinuous in $q\in V_\ast\setminus \ker\vec n \cdot\omega$, then $f$ is discontinuous in $q$. We use this fact in the examples in \Cref{S:discontinuity for loops} and \Cref{S:discontinuity in general} to construct discontinuous functions in $D_{\nu_\omega}(\partial^\perp_{V_0})$.

Nonetheless, for any $f$ with a finite loop expansion the formula \eqref{E:pointwise representation} says that the value of $f$ at a junction point $q$ can be represented by the ratio $\vec n \cdot (f\omega) \slash \vec n \cdot\omega$. This can be used to get an understanding of the Hodge star operator $\star_\omega f = f\omega$ and the fiber representation of $\mathcal H$ from \cite{ebe99, hrt13}: Given $\eta\in\mathcal H$ such that $\eta = f\omega$. Then, by \cite[Appendix D]{ebe99} or \cite[Theorem 2.1]{hrt13}, for every minimal energy-dominant measure $\nu$ there is a measurable field of Hilbert spaces $\{\mathcal H_x\}_{x\in K}$ such that $\mathcal H$ is isometrically isomorphic to $L^2(K, \{\mathcal H_x\}_{x\in K}, \nu)$. Let $\eta_x\in\mathcal H_x$ denote the fiber representation of $\eta$ in $L^2(K, \{\mathcal H_x\}_{x\in K}, \nu_{\omega})$. From the representation $\eta = f\omega$ we get that $\langle \eta_x, \omega_x\rangle_{\mathcal H_x} = f(x)$ $\nu_{\omega}$-almost everywhere and in combination with the Hodge star operator $f= \star_{\omega}\eta$ $\nu_{\omega}$-almost everywhere. In light of \eqref{E:pointwise representation} the fiber representation $\eta_x$ of $\eta$ with respect to the refence measure $\nu_\omega$ may be informally written as
\begin{equation*}
	\langle \eta_x, \omega_x\rangle_{\mathcal H_x} = \frac{\vec n \cdot \eta}{\vec n \cdot \omega} (x).
\end{equation*}
Intuitively speaking, \eqref{E:pointwise representation} says that the fiber representation $\eta_x$ of $\eta$ with respect to $\nu_{\omega}$ can be thought of being the flow of $\eta$ through $x$ rescaled by the flow of $\omega$ through $x$.

Let us now assume that $f$ is the gradient of a potential, that is, $f = \star_{\omega} \partial u$ for some $u\in D_{\nu_\omega}(\Delta_{V_0})$. Then,  \eqref{E:pointwise representation} reads
\begin{equation*}\label{E:pointwise representation partial u}
	\lim_{m\to\infty} \dashint_{K_{wi^m}} \star_\omega\partial u \ d\nu_{\omega} = \frac{\partial_n u(F_wq_i)}{\vec n\cdot \omega(F_wq_i)} \quad\text{for $F_wq_i\notin\ker\vec n \cdot\omega$}.
\end{equation*}
Note that $\langle (\partial u)_x, \omega_x\rangle_{\mathcal H_x} = (\star_{\omega}\partial u)(x)$ as $\partial u = (\star_{\omega}\partial u)\omega$. The fiber representation $(\partial u)_x$ of $\partial u$ with respect to $\nu_\omega$ may be informally written as
\begin{equation*}
	\langle (\partial u)_x, \omega_x\rangle_{\mathcal H_x} = \frac{\partial_n u}{\vec n \cdot \omega} (x).
\end{equation*}
If $\omega = \partial g$, for some $g\in\mathcal F$, then $\vec n \cdot \omega = \partial_n g$. Thus, in light of the total derivatives $\frac{du}{dg}$ from \cite[Theorem 5.4]{hino10} and the relation \eqref{eq:hino derivative}, we have the informal identification
\begin{equation*}
	\langle (\partial u)_x, (\partial g)_x\rangle_{\mathcal H_x} = \frac{\partial_n u}{\partial_n g}(x) = \frac{du}{dg}(x) = (\star_{\partial g}\partial u)(x).
\end{equation*}
Of course, these informal identities are not rigorous as the limit \eqref{E:pointwise representation} only holds on $V_\ast\setminus \ker\vec n \cdot \omega$, which is a $\nu_\omega$-null set. Nonetheless, these formulas may give some hints on the fibers of the first order structure, the Hodge star operator, and their relation to the first order differential operators, normal derivatives\slash parts, and total derivatives.

\begin{remark}\label{R:SG pontwise representative}
	An ideologically similar problem comes from stochastic analysis. In \cite[Theorem 5.1]{nak85}, see also \cite[Theorem 9.1]{hrt13}, it was proven that there is an isomorphism between $\mathcal H$ and $\mathring{\mathcal M}$, the space of martingale additive functionals of finite energy associated with the Hunt process of $(\mathcal E, \mathcal F)$. In particular, for $u\in \mathcal F$ there is a martingale additive functional of finite energy $M^{u}$ and the isomorphism relates $f\partial u\in\mathcal H$ with the stochastic integral $f \bullet M^{u}$, for any real-valued continuous functions $f$. Now, on the Sierpi\'{n}ski gasket, due to \cite{kus89}, we know that the martingale dimension is 1 and thus, there exists a basis martingale $M$ such that every other martingale additive functional of finite energy can be written as a stochastic integral with respect to $M$. Due to the isomorphism between $\mathcal H$ and $\mathring{\mathcal M}$, there is $\omega\in\mathcal H$ such that $M = M^{\omega}$. Moreover, there exists $u'$ such that $M^{u} = u'\bullet M^{\omega}$. Understanding the unknown integrand $u'$ is an important and difficult task in stochastic analysis and is answered in the Gaussian case by the Clark-Ocone formula from Malliavin calculus, see for example \cite[Section 6.3.1]{nua06}. Note, in our formulation $u' = \star_{\omega}\partial u = \langle (\partial u)_x, \omega_x\rangle_{ \mathcal H}$.
\end{remark}
Next, let us discuss a possible implication of \Cref{T:pointwise representation} on the \enquote{limited continuity} of functions $f\in D_{\nu_\omega}(\partial^\perp_{V_0})$ that satisfy \eqref{E:pointwise representation}. From the approximation property of the normal part in \Cref{C:convergence of ntrh baby version} we derive the following convergence of the ratio of normal parts along the sides.
\begin{corollary}\label{C:convergence of trace}
	Let $\omega\in \ker\partial^\ast_{V_0}$ minimal energy-dominant and $f\in D_{\nu_\omega}(\partial^\perp_{V_0})$, both consisting of finitely many loops. Then,
	\begin{equation*}
		\lim_{m\to\infty} \frac{\vec n\cdot(f\omega)(F_{wi^m}q_i)}{\vec n\cdot \omega(F_{wi^m}q_i)} = \frac{\vec n\cdot(f\omega)(F_{w}q_i)}{\vec n\cdot \omega(F_{w}q_i)} \quad\text{for $F_wq_i\notin\ker\vec n \cdot\omega$}.
	\end{equation*}
\end{corollary}
The convergence in \Cref{C:convergence of trace} suggests that the ratio
\begin{equation*}
	R_f(q) := \frac{\vec n\cdot(f\omega)(F_{w}q_i)}{\vec n\cdot \omega(F_{w}q_i)}\quad\text{for $q=F_wq_i \in V_\ast \setminus \ker\vec n \cdot\omega$}
\end{equation*}
can be extended to a continuous function on any metric graph approximation of $\SG$ with vertex set $V_n \setminus \ker\vec n \cdot\omega$. This is similar to \cite[Theorem 5.5]{bhs14}, where the authors show that the Radon-Nikodym derivative $\frac{d\nu_{f, g}}{d\nu}$, for any energy finite $f, g$, and Kusuoka measure $\nu$, restricted to the triangle spanned by the boundary of an arbitrary cell is continuous. The similarity comes from the observation: For $u\in\mathcal F$ and $\omega = \partial h$, where $h\in\mathcal F$, we have
\begin{equation*}
	f := \star_\omega \partial u = \frac{d\nu_{u, h}}{d\nu_h}.
\end{equation*}
In particular, \eqref{E:pointwise representation} is then equivalent to
\begin{equation*}
	\lim_{m\to\infty }\frac{\nu_{u, h}(K_{wi^m})}{\nu_h(K_{wi^m})} = \frac{\partial_n u(F_wq_i)}{\partial_n h(F_wq_i)} \quad\text{for $F_wq_i\notin\ker\partial_n h$}.
\end{equation*}
Thus, the ratio $R_f$ coincides with the \enquote{correct} pointwise definition of the Radon-Nikodym derivative $\frac{d\nu_{u, h}}{d\nu_h}$. However, to make these implications rigorous and therefore prove a variant of \cite[Theorem 5.5]{bhs14}, we need to improve \Cref{C:convergence of trace} (in fact it suffices to improve \eqref{E:simple partialn u approximation}) so that the limit is true for any sequence approximating $F_wq_i$ that lies on the sides of the triangle spanned by $\partial K_w$. We expect that such a result is true but we let it open for future research. However, our examples in \Cref{S:examples} show already certain limits of this approach. The example in \Cref{S:discontinuity for loops} shows that it cannot hold for any $f\in D_{\nu_\omega}(\partial^\perp_{V_0})$ with an infinite loop expansion. The example \Cref{S:discontinuity in general} shows that the limited continuity cannot be improved to a continuity on the whole gasket, thus, it is truly only a \enquote{limited continuity}.

\subsection{Two simple lemmata}
The first lemma we are about to prove allows us to reduce the proof of the general main theorem \Cref{T:general pointwise representation} to the case $q = q_0$ by using the scaling relation \eqref{E:measure scaling} for $\omega\circ F_w$ from \Cref{T:omega Fw}. The second lemma is a representation for harmonic functions on $\SG$ in terms of its normal and tangential derivative in $q_0$.  The latter is well known, see \cite[Lemma 2.7.3]{def}. We just add some simple computations for easy reference.

\begin{lemma}\label{L:spartial u and Fw}
	Let $\omega\in D_{\nu_\omega}(\partial^\ast_{V_0})$ be minimal energy-dominant. Then, $\omega\circ F_w$ is  minimal energy-dominant and for all $u\in\mathcal F$ and $w$ we have
	\begin{equation}
		\star_{\omega} \partial u \circ F_w =  \star_{\omega\circ F_w}\partial(u\circ F_w) \quad\text{$\nu_{\omega\circ F_w}$-almost everywhere}.
	\end{equation}
\end{lemma}
\begin{proof}
	The fact that $\omega\circ F_w$ is minimal energy-dominant follows from the energy measure scaling \eqref{E:measure scaling} as it gives $(\nu_\omega)_{F_w^{-1}} = r^{-|w|}\nu_{\omega\circ F_w}$. Given arbitrary $\varphi\in C(K)$.  Using the energy measure scaling \eqref{E:measure scaling} we have
	\begin{align*}
		\int_{K} (\star_{\omega} \partial u \circ F_w )\cdot \varphi \ d\nu_{\omega\circ F_w} &=\int_{K_w} \star_{\omega}\partial u \cdot (\varphi\circ F_w^{-1}) \ d(\nu_{\omega\circ F_w} )_{F_w} =r^{|w|}\int_{K_w} \star_{\omega}\partial u \cdot (\varphi\circ F_w^{-1}) \ d\nu_\omega \\
		&=r^{|w|} \langle \mathds 1_{K_w}\partial u, (\varphi\circ F_w^{-1})\omega\rangle_{ \mathcal H}=r^{|w|} \int_{K_w}\varphi\circ F_w^{-1} \ d\nu_{\partial u, \omega}\\
		&=r^{|w|} \int_{K}\varphi \ d(\nu_{\partial u, \omega})_{F_w{-1}} =\int_{K}\varphi \ d(\nu_{\partial u\circ F_w, \omega\circ F_w}).
	\end{align*}
	Since $\partial u\circ F_w = \partial(u\circ F_w)$ by definition, we conclude the proof from
	\begin{equation*}
		\int_{K}\varphi \ d(\nu_{\partial u\circ F_w, \omega\circ F_w}) = \langle \varphi \partial(u\circ F_w), \omega\circ F_w\rangle_{ \mathcal H} = \int_{K}  (\star_{\omega\circ F_w} \partial (u\circ F_w))\cdot \varphi \ d\nu_{\omega\circ F_w}.
	\end{equation*}
\end{proof}

\begin{lemma}[Lemma 2.7.3, \cite{def}]\label{L:sym antisym decomposition}
	Let $h_s$ be the harmonic function that is symmetric with respect to $q_0$ with boundary values $(0,1,1)$ and $h_a$ the harmonic function that is anti-symmetric with respect to $q_0$ with boundary values $(0,1,-1)$. Then we can write any harmonic funcion $u\in\mathcal F$ as
	\begin{equation*}
		u = u_s h_s + u_a h_a + c
	\end{equation*}
	where
	\begin{equation*}
		u_s = -\frac 12\partial_n u(q_0), \quad u_a = \frac 12(u(q_1) - u(q_2)), \quad c = u(q_0).
	\end{equation*}
	In particular, we have
	\begin{equation*}
		u(F_{0^n}q_1) = r^{n}u_s  + 5^{-n} u_a \quad\text{and}\quad  u(F_{0^n}q_2) = r^{n}u_s - 5^{-n} u_a,
	\end{equation*}
	and
	\begin{equation*}
		\partial_n u(F_{0^n}q_1) = u_s + 3^{1-n} u_a \quad\text{and}\quad \partial_n u(F_{0^n}q_2) = u_s - 3^{1-n} u_a.
	\end{equation*}
\end{lemma}
\begin{remark}
	Note, $u_a$ agrees with the tangential derivative $\partial_T u(q_0)$, that is, $u_a = \frac 12 \partial_Tu(q_0)$. Thus, in above formula $u$ is represented by its $1$-jet $(u(q_0), \partial_n u(q_0), \partial_T u(q_0))$ as in \cite[Lemma 2.7.3]{def}.
\end{remark}
\begin{proof}
	It is clear that $h_s, h_a$ and $1$ are three linearly independent harmonic functions on $\SG$. Thus we can write any harmonic function as a linear combination of those. On $V_0$ we get
	\begin{equation*}
		u(q_0) = c, \quad u(q_1) = u_s + u_a + c, \quad u(q_2) = u_s - u_a + c.
	\end{equation*}
	A simple calculation using $\partial_n u(q_0) = 2u(q_0) - u(q_1) - u(q_2)$ implies the claimed representation.
	
	For the values in $F_{0^n}q_i$ we see that $h_s(F_{0^n}q_1) = h_s(F_{0^n}q_2) = (3\slash 5)^n$ and $h_a(F_{0^n}q_1) = -h_a(F_{0^n}q_2) = 5^{-n}$. Further, the normal derivatives of $h_s$ and $h_a$ seen from the cell $K_{0^n}$ are
	\begin{equation*}
		\partial_n h_s(F_{0^n}q_1)= \partial_n h_s(F_{0^n}q_2) =  r^{-n}\left(2\cdot \left(\frac{3}{5} \right)^n -\left(\frac{3}{5} \right)^n - 0 \right) = 1 
	\end{equation*}
	and
	\begin{equation*}
		\partial_n h_a(F_{0^n}q_1) = - \partial_n h_a(F_{0^n}q_2) = r^{-n}\left(2\cdot 5^{-n} - (-5^{-n}) - 0 \right) = 3\cdot 3^{-n}.
	\end{equation*}
	From here, we easily conclude the values of $u(F_{0^n}q_i)$ and $\partial_n u(F_{0^n}q_i)$.
\end{proof}

\subsection{Proof of main theorem}\label{S:proof}
Finally, we prove the following more general version of \Cref{T:pointwise representation}.
\begin{theorem}\label{T:general pointwise representation}
	Let $\omega\in \ker\partial^\ast_{V_0}$ be minimal energy-dominant with $\omega = \partial\eta + \Theta$ being its Hodge decomposition \eqref{E:Hodge decomposition}. Fix some $u\in D_{\nu_\omega}(\Delta_{V_0})$ and $q=F_wq_i \in V_\ast$. Assume that there exists $a>0$, $\theta\in (0, \sqrt{3\slash 5})$ and $N\in\mathbb N_0$ such that
	\begin{equation*}
		|\Theta_{w'}| \leq a\theta^{|w'|}  \quad\text{for all $|w'|\geq |w|+N$ with $K_{w'}\subset K_w$.}
	\end{equation*}
	Moreover, let $\chi \in \{0, 1, 2, 3\}$ satisfy
	\begin{equation*}
		\#\{|w'|=|w|+m \mid K_{w'}\subset K_{w} \text{ and } \Theta_{w'} \neq 0 \} \leq \chi^{m} \quad\text{for all $m\geq N$}.
	\end{equation*}
	In any of the following cases
	\begin{enumerate}
		\item $\chi = 0$,
		\item $\chi =1, 2$ and $\theta < (\frac{3}{5})^{3\slash 2} \approx 0.4648$,
		\item $\chi = 3$ and $\theta < \sqrt{\frac 1 5} \approx 0.4472$,
	\end{enumerate}
	we have
	\begin{equation*}
		\lim_{m\to\infty} \dashint_{K_{wi^m}} \star_{\omega}\partial u \ d\nu_\omega = \frac{\partial_n u(F_wq_i)}{\ntrh \omega(F_wq_i)} \quad\text{for $F_wq_i \notin \ker\vec n \cdot \omega$}.
	\end{equation*}
\end{theorem}
If $f\in D_{\nu_\omega}(\partial^\perp_{V_0})$ has a finite loop expansion, then there exists $n$ such that $f\circ F_w$, for any $|w|=n$, can be represented by $\star_{\omega\circ F_w}\partial u_w$ for $u_w\in\mathcal F$. This follows from the Hodge decomposition \eqref{eq L2 hodge}, the basis from \Cref{C:onb}, self-similarity in \Cref{T:omega Fw}, and \Cref{L:spartial u and Fw}. Thus, after a change of measure and the use of the energy measure scaling relation \eqref{E:measure scaling} in  \Cref{T:omega Fw}, we can apply \Cref{T:general pointwise representation} to every $\star_{\omega\circ F_w}\partial u_w$. Reversing the zoom of $f$ to $f\circ F_w$, we arrive at the identity
\begin{equation*}
	\lim_{m\to\infty} \dashint_{K_{wi^m}} f \ d\nu_\omega = \frac{\vec n \cdot (f\omega) (F_wq_i)}{\ntrh \omega(F_wq_i)}.
\end{equation*}
Note, the case $\chi = 0$ means $\Theta$ contains finitely many loops close to $F_wq_i$. Thus, \Cref{T:general pointwise representation} $(i)$ implies \Cref{T:pointwise representation}. The other cases for $\chi\geq 1$ allow not only for infinitely many loops close to $F_wq_i$, but also count the amount of surrounding loops in the $\ker\partial^\perp$-part of $f$ that are infinitely close to $F_wq_i$.

Similar to the implication of \Cref{T:pointwise representation} from \Cref{T:general pointwise representation}, one can extend \Cref{T:general pointwise representation} to $f\in D_{\nu_\omega}(\partial^\perp_{V_0})$ with infinite loop expansion. The crucial step in the proof will be the question if the limit $m\to\infty$, the integral $\dashint_{K_{wi^m}}$, and the series in the $\ker\partial^\perp$-part of $f$ can be interchanged. However, we will not prove this here.

\begin{proof}[Proof of \Cref{T:general pointwise representation}]
	First, let us show that we can reduce it to the case $q = q_i$. Due to symmetry, it is then enough to assume $q = q_0$. Fix a word $w$, $N\in\mathbb N_0$, and $m\geq N$ as in \Cref{T:general pointwise representation}. Since $\nu_{\omega}(K_{wi^m}) = (\nu_\omega)_{F_w^{-1}}(K_{i^m})$ we have
	\begin{equation*}
		\dashint_{K_{wi^m}} \star_{\omega}\partial u \ d\nu_\omega = \frac{\int_{K_{wi^m}} \star_{\omega}\partial u \ d\nu_\omega}{\nu_{\omega}(K_{wi^m})} = \frac{\int_{K_{i^m}} \star_{\omega}\partial u \circ F_w \ d(\nu_\omega)_{F_w^{-1}}}{(\nu_\omega)_{F_w^{-1}}(K_{i^m})} = \dashint_{K_{i^m}} \star_{\omega}\partial u \circ F_w \ d(\nu_\omega)_{F_w^{-1}}.
	\end{equation*}
	Thus, \Cref{T:omega Fw} and \Cref{L:spartial u and Fw} imply
	\begin{equation*}
		\dashint_{K_{wi^m}} \star_{\omega}\partial u \ d\nu_\omega = \dashint_{K_{i^m}} \star_{\omega\circ F_w}\partial (u\circ F_w) \ d\nu_{\omega\circ F_w}.
	\end{equation*}
	Assume that the statement in \Cref{T:general pointwise representation} is true for $q_i$, that is, $w=\emptyset$. Then, we conclude
	\begin{equation*}
		\lim_{m\to\infty}\dashint_{K_{wi^m}} \star_{\omega}\partial u \ d\nu_\omega = \lim_{m\to\infty} \dashint_{K_{i^m}} \star_{\omega\circ F_w}\partial (u\circ F_w) \ d\nu_{\omega\circ F_w}= \frac{\partial_n(u\circ F_w)(q_i)}{\ntrh(\omega\circ F_w)(q_i)}.
	\end{equation*}
	By $\partial_n\varphi(F_wq_i) = r^{-|w|}\partial_n(\varphi \circ F_w)(q_i) $ and $\ntrh(\omega\circ F_w)(q_i) = r^{-|w|}\partial_n\varphi(F_wq_i)$, we deduce the general result. \newline
	
	Hence, we can assume that $q=q_0$. For an arbitrary cell $K_w$ we know that $\mathds 1_{K_w} u \in P_{|w|}\mathcal F$ and so \Cref{P:gauss green for div} implies
	\begin{equation}\label{E:int spartial u}
		\int_{K_w} \star_{\omega}\partial u \ d\nu = \langle \mathds 1_{K_w}\partial u, \omega\rangle_{ \mathcal H} = \sum_{\partial K_{w}}  u (\ntrh \omega).
	\end{equation}
	Similarly, using the Hodge decomposition \eqref{E:Hodge decomposition} of $\omega$ and \Cref{P:gauss green for div} with the fact that $\supp \psi_{w'} = K_{w'}$ and $\psi_{w'} \in P_{|w'|+1}\mathcal F$, we have
	\begin{equation}\label{E:m(Kw)}
		\begin{aligned}
			\nu_\omega(K_w) &= \langle \mathds 1_{K_w}\omega, \omega\rangle_{ \mathcal H} = \langle \mathds 1_{K_w}\partial\eta, \omega\rangle_{ \mathcal H} + \sum_{w'} \Theta_{w'} \langle \mathds 1_{K_w}\partial^{(|w'|+1)}\psi_{w'}, \omega\rangle_{ \mathcal H} \\
			&= \sum_{\partial K_w}\eta (\ntrh \omega)  + \sum_{\substack{|w'| < |w| \\ K_w \subsetneq K_{w'}}} \Theta_{w'} \sum_{\partial K_w} \psi_{w'}(\ntrh\omega) + \sum_{\substack{|w'| \geq |w| \\ K_{w'} \subset K_{w}}} \Theta_{w'} \sum_{i}\sum_{\partial K_{w'i}} \psi_{w'}(\ntrh\omega).
		\end{aligned}
	\end{equation}
	Now, we rewrite each term in \eqref{E:int spartial u} and \eqref{E:m(Kw)} for $w = 0^m$, $m\geq N$. We know that $\ntrh \omega(q_0) = \partial_n\eta(q_0)$ and
	\begin{equation*}
		\ntrh\omega(F_{0^m}q_i) = \partial_n\eta(F_{0^m}q_i) + \sum_{k=0}^{m-1}\Theta_{0^k} \partial_n\psi_{0^k}(F_{0^m}q_i) \quad\text{for $i=1, 2$}.
	\end{equation*} 
	We calculate $\partial_n\psi_{0^k}(F_{0^m}q_1) = -\partial_n\psi_{0^k}(F_{0^m}q_2) = 15\cdot 3^{-m}\cdot 5^k$, see \Cref{S:discontinuity for loops} for details. Since $\omega\in\ker\partial^\ast_{V_0}$, we know that $\eta$ is harmonic. Thus, using \Cref{L:sym antisym decomposition} on $\eta$, we can rearrange the sum in $\eqref{E:int spartial u}$ to get
	\begin{equation}\label{E:Sum u ntr omega}
		\begin{aligned}
			\sum_{\partial K_{0^m}}  u (\ntrh \omega)  &= \frac 12(2u(q_0) - u(F_{0^m}q_1) - u(F_{0^m} q_2))\partial_n \eta(q_0)\\
			&\quad +\frac 32\cdot  3^{-m}(u(F_{0^m} q_1)-u(F_{0^m} q_2))\left((\eta(q_1) - \eta(q_2)) + 10\sum_{k=0}^{m-1}\Theta_{0^k} 5^k \right).
		\end{aligned}
	\end{equation}
	Let us use the notation for the tangential derivative, that is, $\partial_T \eta(q_0) = \eta(q_1) - \eta(q_2)$. Clearly, \eqref{E:Sum u ntr omega} is also true when $u$ is replaced by $\eta$. In this case we can additionally use the harmonicity of $\eta$. Thus, the first term in \eqref{E:m(Kw)} equals
	\begin{equation}\label{E:Sum eta ntr omega}
		\sum_{\partial K_{0^m}}\eta (\ntrh\omega) = \frac 12 r^m \partial_n\eta(q_0)^2 + \frac 32 \cdot 15^{-m}\partial_T \eta(q_0)\left(\partial_T \eta(q_0) + 10\sum_{k=0}^{m-1}\Theta_{0^k} 5^{k}\right).
	\end{equation}
	Next, we treat the second term in \eqref{E:m(Kw)}. In the case $|w'|<m$ such that $K_{0^m} \subset K_{w'}$, we know that $w' = 0^k$ for $k = 0, 1, \ldots, m-1$. Then, $\psi_{0^k}(F_{0^m}q_1) = -\psi_{0^k}(F_{0^m}q_2) = 5^{-(m-k-1)}$. With similar calculations as above, we deduce
	\begin{equation}\label{E:2nd term in m(K0n)}
		\begin{aligned}
			\sum_{\substack{|w'| < |w| \\ K_w \subsetneq K_{w'}}} \Theta_{w'} \sum_{\partial K_w} \psi_{w'}(\ntrh\omega) &=\sum_{k=0}^{m-1} \Theta_{0^k} \sum_{\partial K_{0^m}} \psi_{0^k}(\ntrh\omega)\\
			&=15\cdot 15^{-m}\left(\partial_T \eta(q_0) + 10\sum_{k=0}^{m-1} \Theta_{0^k} 5^k\right)\sum_{k=0}^{m-1} \Theta_{0^k} 5^k.
		\end{aligned}
	\end{equation}
	Hence, first and second term in \eqref{E:m(Kw)} can be rewritten by \eqref{E:Sum eta ntr omega} and \eqref{E:2nd term in m(K0n)} to get
	\begin{equation}\label{E:1+2 of m(K0n)}
		\begin{aligned}
			\sum_{\partial K_{0^m}}&\eta (\ntrh\omega) + \sum_{\substack{|w'| < |w| \\ K_w \subsetneq K_{w'}}} \Theta_{w'} \sum_{\partial K_w} \psi_{w'}(\ntrh\omega)\\
			&=  \frac 12 r^m \partial_n\eta(q_0)^2 + \frac 32 \cdot 15^{-m}\left(\partial_T \eta(q_0) +  10\sum_{k=0}^{m-1} \Theta_{0^k} 5^k\right)^2.
		\end{aligned}
	\end{equation}
	Now, we treat the last term in \eqref{E:m(Kw)}. 
	\begin{figure}[H]
		\centering
		\begin{tikzpicture}[scale=4]
			\draw (0, 0) -- (1, 0) -- (1/2,0.866) -- (0, 0);
			\draw (1/2, 0) -- (3/4, 0.433) -- (1/4, 0.433) -- (1/2, 0);
			\node[above] at (1/2, 0.866) {$a$};
			\node[left] at (0, 0) {$b$};
			\node[right] at (1, 0) {$c$};
			\node[below] at (1/2, 0) {$z$};
			\node[right] at (3/4, 0.433) {$y$};
			\node[left] at (1/4, 0.433) {$x$};
		\end{tikzpicture}
		\caption{An example cell $K_{w'}$.}
		\label{fig:abcgasket}
	\end{figure}
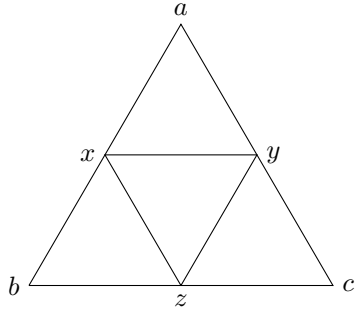
	In the case $|w'| \geq |w|$ with $K_{w'} \subset K_{0^m}$, the values of $\psi_{w'}$ on $\partial K_{w'i}$ are the same as $\psi$ on $V_1$. Assuming the cell $K_{w'}$ looks like in \Cref{fig:abcgasket}, we can use that $\ntrh\omega$ satisfies the matching condition \eqref{E:SG matching condition normal part}, to get
	\begin{equation}\label{E:sum psi ntr omega}
		\sum_i\sum_{\partial K_{w'i}} \psi_{w'} (\ntrh \omega) = 2\left((\swarrow\ntrh\omega)(x) - (\searrow\ntrh\omega)(y) + \frac 12(\leftarrow\ntrh\omega)(z) - \frac 12(\rightarrow\ntrh\omega)(z) \right).
	\end{equation}
	Here, the arrows denote the cell we use to calculate $\ntrh\omega$. For example, $(\rightarrow\ntrh\omega)(z)$ uses the cell left of $z$, that is, $K_{w'1}$.  Clearly, since $x, y, z\in V_{|w'|+1}\setminus V_{|w'|}$ we have
	\begin{equation*}
		\ntrh\omega(p) = \partial_n\eta(p) + \sum_{\substack{|\widetilde w| \leq |w'| \\ K_{w'} \subset K_{\widetilde w}}}\Theta_{\widetilde w} \partial_n\psi_{\widetilde w}(p) \quad\text{for $p = x, y, z$}.
	\end{equation*}
	First, assume we are in the case $\varphi = \eta$ or $\varphi = \psi_{\widetilde w}$ for $|\widetilde w| < |w'|$. This means $\varphi$ is harmonic on $K_{w'}$. Then,
	\begin{equation*}
		\partial_n^\swarrow\varphi(x) - \partial_n^\searrow\varphi(y) = r^{-(|w'|+1)}\left(\varphi(b) - \varphi(c) -2(\varphi(x) - \varphi(y) )\right) = r^{-|w'|}(\varphi(b) - \varphi(c))
	\end{equation*}	
	and
	\begin{equation*}
		\partial_{n}^\rightarrow\varphi(z) - \partial_{n}^\leftarrow\varphi(z) = -r^{-(|w'|+1)}(\varphi(x) - \varphi(y) + \varphi(b) - \varphi(c)) = -2r^{-|w'|}(\varphi(b) - \varphi(c)).
	\end{equation*}
	Thus, the corresponding terms in \eqref{E:sum psi ntr omega} are 0. Now, assume $\widetilde w = w'$. Then we can calculate for $p \in\{x, y, z \}$
	\begin{equation*}
		r^{(|w'|+1)}\partial_n\psi_{w'}(p) = \partial_n(\psi_{w'}\circ F_{w'i})(q_j) = \partial_n(\psi\circ F_i)(q_j) = \begin{cases}
			3 &, p = x = F_{w'0}q_1, \\
			-3 &, p = y = F_{w'0}q_2, \\
			3 &, p = z = F_{w'1}q_2.
		\end{cases}
	\end{equation*}
	Thus,
	\begin{equation*}
		\partial_n^\swarrow\psi_{w'}(x) - \partial_n^\searrow\psi_{w'}(y) + \partial_n^\rightarrow\psi_{w'}(z) =  r^{-(|w'|+1)}(3 - (-3) + 3 ) = 15r^{-|w'|}.
	\end{equation*}
	We conclude for the last term in \eqref{E:m(Kw)} using formula \eqref{E:sum psi ntr omega}:
	\begin{equation}\label{E: 3rd term in m(K0n)}
		\sum_{\substack{|w'| \geq |w| \\ K_{w'} \subset K_{w}}} \Theta_{w'} \sum_{i}\sum_{\partial K_{w'i}} \psi_{w'}(\ntrh\omega) = 30 \sum_{\substack{|w'| \geq |w| \\ K_{w'} \subset K_{w}}} \Theta_{w'}^2 r^{-|w'|}.
	\end{equation}
	In total, $\nu_\omega(K_{0^m})$ can be represented using \eqref{E:m(Kw)} together with \eqref{E:1+2 of m(K0n)} and \eqref{E: 3rd term in m(K0n)} to get
	\begin{equation}\label{E:1+2+3 term in m(K0n)}
		\nu_{\omega}(K_{0^m}) = \frac 12 r^m \partial_n\eta(q_0)^2 + \frac 32 \cdot 15^{-m}\left(\partial_T\eta(q_0)  +  10\sum_{k=0}^{m-1} \Theta_{0^k} 5^k\right)^2+ 30 \sum_{\substack{|w'| \geq m \\ K_{w'} \subset K_{0^m}}} \Theta_{w'}^2 r^{-|w'|}.
	\end{equation}
	Thus, \eqref{E:1+2+3 term in m(K0n)} together with \eqref{E:int spartial u} and \eqref{E:Sum u ntr omega} imply
	\begin{align*}
		&\dashint_{K_{0^m}} \star_{\omega}\partial u d\nu_\omega \\
		&= \frac{\frac 12(2u(q_0) - u(F_0^mq_1) - u(F_0^m q_2))\partial_n \eta(q_0) + \frac 32\cdot  3^{-m}(u(F_0^m q_1)-u(F_0^m q_2))\left(\partial_T\eta(q_0) + 10\sum_{k=0}^{m-1}\Theta_{0^k} 5^k \right)  }{\frac 12 r^m \partial_n\eta(q_0)^2 + \frac 32 \cdot 15^{-m}\left(\partial_T\eta(q_0) + 10\sum_{k=0}^{m-1}\Theta_{0^k} 5^{k} \right)^2 + 30 \sum_{\substack{|w'| \geq m \\ K_{w'} \subset K_{0^m}}} \Theta_{w'}^2 r^{-|w'|} }\\
		&=\left(\frac{r^{-m}(2u(q_0) - u(F_0^mq_1) - u(F_0^m q_2))}{\partial_n\eta(q_0)} +  3\cdot 3^{-m}r^{-m}(u(F_0^m q_1)-u(F_0^m q_2))\frac{\partial_T\eta(q_0) +10\sum_{k=0}^{m-1}\Theta_{0^k} 5^k }{\partial_n\eta(q_0)^2}\right)\\
		&\qquad\cdot \left(1 + 3\cdot r^{-m}\cdot 15^{-m}\frac{\left(\partial_T\eta(q_0) + 10\sum_{k=0}^{m-1}\Theta_{0^k} 5^{k}\right)^2}{\partial_n\eta(q_0)^2} + 60\cdot r^{-m} \frac{ \sum_{\substack{|w'| \geq m \\ K_{w'} \subset K_{0^m}}} \Theta_{w'}^2 r^{-|w'|}}{\partial_n\eta(q_0)^2} \right)^{-1}\\
		&=(\textbf{N}_1 + \textbf{N}_2 + \textbf{N}_3)\cdot (1 + \textbf{D}_1 + \textbf{D}_2)^{-1}.
	\end{align*}
	\normalsize
	The terms $\textbf{N}_2$ and $\textbf{N}_3$ are separated with respect to the sum $\partial_T\eta(q_0) +10\sum_{k=0}^{m-1}\Theta_{0^k} 5^k$. By definition of the normal derivative, we have
	\begin{equation}\label{E:lim N1}
		\lim_{m\to\infty} \textbf{N}_1 = \lim_{m\to\infty} \frac{r^{-m}(2u(q_0) - u(F_0^mq_1) - u(F_0^m q_2))}{\partial_n\eta(q_0)} = \frac{\partial_n u(q_0)}{\partial_n \eta(q_0)}.
	\end{equation}
	It remains to show, that $\textbf{N}_2, \textbf{N}_3$, $\textbf{D}_1$ and $\textbf{D}_2$ converge to 0. In $\textbf{N}_2$ we can use that $u\in D(\Delta_{\nu_{\omega}})$ satisfies the resistance estimate
	\begin{equation}\label{E:resistance estimate}
		|u(F_0^m q_1)-u(F_0^m q_2)| \leq C r^{m\slash 2},
	\end{equation}
	with a constant $C>0$ that depends linearly on $\norm{\Delta_{\nu_{\omega}} u}_{L^2(K_{0^m}, \nu_{\omega})}$. This imples
	\begin{equation}\label{E:lim N2}
		|\textbf{N}_2| \leq C'\cdot 3^{-m}r^{-m}|u(F_0^m q_1)-u(F_0^m q_2)|\leq C \left(\frac{\sqrt 5}{3\sqrt 3}\right)^m \to 0.
	\end{equation}
	For $\textbf{N}_3$ to converge to 0, we just distinguish the following two cases. Either, after some level there are no loops close to $q_0$, that is, $\chi = 0$, or the loops pile up close to $q_0$, that is, $\chi \geq 1$. If $\chi = 0$, then $\sum_{k=0}^{m-1} \Theta_{0^k} 5^k$ is uniformly bounded in $m$ and we can use the same rate as in \eqref{E:lim N2}. In the case $\chi \geq 1$, we use our assumption $|\Theta_{0^k}|\leq a\theta^k$ for some $0 \leq \theta < \sqrt{3\slash 5}$. So, $\left|\sum_{k=0}^{m-1} \Theta_{0^k} 5^k\right| \leq a \frac{(5\theta)^m -1}{5\theta -1}$, and in the case $\theta = 1\slash 5$, we have $|\sum_{k=0}^{m-1} \Theta_{0^k} 5^k| \leq Cm$. Using again the resistance estimate \eqref{E:resistance estimate}, we have
	\begin{equation*}
		|\textbf{N}_3| \leq C3^{-m}r^{-m}r^{-m\slash 2}\left(\max\{(5\theta)^m, m \} +1 \right) =  C\left(\frac{\sqrt 5}{3\sqrt 3}\right)^m\left(\max\{(5\theta)^m, m \} +1 \right).
	\end{equation*}
	As $3^{-m}r^{-m}r^{m\slash 2}(5\theta)^m = \left( \frac{5^{3\slash 2}}{3^{3\slash 2}} \theta \right)^m$, we get
	\begin{equation}\label{E:lim N3}
		\lim_{m\to\infty} \textbf{N}_3  = 0 \quad\text{if $\theta < \left(\frac 35\right)^{3\slash 2}$ and $\chi\geq 1$}.
	\end{equation}	
	Now, we treat the denominator. For $\textbf{D}_1$ we have
	\begin{equation*}
		|\textbf{D}_1| \leq C r^{-m}15^{-m}(1 + \max\{(5\theta)^m, m\}) = C 3^{-2m}(1 + \max\{(5\theta)^m, m\}).
	\end{equation*}
	As $3^{-2m}(5\theta)^m = \left(\frac{5}{3^2} \theta \right)^m$ and $3^2\slash 5 >1$ we always have
	\begin{equation}\label{E:lim D1}
		\lim_{m\to\infty} \textbf D_1 = 0.
	\end{equation} 
	
	For $\textbf{D}_2$, since $K_{0^m}$ contains at most $\chi^{k-m}$ many $(m+k)$-level sets, where the corresponding coefficients of $\Theta$ are not zero, we have
	\begin{equation*}
		r^{-m}\sum_{\substack{|w'| \geq m \\ K_{w'} \subset K_{0^m}}} \Theta_{w'}^2 r^{-|w'|} \leq a^2 r^{-m} \sum_{k=0}^\infty \chi^{k}\theta^{2(k+m)}r^{-(k+m)} = a^2\left(\frac{5}{3}\theta \right)^{2m} \sum_{k=0}^\infty \left(\frac{5}{3}\chi\theta^2\right)^k.
	\end{equation*}
	In particular, in the case $\chi = 0$ we always have $\textbf{D}_2 = 0$. In the remaining cases, we have
	\begin{equation}\label{E:lim D2}
		\lim_{m\to\infty} \textbf{D}_2 = 0 \quad\text{if $\theta < \min\left\{ \frac{3}{5}, \sqrt{\frac{3}{5\chi}} \right\}$ for $\chi \geq 1$.}
	\end{equation}
	Hence, \eqref{E:lim N1}, \eqref{E:lim N2}, \eqref{E:lim N3}, \eqref{E:lim D1} and \eqref{E:lim D2} imply the statement in the corresponding cases.

\end{proof}

\begin{remark}\label{R:tangential derivative in proof}
	Note that the terms $\textbf N_2+\textbf N_3$ and $\textbf D_1$ contain
	\begin{equation*}
		\partial_T\eta(q_0) + 10\sum_{k=0}^{m-1}\Theta_{0^k} 5^{k},
	\end{equation*}
	which converges to the tangential part $\vec t\cdot \omega(q_0)$ for $m\to\infty$, see \Cref{T:tangential part}. In the case $\chi = 0$, this term equals $\vec t\cdot \omega(q_0)$ for all $m\geq N$.
\end{remark}

\begin{remark}\label{R:pointwise bound tangential part}
	From the proof above one can calculate the rate of convergence for \Cref{T:general pointwise representation}. However, it  contains the factor $1\slash |\vec n \cdot\omega(F_wq_i)|$, which is in general very difficult to deal with. If $\omega\in\ker\partial^\ast_{V_0}$ contains only finitely many loops, one can get the following estimate:
	\begin{equation*}
		\begin{aligned}
			\left|\frac{\partial_n u(F_wq_i)}{\ntrh \omega(F_wq_i)}\right.&\left. - \dashint_{K_{wi^m}}\star_\omega\partial u \ d\nu_\omega \right| \\
			&\lesssim r^{\frac{|w|+m}{2}} + 3^{-|w|} \left(\frac{\sqrt{5}}{3\sqrt 3 } \right)^m \frac{|\vec t\cdot \omega(F_wq_i)|}{|\ntrh\omega(F_wq_i)|} + 3^{-2(|w|+m)} \frac{|\partial_n u(F_wq_i)|}{|\ntrh \omega(F_wq_i)|} \frac{\vec t\cdot\omega(F_wq_i)^2}{\ntrh\omega(F_wq_i)^2}.
		\end{aligned}
	\end{equation*}
\end{remark}

\section{Examples for discontinuous functions in the domain}\label{S:examples}
Due to \Cref{T:pointwise representation} we know, if $f\in D_{\nu_\omega}(\partial^\perp_{V_0})$ is continuous, then
\begin{equation*}
	f(q) = \frac{\vec n\cdot(f\omega)(q)}{\vec n \cdot \omega(q)} \quad\text{for $q\in V_\ast\setminus \ker\vec n\cdot\omega$}.
\end{equation*}
In the following, we construct examples where the continuity of the ratio $\vec n\cdot(f\omega)\slash \vec n \cdot \omega$ fails. More precisely, we provide the following examples.

In \Cref{S:discontinuity on cut}, for a special $\omega$, we extend \Cref{P:PnF in domain of div}, that is, $u\in P_n\mathcal F$ is in $D_{\nu_\omega}(\partial^\perp_{V_0})$ if and only if $u$ is continuous on $V_n\setminus(\ker\vec n\cdot\omega\cup V_0)$. We show for $\omega = \partial h$, where $h\in\mathcal F$ is the harmonic function with boundary values $(1, 0, 0)$, that $D_{\nu_\omega}(\partial^\perp_{V_0})$ contains arbitrary sums of functions $u= u_L + u_R$, where $u_L$ and $u_R$ are energy finite with respect to the left and right half of $\SG$, respectively. Note that in this example we have $\partial_n h = 0$ along the vertical line spanned by $q_0$ and $F_1q_2$.

In \Cref{S:discontinuity for loops} we study elements $\Phi\in \ker\partial^\perp$, in particular those who have infinitely many loops. We show that the growth rate of the coefficients $\Phi_w$ has to be smaller than $(3\slash 5)^{|w|}$, so that one can expect continuity of the ratio $\vec n\cdot(\Phi\omega)\slash \vec n\cdot \omega$ along the sides as in \Cref{C:convergence of trace}.

In \Cref{S:discontinuity in general} we show that continuity for general elements $f\in D_{\nu_\omega}(\partial^\perp_{V_0})$ is a more delicate problem when approaching $q$ along a vertical line, instead of the sides. Here we provide an example for discontinuous $\Phi\in\ker\partial^\perp$ with infinitely many and also finitely many loops. Also, we provide a discontinuous $\star_{\omega}\partial u$ for $u\in\mathcal F$ harmonic, where convergence along the vertical line fails.

\subsection{Discontinuity due to separation}\label{S:discontinuity on cut}
As explained in \Cref{R:discontinuity on cut}, from \Cref{P:PnF in domain of div} we know that $D_{\nu_\omega}(\partial^\perp_{V_0})$ is usually larger than $\mathcal F$ as it contains piecewise energy finite functions that may be only continuous on $V_\ast\setminus (\ker\ntrh\omega\cup V_0)$.

Consider the harmonic function $h\in\mathcal F$ given by the boundary values $(1, 0, 0)$, that is, $h(q_0) = 1$ and $h(q_1) = h(q_2) = 0$. Then $\omega := \partial h \in \ker\partial^\ast_{V_0}$ is minimal energy-dominant due to \cite{hino10}. Since $h$ is symmetric with respect to $q_0$, it is easy to see that 
\begin{equation*}
	\{F_{0^n1}q_2 \mid n\in \mathbb N_0 \} \subset \ker\partial_n h = \ker\ntrh\omega.
\end{equation*}
For this choice of $\omega$ we are going to prove that $D_{\nu_\omega}(\partial^\perp_{V_0})$ contains all functions that are the sum of two arbitary functions that are energy finite on the left and right half of $\SG$, but may be discontinuous along all points along the vertical line in $\SG$.

Let us denote by $W_L$ and $W_R$ the set of words representing the left and right half of $\SG$: 
\begin{equation*}
	W_L = \{1, 01, 001, 0001, \ldots \} \quad\text{and}\quad W_R = \{2, 02, 002, 0002, \ldots \}.
\end{equation*}
Denote $\SG_L := K_L := \bigcup_{w\in W_L} K_w$ and $\SG_R := K_R := \bigcup_{w\in W_R} K_w$. Note that $W_L\cup W_R$ is a partition of $\SG$. Let us define the space of energy finite functions on $\SG_L$ and $\SG_R$ via
\begin{equation*}
	\mathcal F_L := \left\{u \in C(\SG_L) : \sum_{w\in W_L} \mathcal E_{K_w}(u)<\infty \right\} \quad\text{and}\quad \mathcal F_R := \left\{u \in C(\SG_R) : \sum_{w\in W_R} \mathcal E_{K_w}(u)<\infty \right\}.
\end{equation*}
We now show that $\mathcal F_L \oplus \mathcal F_R \subset D_{\nu_\omega}(\partial^\perp_{V_0})$. Clearly, $\mathcal F_L \oplus \mathcal F_R$ contains functions that are discontinuous along the vertical line $\{F_{0^n1}q_2 \mid n\in \mathbb N_0 \}$.

\begin{proposition}
	Let $h\in\mathcal F$ be the harmonic function with boundary values $(1, 0, 0)$ and $\omega:=\partial h$.	Given $u_L \in \mathcal F_L$, $u_R\in\mathcal F_R$ and define $u := \mathds 1_{K_L} u_L + \mathds 1_{K_R}u_R$. Then $u\in D_{\nu_h}(\partial^\perp_{V_0})$.
\end{proposition}
\begin{proof}
	Clearly, $u\in L^2(K, \nu_h)$. Without loss of generality we consider energy finite extensions of $u_L$ and $u_R$ in $\mathcal F$. Fix arbitrary $\varphi\in \mathcal F$. We do the calculation for $u_L$, as it is analogous for $u_R$. We have
	\begin{align*}
		(\mathds 1_{K_L}u_L, \star_{\omega}\partial\varphi)_{L^2(K, \nu_h)} &= \sum_{w\in W_L} \langle \mathds 1_{K_w}u_L\partial \varphi, \partial h\rangle_{ \mathcal H} \\
		&=\sum_{w\in W_L} \langle \mathds 1_{K_w}\partial (u_L\varphi), \partial h\rangle_{ \mathcal H} - \langle \mathds 1_{K_w}\partial u_L, \varphi\partial h\rangle_{ \mathcal H} \\
		&=\sum_{w\in W_L}\mathcal E_{K_w}(u_L\varphi, h) - (\sum_{w\in W_L} \mathds 1_{K_w}\star_{\omega}\partial u_L, \varphi)_{L^2(K, \nu_h)}.
	\end{align*}
	Thus, if $\sum_{w\in W_L}\mathcal E_{K_w}(u_L\varphi, h) = 0$, then it is clear that $\partial^\perp_{V_0} u = \mathds 1_{K_L}\star_{\omega}\partial u_L + \mathds 1_{K_R}\star_{\omega}\partial u_R$ with
	\begin{equation*}
		\norm{\partial^\perp_{V_0} u}_{L^2(K, \nu_h)} = \sqrt{\sum_{w\in W_L\cup W_R} \mathcal E_{K_w}(u)} <\infty.
	\end{equation*}
	The harmonicity of $h$ implies
	\begin{equation*}
		\mathcal E_{K_w}(u_L\varphi, h) = \sum_{\partial K_w} u_L \varphi \partial_n h.
	\end{equation*}
	Because of our particular choice of $h$, we know that $\partial_n h = 0$ on $\SG_L\cap \SG_R\setminus \{q_0\}$. Thus, after summing over all $w\in W_L$, only $q_0, q_1\in V_0$ and all the vertices $q_0^{(n)} := F_{0^n1}q_0$ remain. We get
	\begin{equation*}
		\sum_{w\in W_L}\mathcal E_{K_w}(u_L\varphi, h) = \sum_{w\in W_L : |w| > n} \mathcal E_{K_w}(u_L\varphi, h) + u_L(q_1)\varphi(q_1)\partial_n h(q_1) + \sum_{k=1}^n u_L(q_0^{(k)} )\varphi(q_0^{(k)} )\partial_nh(q_0^{(k)}).
	\end{equation*}
	Here, we use that $u_L$ and $\varphi$ are continuous in $q_0^{(k)}$ and $\partial_n h$ satisfies the matching condition \eqref{E:matching condition}. Thus, the right sum is a telescoping series, implying for any $n$
	\begin{equation*}
		\sum_{w\in W_L}\mathcal E_{K_w}(u_L\varphi, h) = \sum_{w\in W_L : |w| > n} \mathcal E_{K_w}(u_L\varphi, h) + u_L(q_1)\varphi(q_1)\partial_n h(q_1) +  u_L(q_0^{(n)} )\varphi(q_0^{(n)} )\partial_nh(q_0^{(n)}).
	\end{equation*}
	Now consider $n\to\infty$. We know that $\sum_{w\in W_L : |w| > n} \mathcal E_{K_w}(u_L\varphi, h)$ converges to 0, as $u_L, \varphi$ and $h$ are energy finite with respect to $W_L$. Due to continuity, we see that $u_L(q_0^{(n)})\varphi(q_0^{(n)})\to u_L(q_0)\varphi(q_0)$. For $\partial_n h$ we use that $\partial_n h(q_0^{(n)}) \to \frac 12 \partial_n h(q_0)$ from \cite[Lemma 4.3]{ls14}. Hence
	\begin{equation*}
		\sum_{w\in W_L}\mathcal E_{K_w}(u_L\varphi, h) = u_L(q_1)\varphi(q_1)\partial_n h(q_1) + \frac 12 u_L(q_0 )\varphi(q_0 )\partial_nh(q_0).
	\end{equation*}
	If we now test with $\varphi\in\mathcal F_{V_0}$, the right-hand side vanishes. As the same is true for $u_R$, we conclude $u\in D(\partial^\perp_{V_0})$.
\end{proof}

\subsection{Discontinuity along the sides}\label{S:discontinuity for loops}
Let $h\in\mathcal F$ be any nonconstant harmonic function and consider
\begin{equation*}
	\Phi = \sum_{w} \Phi_w \star_{\partial h} \partial\psi_w,
\end{equation*}
with $\Phi_w$ as in \Cref{C:onb} so that $\Phi\in\ker\partial^\perp$. We study the convergence of the ratio
\begin{equation*}
	R\Phi(F_wq_i) := \frac{\vec n\cdot (\Phi\partial h)(F_wq_i)}{\vec n \cdot \partial h(F_wq_i)} = \sum_{w} \Phi_w\frac{\partial_n\psi_w(F_wq_i)}{\partial_n h(F_wq_i)},
\end{equation*}
when approaching $q_0$ along the sides. In particular, we are interested in the case when the loops pile up close to $q_0$. For simplicity, we only consider loops related to $q_0$:
\begin{equation}\label{E:only loops close to q0}
	\Phi_w \neq 0 \quad\text{iff}\quad w = 0^k, k\in\mathbb N_0.
\end{equation}
Moreover, we assume that there exists $a > 0$ and $\phi\in (-\sqrt{3\slash 5}, \sqrt{3\slash 5})$ such that
\begin{equation}\label{E:loops close to q0}
	\Phi_{0^k} = a\phi^k, k\in\mathbb N_0.
\end{equation}
We have the following result on the convergence of the ratio.

\begin{proposition}\label{P:discontinuity along sides using loops}
	Let $h\in\mathcal F$ be nonconstant harmonic such that $\partial_n h(q_0)\neq 0$ and $\Phi\in\ker\partial^\perp$ satisfying \eqref{E:only loops close to q0} and \eqref{E:loops close to q0}. Then,
	\begin{equation*}
		\lim_{n\to\infty}\frac{\vec n\cdot (\Phi\partial h)(F_{0^n}q_1)}{\vec n \cdot \partial h(F_{0^n}q_1)} =\begin{cases}
			\text{divergent} &\text{if $-\sqrt{3\slash 5} <\phi \leq -3\slash 5$,} \\
			0 &\text{if $|\phi| < 3\slash 5$,}\\
			-\frac{30 }{\partial_n h(q_0)} \frac{a}{5\phi -1} = -15\frac{a}{\partial_n h(q_0)} &\text{if $\phi = 3\slash 5$,}\\
			-\sgn  \frac{a}{\partial_n h(q_0)} \cdot \infty &\text{if $3\slash 5 < \phi < \sqrt{3\slash 5}$,}
		\end{cases}
	\end{equation*}
	and
	\begin{equation*}
		\lim_{n\to\infty}\frac{\vec n\cdot (\Phi\partial h)(F_{0^n}q_2)}{\vec n \cdot \partial h(F_{0^n}q_2)}  =\begin{cases}
			\text{divergent} &\text{if $-\sqrt{3\slash 5} <\phi \leq -3\slash 5$,} \\
			0 &\text{if $|\phi| < 3\slash 5$,}\\
			\frac{30 }{\partial_n h(q_0)} \frac{a}{5\phi -1} = 15\frac{a}{\partial_n h(q_0)} &\text{if $\phi = 3\slash 5$,}\\
			\sgn \frac{a}{\partial_n h(q_0)} \cdot \infty &\text{if $3\slash 5 < \phi < \sqrt{3\slash 5}$.}
		\end{cases}
	\end{equation*}
\end{proposition}
Note that for the ratio at $q_0$ we have $R\Phi(q_0) = 0$. Thus, $R\Phi$ converges along the sides to $R\Phi(q_0)$ only if the growth of the coefficients $\Phi_w$ is smaller than $(3\slash 5)^{|w|}$ close to $q_0$. If the growth of the coefficients is larger than $(3\slash 5)^{|w|}$, continuity of $R\Phi$ can not be expected. 

\begin{remark}\label{R:discontinuity}
	In \Cref{P:discontinuity along sides using loops}, condition \eqref{E:only loops close to q0} is not a restriction, as $\partial_n \psi_w(F_{0^n}q_i) = 0$, $i=1, 2$, for all $w\neq 0^k$, $k\in\mathbb N_0$.
\end{remark}

\begin{proof}[Proof of \Cref{P:discontinuity along sides using loops}]
	As $\supp\psi_w = K_w$, we know that
	\begin{equation*}
		R\Phi(F_{0^n} q_i) = \sum_{k=0}^{n-1} \Phi_{0^k} \frac{\partial_n\psi_{0^k}(F_{0^n}q_i)}{\partial_n h(F_{0^n}q_i)} \quad\text{for $i =1, 2$}.
	\end{equation*}
	Let us first calculate $\partial_n\psi_{0^k}(F_{0^n}q_i)$. We use that $\widetilde\psi := \psi\circ F_0$ is the harmonic function with boundary values $(0, 1, -1)$. We have
	\begin{equation*}
		\widetilde \psi(F_{0^n}q_1) = -\widetilde\psi(F_{0^n}q_2) = 5^{-n} \quad\text{and}\quad \partial_n \widetilde\psi(F_{0^n}q_1) = - \partial_n \widetilde\psi(F_{0^n}q_2) = 3^{1-n}.
	\end{equation*}
	Thus,
	\begin{equation}\label{E:normal derivative psi}
		\partial_n \psi_{0^k}(F_{0^n} q_1) = r^{-n}\partial_n(\psi_{0^k}\circ F_{0^n})(q_1) = r^{-n}\partial_n(\widetilde \psi \circ F_{0^{n-k-1}})(q_1) = r^{-k-1}\partial_n \widetilde\psi(F_{0^{n-k-1}}q_1) = 15\cdot \frac{5^k}{3^n}
	\end{equation}
	and 
	\begin{equation*}
		\partial_n \psi_{0^k}(F_{0^n} q_2) =  r^{-k-1}\partial_n \widetilde\psi(F_{0^{n-k-1}}q_2) = -r^{-k-1}\partial_n \widetilde\psi(F_{0^{n-k-1}}q_1) = -15\cdot \frac{5^k}{3^n}.
	\end{equation*}
	Using \Cref{L:sym antisym decomposition} for $\partial_n h(F_{0^n}q_1)$ we have
	\begin{equation*}
		R\Phi(F_{0^n}q_1) = \sum_{k=0}^{n-1} \Phi_{0^k}\frac{15\cdot \frac{5^k}{3^n}}{-\frac 12\partial_n h(q_0) + \frac 12 3^{-n}\partial_T h(q_0)} = \frac{-30}{\partial_n h(q_0) - 3^{-n}\partial_T h(q_0)} \cdot 3^{-n}\sum_{k=0}^{n-1}\Phi_{0^k} 5^k.
	\end{equation*}
	Due to our assumption $\Phi_{0^k} = a\phi^k$, in the case $\phi \neq \pm 1\slash 5$ we have
	\begin{equation*}
		3^{-n}\sum_{k=0}^{n-1}\Phi_{0^k} 5^k = \frac{a}{5\phi -1}\left(\left(\frac{5\phi}{3} \right)^n - 3^{-n} \right).
	\end{equation*}
	In the case $\phi = \pm 1\slash 5$, the sum is bounded by $n$ and $R\Phi(F_{0^n} q_1)$ converges to $0$. Hence, 
	\begin{equation*}
		\lim_{n\to\infty} R\Phi(F_{0^n}q_1) =\begin{cases}
			\text{divergent} &\text{if $-\sqrt{3\slash 5} <\phi \leq -3\slash 5$,} \\
			0 &\text{if $|\phi| < 3\slash 5$,}\\
			-\frac{30 }{\partial_n h(q_0)} \frac{a}{5\phi -1} = -15\frac{a}{\partial_n h(q_0)} &\text{if $\phi = 3\slash 5$,}\\
			-\sgn  \frac{a}{\partial_n h(q_0)} \cdot \infty &\text{if $3\slash 5 < \phi < \sqrt{3\slash 5}$.}
		\end{cases}
	\end{equation*}
	Similarly, using that $\partial_n\psi_{0^k}(F_{0^n}q_2) = - \partial_n\psi_{0^k}(F_{0^n}q_1)$ and $\partial_n h(F_{0^n}q_2) = -\frac 12\partial_n h(q_0) - \frac 12 3^{-n}\partial_Th(q_0)$, we have
	\begin{equation*}
		R\Phi(F_{0^n}q_2) = \frac{30}{\partial_n h(q_0) + 3^{-n}\partial_T h(q_0)} \cdot 3^{-n} \sum_{k=0}^{n-1}\Phi_{0^k} 5^k
	\end{equation*}
	and therefore
	\begin{equation*}
		\lim_{n\to\infty} R\Phi(F_{0^n}q_2) =\begin{cases}
			\text{divergent} &\text{if $-\sqrt{3\slash 5} <\phi \leq -3\slash 5$,} \\
			0 &\text{if $|\phi| < 3\slash 5$,}\\
			\frac{30 }{\partial_n h(q_0)} \frac{a}{5\phi -1} = 15\frac{a}{\partial_n h(q_0)} &\text{if $\phi = 3\slash 5$,}\\
			\sgn  \frac{a}{\partial_n h(q_0)} \cdot \infty &\text{if $3\slash 5 < \phi < \sqrt{3\slash 5}$.}
		\end{cases}
	\end{equation*}
\end{proof}

\subsection{Discontinuity along the vertical line} \label{S:discontinuity in general}
For a nonconstant harmonic function $h\in\mathcal F$ and $f\in D_{\nu_h}(\partial^\perp_{V_0})$, with $L^2$-Hodge decomposition $f = \star_{\partial h}\partial u + \Phi$, we want to study the convergence of the ratio
\begin{equation*}
	Rf(F_wq_i) := \frac{\vec n \cdot (f\partial h)(F_wq_i)}{\vec n \cdot \partial h(F_wq_i)} = \frac{\partial_n u (F_wq_i)}{\partial_n h(F_wq_i)} + \sum_{w'} \Phi_w \frac{\partial_n \psi_{w'}(F_wq_i)}{\partial_n h(F_wq_i)},
\end{equation*}
when approaching $q_0$ along the vertical line $z_{n+1} :=  F_{0^n1}q_2$. We show that this problem is more delicate than the convergence along the sides.

\begin{proposition}\label{P:discontinuity along middle line}
	Let $h\in\mathcal F$ be nonconstant harmonic such that $\partial_T h(q_0) \neq 0$ and $f\in D_{\nu_h}(\partial^\perp_{V_0})$ with $L^2$-Hodge decomposition \eqref{eq L2 hodge} given by  $f = \star_{\partial h} \partial u + \Phi$.
	\begin{enumerate}
		\item If $u\in \mathcal F$ is harmonic, then
		\begin{equation*}
			\lim_{m\to\infty} \frac{\vec n \cdot \partial u(F_{0^n1}q_2)}{\vec n \cdot \partial h(F_{0^n1}q_2)} = \frac{\partial_T u(q_0)}{\partial_Th(q_0)}.
		\end{equation*}
		\item If the coefficients of $\Phi$ satisfy \eqref{E:only loops close to q0} and \eqref{E:loops close to q0}, then
		\begin{equation*}
			\lim_{n\to\infty} \frac{\vec n\cdot(\Phi\partial h)(F_{0^n1}q_2)}{\vec n \cdot \partial h(F_{0^n1}q_2)} = \begin{cases}
				\text{divergent} &\text{if $-\sqrt{3\slash 5} <\phi \leq -1\slash 5$,} \\
				\frac{10 a}{\partial_T^{\leftarrow} h(q_0)} \frac{1}{1-5\phi} &\text{if $|\phi| < 1\slash 5$,}\\
				\sgn \frac{a}{\partial_T^{\leftarrow} h(q_0)} \cdot \infty &\text{if $\phi \geq 1\slash 5$.}
			\end{cases}
		\end{equation*}
		\item If $u\in\mathcal F$ is harmonic and the tangential part $\vec t\cdot(\Phi\partial h)$ exists (see \Cref{T:tangential part}), then
		\begin{equation*}
			\lim_{m\to\infty} \frac{\vec n \cdot (f\partial h)(F_{0^n1}q_2)}{\vec n \cdot \partial h(F_{0^n1}q_2)} = \frac{\vec t \cdot (f\partial h)(q_0)}{\vec t \cdot \partial h(q_0)}.
		\end{equation*}
	\end{enumerate}
\end{proposition}
From \Cref{P:discontinuity along middle line} $(i)$ we see that continuity in $q_0$ of $f = \star_{\partial h} \partial u$ for $u\in\mathcal F$ harmonic can only be given if 
\begin{equation*}
	\partial_T u(q_0) = \frac{\partial_n u(q_0)}{\partial_n h(q_0) } \partial_T h(q_0).
\end{equation*}
It is a simple task to construct $h$ and $u$ such that this identity fails. Thus for $u\in \mathcal F$ harmonic, $\star_{\partial h}\partial u$ can not be expected to be continuous. On the contrary, one can solve for harmonic $u\in\mathcal F$ such that above equality holds for each $q_i \in V_0$. Assuming $\partial_T h(q_i) \neq 0$ for all $i$, the space of such solutions $u$ is two dimensional with
\begin{align*}
	u(q_0) =-\left[ \frac{\partial_T h(q_1)}{\partial_T h(q_0)} u(q_1) + \frac{\partial_T h(q_2)}{\partial_T h(q_0)} u(q_2) \right].
\end{align*}
However, it is unclear if for this choice of $u$ the ratio $R(\star_{\partial h}\partial u)$ is continuous on $V_\ast$.

For infinite loops $\Phi\in\ker\partial^\perp$ the statement in \Cref{P:discontinuity along middle line} $(ii)$, in some sense, extends the previous example in \Cref{P:discontinuity along sides using loops}. It shows that for infinite loops satisfying \eqref{E:only loops close to q0} and \eqref{E:loops close to q0} one can not expect convergence along the vertical line. Note that case $(iii)$ contains the case of finitely many loops.

From \Cref{P:discontinuity along middle line} $(iii)$ we see that, to expect continuity of the ratio $Rf$, one needs some relation between the ratio of normal parts and the ratio of tangential parts. Similar to the comment on the first case $(i)$, one may construct $u$ and $\Phi_w$ such that these two ratios actually agree in $q_0$. Again, continuity of the corresponding function $f$ is unclear.

\begin{remark} \
	\begin{enumerate}
		\item For harmonic functions $h$ the condition $\partial_T h(q_0) \neq 0$ is equivalent to $z_{n+1}\notin \ker\partial_n h$ for some $n\in\mathbb N_0$. This is immediately implied by \Cref{L:normal derivate at z}.
		\item A result similar to \Cref{P:discontinuity along middle line} $(i)$ for more general $u\in D_{\nu_h}(\Delta_{V_0})$ seems difficult due to the tangential derivative. However, such a generalization would be interesting due to the $L^2$-Hodge decomposition \eqref{eq L2 hodge}, which, for $f\in D_{\nu_h}(\partial^\perp_{V_0})$ with $L^2$-Hodge decomposition $f = \star_{\partial h}\partial u + \Phi$ implies $u\in D_{\nu_h}(\Delta_{V_0})$.
	\end{enumerate}
\end{remark}
For the proof of \eqref{P:discontinuity along middle line} we need the following lemma that calculates the values of the normal derivative of a harmonic function along the vertical line.
\begin{lemma}\label{L:normal derivate at z}
	Given $|w| = n$ and consider the cell $K_w$ with $V_{n+1}\cap K_w$ as in \Cref{fig:abcgasket}. If $u : K \to \mathbb R$ is harmonic in $K_w$, then
	\begin{equation*}
		\partial_n^{\rightarrow} u(z) = -3^{-|w|}\partial_T^{\leftarrow} u(a).
	\end{equation*}
\end{lemma}
\begin{proof}
	First, assume $n=0$ and $K_w = K$. Then, using the harmonicity of $u$ we can simply calculate
	\begin{equation*}
		\partial_n^{\rightarrow}u(F_1q_2) = r^{-1}(2u(F_1q_2) - u(F_1q_0) - u(q_1)) = -(u(q_1) - u(q_2)) = -\partial_T^{\leftarrow} u(q_0).
	\end{equation*}
	Now, for arbitrary $n$ and $w$ we use the scaling of $\partial_n$ and $\partial_T$ to get
	\begin{equation*}
		\partial_n^{\rightarrow}u(z) = r^{-|w|}\partial_n (u\circ F_w)(F_1q_2) = -r^{-|w|} \partial_T^{\leftarrow} (u\circ F_w)(q_0) = -3^{-|w|} \partial_T^{\leftarrow} u(a).
	\end{equation*}
\end{proof}

\begin{proof}[Proof of \Cref{P:discontinuity along middle line}]
	The first statemenet $(i)$ follows immediately from \Cref{L:normal derivate at z} as $u$ and $h$ are harmonic and $\vec n\cdot \partial u = \partial_n u$. 
	
	For $(ii)$, assume that $\Phi$ satisfies \eqref{E:only loops close to q0} and \eqref{E:loops close to q0}.	As $\supp \psi_w = K_w$, we have
	\begin{equation*}
		R\Phi(z_{n+1}) = \sum_{k=0}^n \Phi_{0^k} \frac{\partial_n\psi_{0^k}(z_{n+1}) }{\partial_n h(z_{n+1}) }.
	\end{equation*}
	To calculate $\partial_n\psi_{0^k}(z_{n+1})$ we use that $\sum_{\partial K_w} \partial_n u = 0$ if $u$ is harmonic in $K_w$. Thus,
	\begin{equation*}
		\partial_n^{\rightarrow}\psi_{0^k}(z_{n+1}) = -(\partial_n^{\nearrow}\psi_{0^k}(F_{0^{n}}q_0) + \partial_n^{\swarrow}\psi_{0^k}(F_{0^n}q_1)) = \partial_n^{\swarrow}\psi_{0^k}(F_{0^{n+1}}q_1) - \partial_n^{\swarrow}\psi_{0^k}(F_{0^n}q_1).
	\end{equation*}
	By \eqref{E:normal derivative psi} we have $\partial_n^{\rightarrow}\psi_{0^k}(z_{n+1}) = -10 \cdot 5^k 3^{-n}$. For $\partial_n h(z_{n+1})$ we use \Cref{L:normal derivate at z}. Thus,
	\begin{equation*}
		R\Phi(z_{n+1}) = \sum_{k=0}^n\Phi_{0^k}\frac{-10 \cdot 5^k 3^{-n}}{-3^{-n}\partial_T^{\leftarrow}h(q_0)} = \frac{10}{\partial_T^{\leftarrow} h(q_0)} \sum_{k=0}^n \Phi_{0^k}5^k = \frac{10 a}{\partial_T^{\leftarrow} h(q_0)} \sum_{k=0}^n (5\phi)^k.
	\end{equation*}
	Hence,
	\begin{equation*}
		\lim_{n\to\infty} R\Phi(z_{n+1}) = \begin{cases}
			\text{divergent} &\text{if $-\sqrt{3\slash 5} <\phi \leq -1\slash 5$,} \\
			\frac{10 a}{\partial_T^{\leftarrow} h(q_0)} \frac{1}{1-5\phi} &\text{if $|\phi| < 1\slash 5$,}\\
			\sgn \frac{a}{\partial_T^{\leftarrow} h(q_0)} \cdot \infty &\text{if $\phi \geq 1\slash 5$.}
		\end{cases}
	\end{equation*}
	
	At last we show $(iii)$. Assume that $u$ is harmonic and the tangential part of $\Phi\partial h$ exists in $q_0$. Due to linearity, we have
	\begin{equation*}
		\frac{\vec n \cdot (f\partial h)}{\vec n \cdot \partial h} = \frac{\vec n \cdot \partial u}{\vec n \cdot \partial h} + \frac{\vec n \cdot (\Phi\partial h)}{\vec n \cdot \partial h}.
	\end{equation*}
	For the ratio involving $u$ we use $(i)$. For the other ratio, as in the case $(ii)$ with \Cref{R:discontinuity} in mind, we use \Cref{L:normal derivate at z} to get
	\begin{equation*}
		\vec n \cdot (\Phi\partial h)(z_{n+1}) = \sum_{k=0}^n \Phi_{0^k} \partial_n \psi_{0^k}(z_{n+1}) = \pm 3^{-n}\sum_{k=0}^n \Phi_{0^k} \partial_T\psi_{0^k}(q_0). 
	\end{equation*}
	Then,  \Cref{T:tangential part} implies $\vec n \cdot (\Phi\partial h)(z_{n+1}) = \pm 3^{-n}\vec t \cdot (\Phi\partial h)(q_0)$. Note, the sign depends on the orientation we choose to calculate the normal and tangential derivative. However, using \Cref{L:normal derivate at z} with the same orientation of normal and tangential derivatives, we deduce
	\begin{equation*}
		\frac{\vec n \cdot (\Phi\partial h)(z_{n+1}) }{\vec n \cdot \partial h(z_{n+1})} = \frac{ \pm 3^{-n}\vec t \cdot (\Phi\partial h)(q_0)}{\pm 3^{-n}\partial_T h(q_0)} =  \frac{ \vec t \cdot (\Phi\partial h)(q_0)}{\partial_T h(q_0)}.
	\end{equation*}
	From here we easily conclude $(iii)$.
\end{proof}


\end{document}